\documentclass[11pt,
]{article}

\usepackage{amssymb}
\usepackage{amsfonts,amstext,amsmath,amsthm,
latexsym,mathrsfs,amsbsy,mparhack}
\usepackage[bbgreekl]{mathbbol}

\usepackage{graphicx}
\usepackage{color}

\usepackage{marginnote}

\usepackage[russian,english,greek]{babel}
\languageattribute{greek}{polutoniko}

\usepackage{makeidx}
\makeindex

\usepackage
[
citecolor=blue
linkcolor=yellow
]
{hyperref}

\input{s.sty}

\begin{document}

\selectlanguage{english}

\title{Models of set theory in which separation
theorem fails}

\author{Vladimir Kanovei\thanks
{Contact author. IITP RAS, 
Bolshoy Karetny per.\ 19, build.1, Moscow 127051, Russia. 
{\tt kanovei@googlemail.com}. 
Partial support of RFBR grant 17-01-00705 acknowledged.}
\and
Vassily Lyubetsky\thanks
{IITP RAS, 
Bolshoy Karetny per.\ 19, build.1, Moscow 127051, Russia. 
{\tt kanovei@googlemail.com}. 
Partial support of RFBR grant 18-29-13037 acknowledged.}
}

\maketitle

\begin{abstract}
We make use of a finite support product of 
the Jensen minimal forcing to 
define a model of set theory 
in which the separation theorem fails 
for projective classes $\fs1n$ and $\fp1n$, for 
a given $n\ge3$.



\vspace*{-6ex}

{\scriptsize%
\def\contentsname{}%
\tableofcontents
}

\end{abstract}

\np


\parf{Introduction}
\las{int}

The separation problem was introduced in 
descriptive set theory by Luzin~\cite{lbook}. 
In particular, Luzin asked whether 
(in modern notation for projective classes):
\ben
\Renu
\itlb{int1}%
any pair of disjoint $\fs1n$ sets of reals can be 
separated by a $\fd1n$ set,  

\itlb{int2}%
the remainders of two $\fs1n$ sets, 
obtained by the removal of their intersection, 
can be separated by disjoint $\fp1n$ sets, and 

\itlb{int3}%
there are two disjoint $\fp1n$ sets not separable 
\imar{snos lusepa}
by a $\fd1n$ set.
\een
Luzin underlined the importance and difficulty of 
these problems.\snos
{\label{lusepa}\mtho%
L'un des probl\'emes les plus importants de la th\'eorie
des ensembles projectifs et qui attend encore sa solution, 
est celui de leur \rit{s\'eparabilit\'e}.
On sait que deux ensembles analytiques quelconques sans point
commun sont toujours s\'eparables $\text{B}$. 
Il serait tr\'es important de d\'emontrer
que deux ensembles ($\text{A}_n$) quelconques sans point commun sont
s\'eparables ($\text{B}_n$).
De m\^eme, nous savons que si l'on supprime la partie commune 
\`a deux ensembles analytiques, les parties restantes sont 
s\'eparables au moyen de deux compl\'ementaires analytiques. 
La question se pose
naturellement de savoir si ce principe subsiste quand on 
remplace les  ensembles analytiques par ($\text{A}_n$) et les 
compl\'ementaires analytiques
par ($\text{CA}_n$). 
C'est un probl\'eme qui m\'erite d'attirer l'attention des
analystes malgr\'e sa difficult\'e.
D'ailleurs, il importe de savoir s'il existe deux ensembles 
($\text{CA}_n$) qui ne soient pas s\'eparables ($\text{B}_n$).
(Lusin~\cite[p.~289]{lbook}.) 
} 
Novikov characterized the separation problem as one  
of three main problems of descriptive set theory in 
\cite{nov1951e}, along with the measurability problem 
for $\fs12$ sets and the cardinality problem 
for $\fp11$ sets. 
(See \eg\ \cite{kl1e} on the two latter poblems.)

The problem is well known in descriptive set theory. 
In modern terms 
(see Moschovakis~\cite{mDST}, Kechris~\cite{dst}), 
the \rit{(first) separation theorem} for a class $\Ga$ 
of pointsets (sets in Polish spaces) is the claim that 
any two disjoint sets in $\Ga$ (in the same space) 
can be separated by a set in $\Ga\cap\dop\Ga$, where 
$\dop\Ga$ is the class of complements of \dd\Ga sets.  
The \rit{second separation theorem} for $\Ga$ 
claims that if $X,Y$ are sets in $\Ga$ 
(in the same space) 
then the sets $X'=X\bez Y$ and $Y'=Y\bez X$ 
are separable by two disjoint sets in $\dop\Ga$. 
Thus the content of 
\ref{int1}, \ref{int2}, \ref{int3} 
is this:\vom 

--- \ does the (first) separation theorem hold for $\fs1n$?\vom 

--- \
does the second separation theorem hold for $\fs1n$?\vom 

--- \
does the (first) separation theorem fail for $\fp1n$?\vom

\noi
Both separation theorems hold for $\fs11$ by 
Luzin \cite{lus:ea,lbook}, but fail for $\fp11$ 
by Novikov~\cite{nov1931}, and these results were 
known before the publication of the 
(French original) of \cite{lbook} in 1930.
Somewhat later, it was established 
by Novikov~\cite{nov1935} that the picture changes 
at the second projective level: 
both separation theorems hold for $\fp12$ 
but fail for $\fs12$. 

In the same time Kuratowski~\cite{kursep} proved 
the \rit{reduction theorem} for $\fs12$, that is, 
if $X,Y$ are sets in $\fs12$ then there exist disjoint 
sets $X'\sq X$ and $Y'\sq Y$ in the same class $\fs12$, 
with the same union $X'\cup Y'=X\cup Y$. 
Kuratowski also observed that Luzin's arguments in the 
proof of the separation theorem for $\fs11$ yield 
reduction for $\fp11$. 
Generally, 
if reduction holds for a projective class $\Ga$ then 
both separation theorems hold for the dual class 
$\dop\Ga$.

Generally, by classical studies, reduction 
holds for projective classes $\fp11$, $\fs12$ 
and fails for $\fs11$, $\fp12$, while the 
separation theorems hold for $\fs11$, $\fp12$ 
and fail for $\fp11$, $\fs12$. 
Note the inversion between the 1st and 2nd levels 
of the hierarchy.

As for the higher levels of projective hierarchy, 
all attempts made in classical descriptive set theory 
to solve the separation/reduction problems above the 
2nd level did not work, 
until some additional set theoretic axioms were 
commended. 
In particular, by Novikov~\cite{nov1951e} 
(see also Addison~\cite{add1}), 
G\"odel's \rit{axiom of constructibility} 
$\rV=\rL$ implies that, for any $n\ge 3$, 
reduction holds for $\fs1n$ 
and fails for $\fp1n$, while the 
separation theorems hold for $\fp1n$ 
and fail for $\fs1n$ --- pretty similar to 2nd level. 
On the contrary, 
by Addison and Moschovakis \cite{addmos}
and Martin~\cite{martAD}, 
the \rit{axiom of projective determinacy} $\PD$  
implies that, for any $m\ge 1$, 
reduction 
holds for projective classes $\fp1{2m+1}$, $\fs1{2m+2}$ 
and fails for $\fs1{2m+1}$, $\fp1{2m+2}$, while the 
separation theorems hold for $\fs1{2m+1}$, $\fp1{2m+2}$ 
and fail for $\fp1{2m+1}$, $\fs1{2m+2}$ --- 
pretty similar to what happens at the 1st and 2nd level 
corresponding to $n=0$ in this scheme. 
Moreover, by Steel ~\cite{steel_detsep}, 
it is true under the full 
\rit{axiom of determinacy} $\AD$, 
that if $\Ga$ is a class of pointsets closed under 
some simple operations and not self-dual 
(that is, $\Ga\ne\dop\Ga$), 
then reduction holds for exactly one of the classes 
$\Ga,\dop\Ga$, and the separation theorems hold 
for the other one. 
Conversely, 
Steel~\cite{steel_core} proved that a more special 
form of $\fp13$ separation implies otherwise 
impossible connections between some determinacy 
hypotheses.
See also Hauser and Schindler \cite{hashi} for other 
relevant results. 

These achievements still leave open important questions 
about the status of the separation theorems for higher 
projective classes. 
For instance the following: 

\bvo
[Mathias \cite{matsur} for $n=3$]
\lam{vo1}
Given a number $n\ge3$, is it consistent with $\ZFC$ 
that the (first) separation theorem fails for  
both $\fs1n$ and $\fp1n$? \qed
\evo

Harrington   
solved the problem in the positive by means 
of a generic extension of $\rL$ in which 
the (first) separation theorem fails for  
both $\fs13$ and $\fp13$. 
The solution was obtained by the technique of 
almost-disjoint forcing of \cite{jsad}, 
and was sketched in unpublished notes 
\cite[Part B]{h74}. 
The result itself was 
mentioned, with a reference to Harrington,
\eg\ in 
Moschovakis~\cite[5B.3]{mDST}. 
Harrington also suggested in \cite{h74}  
some substantial changes in the construction  
of the generic extension, intended  to get 
the failure of separation for  
both $\fs1n$ and $\fp1n$ for a given $n>3$, 
or even for all $n$, but 
such a generalization has never materialized in 
detail 
(albeit mentioned in \cite{hin,matsur,samiPHD}).

Our goal here is to prove the next theorem, 
which solves Problem~\ref{vo1} in the positive for any 
given $n>3$, albeit by a method different from the 
one used~in~\cite{h74}. 

\bte
\lam{Tsep}
Let\/ $\nn\ge3$. 
It is true, in a suitable generic extension of\/ $\rL$, 
that 
\ben
\renu
\itlb{Tsep1}%
there is a pair of 
disjoint lightface\/ $\ip1{\nn}$ sets\/ $X,Y\sq\dn$, not 
separable by disjoint\/ $\fs1{\nn}$ sets, and hence\/ 
separation fails for both\/ $\fp1{\nn}$ and\/ 
$\ip1{\nn}\,;$ 

\itlb{Tsep2}%
there is a pair of 
disjoint lightface\/ $\is1{\nn}$ sets\/ $X,Y\sq\dn$, not 
separable by disjoint\/ $\fp1{\nn}$ sets, and hence\/ 
separation fails for both\/ $\fs1{\nn}$ and\/ 
$\is1{\nn}$.
\een
\ete

\parf{Outline of the proof} 
\las{out}


Given $\nn\ge3$, our plan is to define a sequence of  
forcing notions $\dP_\xi$, $\xi<\omi$ in $\rL$, whose 
finite-support product $\dP=\prod_\xi\dP_\xi$ satisfies CCC 
and adjoins a sequence of generic reals 
$x_\xi\in\dn$, 
that are independent of each other in the sense 
that 
\ben
\Renu
\itlb{out1}
if 
$\et<\omi$, 
then (a) the submodel 
$\rL[\sis{x_\xi}{\xi\ne\et}]$ contains no reals 
$\dP_\et$-generic over $\rL$, and moreover, 
(b) $x_\et$ is 
the only real in $\rL[\sis{x_\xi}{\xi<\omi}]$, 
$\dP_\et$-generic over $\rL$,
\een
and have the following definability property:
\ben
\Renu
\atc
\itlb{out2}
the relation 
``$x\in\dn$ is a real $\dP_\xi$-generic over $\rL$'' 
(of arguments $x,\xi$)  
is $\ip1{\nn-1}$ in the whole extension and 
any its submodel. 
\een
Then we generically split $\omi$ into three 
unbounded sets $\omi=\Om_1\cup\Om_2\cup\Om_3$, 
define 
$\Da=\ens{2\nu}{\nu\in \Om_1\cup \Om_3} 
\cup\ens{2\nu+1}{\nu\in \Om_2\cup \Om_3}$, 
and prove that $\Om_1$ and $\Om_2$ are disjoint 
$\ip1{\nn}$ sets not separable by disjoint 
$\fs1{\nn}$ sets in the model 
$M=\rL[\sis{x_\xi}{\xi\in\Da}]$. 
Indeed by \ref{out1} we have 
$$
\Om_1=\ens{\nu<\omi}
{\neg\:\sus x\,
(x\text{ is $\dP_{2\nu+1}$-generic over $\rL$})}
$$
in $M$, so $\Om_1$ is $\ip1{\nn}$ in $M$ by \ref{out2}, 
and accordingly so is $\Om_2$.
The non-separability claim involves the following 
crucial property of \dd\dP generic extensions:
\ben
\Renu
\atc
\atc
\itlb{out3}
if $X\in\rL$, $X\sq\omi$ is unbounded in\/ 
$\omi$, and $G\sq\fP$ is 
\dd\fP generic over $\rL$  
then\/ $\rL[\sis{x_\xi}{\xi\in X}]$ is an elementary
submodel of 
$\rL[G]$ \poo\ all  $\is1{\nn-1}$ formulas. 
\een
Each factor forcing $\dP_\xi$ in this scheme is a clone 
of Jensen's minimal forcing, 
defined in \cite{jenmin}  
(\rit{Jensen's forcing} below, for the sake of brevity, 
see also \cite[28A]{jechmill} on this forcing). 
in particular, it consists of \rit{perfect trees} in $\bse.$
The idea to use finite-support products of Jensen's 
forcing in order to obtain models with different 
definability effects belongs to Enayat \cite{ena}. 
It was exploited to obtain generic models 
with: countable non-empty $\ip12$ sets 
(even $\Eo$-classes) 
with no $\od$ elements \cite{kl22,kl27e},  
a countable $\ip12$ Groszek -- Laver pair 
\cite{kl25}, 
planar $\ip12$ sets with countable 
cross-sections and OD-non-uniformizable \cite{kl28,kl30e}, 
and also a model where the separation 
theorem fails for both $\fs13$ and $\fp13$ \cite{kl28}.

The latter result corresponds to the case $\nn=3$ of 
Theorem~\ref{Tsep}, in which case \ref{out3} 
is immediately true by Shoenfield.
On the other hand, conditions similar to \ref{out1}, 
\ref{out2} for $\nn=3$, are involved in the   
forcing constructions in \cite{kl25,kl22,kl28,kl27e,kl30e}, 
and in \cite{jenmin} itself, where a CCC forcing 
$\dJ\in\rL$ is defined to add a real $a\in\dn$  so that  
$a$ is the only \dd\dJ generic real in $\rL[a]$, 
and ``being a \dd\dJ generic real'' is $\ip12$. 
These properties are implied by a special construction 
of $\dJ=\bigcup_{\al<\omi}\dJ_\al$ in $\rL$ from 
countable sets $\dJ_\al$ of perfect trees. 
The construction can be viewed as a maximal branch in 
a certain mega-tree, say $\cP$, whose nodes are countable sets 
of perfect trees, and each $\dJ_\al$ is chosen just 
as the \dd\lel least appropriate extension.
The complexity of this construction is $\id12$ in the 
codes, that leads to the $\ip12$ definability of 
being generic, while a suitable character of extension 
in the mega-tree allows to ``kill'' all possible 
competitors of $a$ to be $\dJ$-generic. 

Pretty similar ideas and constructions work in the 
mentioned papers, in particular in \cite{kl28}, where 
a model is defined in which $\fp13$-separation fails.  

A method of reproducing some generic counterexamples, 
originally defined on 2nd and 3rd projective level,
at any given higher projective level $\nn$, 
was introduced 
by Harrington~\cite{h74} on the base of the 
almost-disjoint forcing \cite{jsad}, 
and independently in \cite{k78} on the base of Jensen's 
forcing of \cite{jenmin}. 
In the terms above, the method requires to define 
a maximal branch in $\cP$ that 
intersects all dense sets in $\cP$ of descriptive 
complexity $\nn$ 
(or $\nn+c$, where $c$ is a small entire constant 
depending on the nature of the problem).  
The method was recently applied to get models 
in which, for a given $\nn\ge2$, there exists: 
\ben
\aenu
\itlb{1out}
a $\ip1\nn$ \dd\Eo equivalence class containing no 
$\od$ elements, while 
every countable \dd{\is1{\nn}}set of reals 
contains only 
$\od$ reals \cite{kl34}, 

\itlb{2out}
a $\varPi^1_\nn$ singleton $\ans a$ such that 
$a$ codes a cofinal map 
$f:\omega\to\omega_1^\rL$ minimal over $\rL$,
while every $\varSigma^1_\nn$
set $X\sq\om$ is constructible \cite{kl36},

\itlb{3out}
a planar non-ROD-uniformizable $\ip1{\nn}$ set, whose 
all vertical cross-sections are $\Eo$-classes, 
while all planar $\fs1\nn$ sets with countable  
cross-sections are $\fd1{\nn+1}$-uniformizable \cite{kl38}.
\een
Here the method is used to prove Theorem~\ref{Tsep}.

Sections \ref{tre} to \ref{jex}: perfect trees 
in $\bse,$ perfect tree forcing notions, 
multitrees (finite products of trees), 
multiforcings (countable products of forcings), 
splitting, 
%
the refinement relation, 
generic refinements by Jensen's splitting construction.

Sections \ref{pres} to \ref{comb}: 
properties of generic refinements, sealing dense sets, 
sealing real names, and applications to generic extensions. 

Sections \ref{incS} to \ref{bpf}: 
we define the set $\vmf$ of all countable sequences $\vjpi$ 
of small \muf s, increasing in the sense of the refinement 
relation. 
Arguing in $\rL$, we define 
a $\id1{\nn-1}$ (in the codes) 
maximal branch $\vjPi$ in $\vmf$, 
which \se s all $\is1{\nn-2}$ sets in $\vmf$, 
where $\nn$ is the number in Theorem~\ref{Tsep}, and 
$\vjpi\in\vmf$ \rit{\se s} a set $W\sq\vmf$ if either 
$\vjpi\in W$ or no extension of $\vjpi$ in $\vmf$ 
belongs to $W.$ 
The forcing notion $\fP$ for Theorem~\ref{Tsep} is 
a derivate of $\vjPi$.

Sections \ref{kge} to \ref{nspf}: 
we show that $\fP$ satisfies \ref{out1} and \ref{out2}.

Sections \ref{auxA} to \ref{ee}: 
to achieve \ref{out3}, we develop an auxiliary forcing 
notion $\fo$, that approximates the truth in 
$\fP$-generic extensions for 
$\is1{\nn-1}$-formulas and below, so that the relation 
$\fo$ restricted to any class $\is1m$ or $\ip1m$, 
$m\ge2$, is $\is1m$, resp., $\ip1m$. 
Using the invariance of $\fo$ under certain 
transformations (while $\fP$ is not invariant!), 
we accomplish the proof of \ref{out3} and Theorem~\ref{Tsep}.

\parf{Trees and perfect-tree forcing notions} 
\las{tre}

Let $\bse$ be the set of all strings (finite sequences) 
\index{string}%
of numbers $0,1$.
If $t\in\bse$ and $i=0,1$ then 
$t\we i$ is the extension of $t$ by $i$. 
If $s,t\in\bse$ then $s\sq t$ means that $t$ extends $s$, while 
$s\su t$ means proper extension. 
$\lh t$ is the length of $t$,  
\index{zzlhs@$\lh t$, the length}%
and $2^n=\ens{t\in\bse}{\lh t=n}$ (strings of length $n$).%

A set $T\sq\bse$ is a \rit{tree} iff 
\index{tree}%
for any strings $s\su t$ in $\bse$, if $t\in T$ then $s\in T$.
Thus every non-empty tree $T\sq\bse$ contains the 
\rit{empty string} $\La$. 
\index{zzLa@$\La$, the empty string}%

If $T\sq\bse$ is a tree and $s\in T$ then put 
$T\ret s=\ens{t\in T}{s\sq t\lor t\sq s}$. 
\index{zzT/s@$T\ret s$}%

\bdf
\lam{body}
$\pet$ is the set of all \rit{perfect} trees 
\index{zzPT@$\pet$, perfect trees}%
$\pu\ne T\sq \bse$. 
\imar{pet}%
Thus a tree $\pu\ne T\sq\bse$ belongs to $\pet$ iff 
it has no endpoints and no isolated branches. 
%
If $T\in\pet$ then define a perfect set  
$$
[T]=\ens{a\in\dn}{\kaz n\,(a\res n\in T)}\sq\dn.  
$$  
\index{zzT]@$[T]$}%
Trees $T,S\in\pet$ are \rit{almost disjoint}, 
\rit{\ad} for brevity, iff 
\index{tree!almost disjoint}%
\index{tree!almost disjoint@\ad\ trees}%
\index{almost disjoint, \ad!\ad\ (trees)}%
the intersection $S\cap T$ is finite; 
this is equivalent to just $[S]\cap[T]=\pu$.
A set $\dA\sq\pet$ is an \rit{antichain} iff any 
two trees 
\index{antichain}
$T\ne T'$ in $\dA$ are \ad.
\edf

We'll consider pairs of the form 
$\ang{n,T}$, where $n<\om$ and $T\in\pet$. 
Following \cite{abr}, the set $\om\ti\pet$ of 
such pairs is ordered by a special relation $\lec$ 
so that $\ang{n,T}\lec \ang{m,S}$ 
(reads: $\ang{n,T}$ \rit{extends} $\ang{m,S}$) 
iff $m\le n$, $T\sq S$, and 
$T\cap2^m=S\cap2^m$.\snos
{This definition does not explicitly contain any 
splitting condition. 
This is why one needs the genericity condition in 
Lemma~\ref{fus}. 
An earlier definition in \cite{bala}
stipulates that 
for any $s\in S\cap2^m$ 
there exist two strings $s'\ne s''$ in $T\cap2^{n}$ 
such that $s\su s'$ and $s\su s''$. 
With such an ordering, Lemma~\ref{fus} holds without 
the genericity condition.}
%
The role of the number $m$ in a pair $\ang{m,S}$
is to preserve  
the value $S\cap2^m$ under 
$\lec$-extensions. 

The implication $m>n\imp\ang{m,T}\lec \ang{n,T}$
(the same $T$!) always holds, but 
$S\sq T\imp\ang{n,S}\lec \ang{n,T}$ 
is not necessarily true:  
we also need $T\cap2^n=S\cap2^n$. 


\ble
[see \cite{abr}]
\lam{fus}
Let\/ 
$\ldots \lec\ang{n_2,T_2}\lec\ang{n_1,T_1}\lec\ang{n_0,T_0}$ 
be a decreasing sequence in\/ $\om\ti\pet$, 
with\/ $n_0<n_1<n_2<\ldots$ strictly, minimally generic 
in the sense that it meets every set of the form
$$
D_t=\ens{\ang{n,T}\in\om\ti\pet}
{t\nin T\lor \sus s\in T\,
(t\sq s\land s\we0,s\we1\in T)}\,,\;
t\in\bse.
$$ 
Then\/ $T=\bigcap_nT_n\in\pet$, and 
if\/ $i<\om$ then\/ 
$\ang{n_i,T}\lec\ang{n_i,T_i}$.\qed
\ele

\bdf
\lam{ptfn}
Let an \rit{arboreal forcing} 
be any set 
$\dP\sq\pet$ such that if $u\in T\in\dP$ 
then $T\ret u\in \dP$.
Let $\ptf$ be the set of all such sets $\dP$.
A forcing $\dP\in\ptf$ is:  
\bde
\item[\rit{regular}\rm,] 
if for any $S,T\in\dP$, the intersection  
\index{arboreal forcing, $\ptf$!regular}%
\index{regular}%
$[S]\cap [T]$ is clopen in $[S]$ or clopen in $[T]$ 
(or clopen in both $[S]$ and $[T]$);

\item[\rit{special}\rm,]  
if there is 
\index{arboreal forcing, $\ptf$!special}%
\index{special}%
a finite or countable antichain $\dA\sq\dP$ 
such that $\dP=\ens{T\ret s}{s\in T\in \dA}$ --- 
the antichain $\dA$ 
is unique in this case, and the forcing $\dP$ itself
is obviously regular.
\qed
\ede
\eDf

\bpri
\lam{cloL}
If $s\in\bse$ then the 
tree $T[s]=\ens{t\in\bse}{s\sq t\lor t\sq s}$ 
\index{zzTs]@$T[s]$}%
belongs to $\pet$ and  
$T[s]=\req{(\bse)}s\,,\:\kaz s$. 
The set $\dpo=\ens{T[s]}{s\in\bse}$ 
\index{arboreal forcing!Cohen forcing, $\dpo$}%
\index{Cohen forcing, $\dpo$}%
\index{zzPcoh@$\dpo$}%
(the Cohen forcing) 
is a regular and special arboreal forcing notion.
\vyk{
The set $\dpop$ of all finite (non-empty) unions of 
trees in $\dpo$
\index{zzPcohp@$\dpop$}%
is a regular non-special arboreal forcing notion.
}
\epri

Any set $\dP\in\ptf$ can be considered as a forcing notion 
(if $T\sq T'$ then $T$ is a stronger condition); 
such a forcing $\dP$ obviously adds a real in $\dn$. 
Lemma~\ref{regfn} below implies that the compatibility in 
\rit{regular} forcing notions is absolute.

To carry out splitting constructions, 
as in Lemma~\ref{fus}, 
over a forcing $\dP\in\ptf$, we make use of  
a bigger forcing notion $\clg\dP\in\ptf$, 
that consists of all finite 
unions of trees in $\dP$. 
\index{zzUfinP@$\clg\dP$}%
\index{tree!finite unions}%
Then   
$\dP$ is dense in $\clg\dP$, so the forcing 
properties of both sets coincide. 
Yet $\clg\dP$ is more flexible \poo\ tree constructions.

\ble
\lam{regfn}
Assume that\/ $\dP\in\ptf$ 
is regular and\/ $S,T\in\dP$ are not \ad. 
Then\/
$S\cap T\in\clg\dP$, hence the trees\/ $S,T$
are compatible in\/ $\dP$.
\qed
\ele

\ble
\lam{nq}
Let\/ $\dP\in\ptf$ and\/ $S,T\in\clg\dP$, 
$u\in S$, $n=\lh u$, $T\sq S\ret s$. 
Then the tree\/ 
$S'=T\cup\bigcup_{v\in S\cap2^n,v\ne u}S\ret v$
belongs to\/ $\clg\dP$, $\ang{n,S'}\lec\ang{n,S}$, 
$S'\ret u=T$, and\/ 
$S'\ret v=S\ret v$ whenever\/ 
$v\in S$, $\lh v=n$, $v\ne u$.
\qed
\ele


\bcor
\lam{suz}
Assume that\/ $\dP,\dP'\in\ptf$. 
Then
\ben
\renu
\itlb{suz0}%
if\/ $n<\om$ and $T\in\clg\dP$, then 
there is a tree\/ $S\in\clg\dP$ such that\/ 
$\ang{n,S}\lec\ang{n,T}$ and\/ $S\ret t\in\dP$ 
{\rm(not just $\in\clg\dP$!)}
for all\/ $t\in2^n\cap S\,,$

\itlb{suz1}%
if\/ $T\in\dP$ and\/ $T'\in\dP'$, then 
there are trees\/ $S\in\dP$, $S'\in\dP'$ such that\/ 
$S\sq T$, $S'\sq T'$, and\/ $[S]\cap[S']=\pu\,;$.

\itlb{suz2}%
if\/ $n<\om$ and\/ $T\in\clg\dP$, $T'\in\clg{\dP'}$, then 
there exist trees\/ $S\in\clg\dP$, $S'\in\clg{\dP'}$ 
\st\/ 
$\ang{n,S}\lec\ang{n,T}$, 
$\ang{n,S'}\lec\ang{n,T'}$, 
and\/ $[S]\cap[S']=\pu$.
\een
\ecor
\bpf
\ref{suz1}
If $T=T'$ then pick a pair of strings 
$u\ne v$ in $T=T'$ with $\lh u=\lh v$, and let 
$S=T'\ret{u}$, $S'=T'\ret{v}$.
If say $T\not\sq T'$ then let $u\in T\bez T'$, 
$S=T\ret{u}$, and simply $S'=T'$.
To prove \ref{suz2} iterate \ref{suz1} and make use of 
Lemma~\ref{nq}.
\epf

\parf{\Muf s and \mut s}   
\las{mul}

Call a {\ubf\muf} 
\index{multiforcing, $\mfp$}%
\index{zzMF@$\mfp$}%
any map $\jpi:\abc\jpi\to\ptf$, where 
$\abc\jpi=\dom\jpi\sq\omi$. 
Let $\mfp$ be the collection of all \muf s. 
\kmar{mfp} 
Every $\jpi\in\mfp$ will be typically presented as 
an indexed set $\jpi=\sis{\dP_\xi}{\xi\in\abs\jpi}$, 
where $\dP_\xi\in\ptf$ for all $\xi\in\abs\jpi$,
so that each set $\dP_\xi=\djpi\xi=\jpi(\xi)$, 
\kmar{djpi xi}%
\index{zzPpixi@$\djpi\xi$}%
$\xi\in\abc\jpi$, 
is an arboreal forcing notion. 
Such a $\jpi$ is:
\bit
\item[$-$]
\rit{small}, if both $\abc\jpi$ and each 
\index{multiforcing!small}%
\index{small}%
forcing $\djpi\xi$, $\xi\in\abc\jpi$, are countable;

\item[$-$]
\rit{\mdi}, 
\index{multiforcing!special}%
\index{special}%
if each $\djpi\xi$ is special in the sense of 
Definition~\ref{ptfn};

\item[$-$] 
\rit{\mre}, if each $\djpi\xi$ is regular, 
\index{multiforcing!regular}%
\index{regular}%
in the sense of Definition~\ref{ptfn}.
\eit

Let a \rit{\mut} be any function  
\index{multitree, $\md$}%
$\zp:\abs\zp\to\pet$, 
with a finite 
\rit{support} 
$\abs \zp=\dom\zp$; 
$\md$ will be the collection of all \mut s.
\imar{md}%
\index{zzMT@$\md$, all \mut s}%
Every $\zp\in\md$ will be typically presented as 
an indexed set $\zp=\sis{\zc\zp\xi}{\xi\in\abs\zp}$, 
\imar{zc zp xi}%
\index{zzTpxi@$\zc\zp\xi$}%
where $\zc\zp\xi=\zp(\xi)\in\pet$ for all $\xi\in\abs\zp$.

Let 
$\jpi=\sis{\dP_\xi}{\xi\in\abs\jpi}$ 
be a \muf. 

A \rit{\dd\jpi\mut} is any $\zp\in\md$   
\index{multitree!\dd\jpi\mut}%
such that $\abs \zp\sq\abs\jpi$, 
and if $\xi\in\abs\zp$ then the tree  
$\zp(\xi)=\zc\zp\xi$ belongs to $\dP_\xi$. 
The set $\mt\jpi$  of all \dd\jpi\mut s
\index{zzMTpi@$\mt\jpi$, all $\jpi$-\mut s}%
is equal to the finite support product 
$\prod_{\xi\in\abs\jpi}\dP_\xi$, and if $\zp\in\mt\jpi$ 
then the set
$$
[\zp]=\ens{x\in(\dn)^{\abs\jpi}}
{\kaz \xi\in\abs\zp\,(x(\xi)\in [\zc\zp\xi])} 
$$
is a cofinite-dimensional perfect cube 
in $(\dn)^{\abs\jpi}$.
We order $\md$ and each $\mt\jpi$ componentwise: 
$\zq\leq\zp$ 
($\zq$ is stronger that $\zq$) 
iff $\abs\zp\sq\abs\zq$ and 
$\zc\zq\xi\sq\zc\zp\xi$ for all $\xi\in\abs\zp$; 
this is equivalent to $[\zq]\sq[\zp]$.
The empty \mut\ $\jLa$ defined by $\abs\jLa=\pu$, 
belongs to $\mt\jpi$ and 
is the weakest condition. 

\bdf
\lam{dsad}
Multitrees $\zp,\zq\in\mt\jpi$ are 
\rit{somewhere almost disjoint} (\sad) 
\index{somewhere almost disjoint, \sad}%
\index{sad@\sad}%
if there is $\xi\in\abs\zp\cap\abs\zq$ such that 
the trees $\zc\zp\xi$ and $\zc\zq\xi$ are \ad. 
Being \sad\ is equivalent to $[\zp]\cap[\zq]=\pu$, 
and, in the case of regular \muf s $\jpi$, 
equivalent to the \rit{incompatibility} in $\mt\jpi$ 
by the following result.
\edf

\bcor
[of Lemma~\ref{regfn}]
\lam{lsad}
Assume that\/ $\jpi$ is a regular \muf\ and\/ 
$\zp,\zq\in\mt\jpi$ 
are not\/ \sad. 
Then there is a finite set\/ 
$R\sq\mt\jpi$ such that\/ 
$[\zp]\cap[\zq]=\bigcup_{\zr\in R}[\zr]$.
Therefore $\zp,\zq$ are compatible in\/ $\mt\jpi$, 
that is, there is a \mut\/ $\zr\in\mt\jpi$ 
satisfying\/ 
$\zr\leq\zp$ and\/ $\zr\leq\zq$.\qed
\ecor

\bdf
\lam{kw}
The \rit{componentwise union} of \muf s 
$\jpi,\jqo$ is a \muf\ $\jpi\kw\jqo$ satisfying 
\kmar{kw}
$\abc{(\jpi\kw\jqo)}=\abc\jpi\cup\abc\jqo$ and 
\index{componentwise union!piUqu@$\jpi\kw\jqo$}%
\index{zzzpiUqu@$\jpi\kw\jqo$}
$$
(\jpi\kw\jqo)(\xi)\;=\;
\jpi(\xi)\;\text{ or }\; 
\jqo(\xi)\;\text{ or }\;
\jpi(\xi)\cup\jqo(\xi)
$$
in cases resp.\ $\xi\in\abc\jpi\bez\abc\jqo$, 
$\xi\in\abc\jqo\bez\abc\jpi$, 
$\xi\in\abc\jqo\cap\abc\jpi$.

If $\vjpi=\sis{\nor\jpi\al}{\al<\la}$ is a sequence 
in $\mfp$ then 
define 
$\jpi=\bkw\vjpi=\bkw_{\al<\la}\nor\jpi\al\in\mfp$ so that 
\kmar{bkw}%
\index{componentwise union!Ucwb@$\bkw\vjpi=
\bkw_{\al<\la}\nor\jpi\al$}%
\index{zzUcwb@$\bkw_{\al<\la}\nor\jpi\al$}%
$\abc\jpi=\bigcup_{\al<\la}\abc{\nor\jpi\al}$ and 
$\jpi(\xi)=
\bigcup_{\al<\la,\:\xi\in\abc{\nor\jpi\al}}
{\nor\jpi\al}(\xi)$ for  
$\xi\in\abc\jpi$.
\edf

\bre
\lam{adds}
Any forcing of the form $\mt\jpi$, 
where $\jpi=\sis{\dP_\xi}{\xi\in\abs\jpi}\in\mfp$, 
adds a generic sequence    
$\sis{x_{\xi}}{\xi\in\abs\jpi}$, where each 
$x_{\xi}=x_\xi[G]\in\dn$ 
is a \dd{\dP_\xi}generic real. 
Reals of the form $x_{\xi}[G]$ will be called 
\rit{principal generic reals\/} in $\rV[G]$. 
\index{reals!principal generic reals, $x_{\xi}[G]$}%
\index{principal generic reals, $x_{\xi}[G]$}%
\index{zzxxiG@$x_{\xi}[G]$}%
\ere


\parf{Refining arboreal forcings} 
\las{em}  

If $T\in\pet$ (a tree) and $D\sq\pet$ then 
$X\sqf\bigcup D$ will mean 
\index{zzsqf@$\sqf$}%
that there is a finite set $D'\sq D$ such that 
$T\sq\bigcup D'$, or equivalently $[T]\sq\bigcup_{S\in D'}[S]$.

\bdf
\lam{fm}
Let $\dP,\dQ\in\ptf$ be arboreal forcing 
notions. 
Say that $\dQ$ is a \rit{refinement} of $\dP$ 
\index{arboreal forcing, $\ptf$!refinement, $\dP\ssq\dQ$}
\index{forcing!refinement, $\dP\ssq\dQ$}
\index{zzPssqQ@$\dP\ssq\dQ$}%
(symbolically $\dP\ssq\dQ$) if 
\ben
\nenu
\itlb{fm1}%
the set 
$\dQ$ is dense\snos
{\label{fnden}%
If $\dP\sq \dR\sq\pet$ then, as usual, $\dP$ is 
1) 
\rit{dense} in $\dR$ iff 
\index{set!dense}%
\index{dense}%
$\kaz T\in \dR\,\sus S\in \dP\,(S\sq T)$, 
2) 
\rit{open dense} in $\dR$ iff in addition 
\index{set!open dense}%
\index{open dense}%
$\kaz T\in \dR\,\kaz S\in \dP\,(T\sq S\imp T\in\dP)$, 
and 3) 
\rit{pre-dense} in $\dR$ iff the derived set  
\index{set!pre-dense}%
\index{pre-dense}%
$\dP'=\ens{T\in\dR}{\sus S\in\dP(T\sq S)}$ 
is dense in $\dR$.
} 
in $\dP\cup\dQ$: 
if $T\in\dP$ then $\sus Q\in\dQ\,(Q\sq T)$;


\itlb{fm3}%
if $Q\in\dQ$ 
then $Q\sqf\bigcup\dP$;

\itlb{fm4}%
if $Q\in\dQ$ and $T\in\dP$ then $[Q]\cap[T]$ is 
clopen in $[Q]$ and $T\not\sq Q$.\qed 
\een
\eDf

\ble
\label{pqr}
\ben
\renu
\itlb{pqr0}%
If\/ $\dP\ssq\dQ$ and\/ $S\in\dP$, $T\in\dQ$, then\/ 
$[S]\cap[T]$ is meager in\/ $[S]$, 
therefore\/ $\dP\cap\dQ=\pu$ and\/ $\dQ$ is 
open dense in\/ $\dP\cup\dQ\;;$

\itlb{pqr1}%
if\/ $\dP\ssq\dQ\ssq\dR$ then\/ $\dP\ssq\dR$,
thus\/ $\ssq$ is a strict partial order$;$  

\vyk{
\itlb{pqr4}%
if\/ $\dQ\in\ptf$ is special, $R\in\pet$, 
$R\sq \baQ\in\dQ$, then there is a special forcing\/ 
$\dR\in\ptf$ containing\/ $R$ and such that\/ 
$\dP\ssq\dQ$ implies\/ 
$\dP\ssq\dR\;;$
}

\itlb{pqr2}%
if\/ $\sis{\dP_\al}{\al<\la}$  
is a\/ \dd\ssq increasing sequence in\/ \ptf\ and\/  
$0<\mu<\la$ then\/
$\dP=\bigcup_{\al<\mu}\dP_\al\ssq
\dQ=\bigcup_{\mu\le\al<\la}\dP_\al\;;$

\itlb{pqr3}%
if\/ $\sis{\dP_\al}{\al<\la}$  
is a\/ \dd\ssq increasing sequence in\/ \ptf\ and each\/  
$\dP_\al$ is special then\/ 
$\dP=\bigcup_{\al<\la}\dP_\al\in\ptf$, 
$\dP$ is regular, 
and each\/ $\dP_\ga$ is pre-dense in\/ 
$\dP$.
\een                                             
\ele
\bpf
\ref{pqr0}
Otherwise there is a string $u\in S$ such that 
$S\ret u\sq[T]\cap[S]$. 
But $S\ret u\in\dP$, which contradicts to 
\ref{fm}\ref{fm4}. 

\ref{pqr1}, \ref{pqr2}
Make use of \ref{pqr0} to establish \ref{fm}\ref{fm4}.

\ref{pqr3}
To check the regularity, let 
$S\in\dP_\al$, $T\in\dP_\ba$, $\al\le\ba$. 
If $\al=\ba$ then, as $\dP_\al$ is special, the trees 
$S,T$ are either \evd\ or \dd\sq comparable by
Lemma~\ref{regfn}.
If $\al<\ba$ then $[S]\cap[T]$ is clopen in $[T]$
by \ref{fm}\ref{fm4}. 

To check the pre-density, 
let $S\in\dP_\al$, $\al\ne\ga$. 
If $\al<\ga$ then by \ref{fm}\ref{fm1} there is a tree 
$T\in\dP_\ga$, $T\sq S$. 
Now let $\ga<\al$. 
Then $S\sqf\bigcup\dP_\ga$ by \ref{fm}\ref{fm3}, hence  
there is a tree $T\in\dP_\ga$ such that $[S]\cap[T]\ne\pu$. 
However $[S]\cap[T]$ is clopen in $[S]$ by \ref{fm}\ref{fm4}. 
Therefore $S\ret u\sq T$ for a string $u\in S$.
Finally $S\ret u\in \dP_\al$ since $\dP_\al\in\ptf$.
\epf         

Note that if $\dP,\dQ\in\ptf$ and $\dP\ssq\dQ$ then a 
dense set $D\sq \dP$ is not necessarily dense or even 
pre-dense in $\dP\cup\dQ$. 
Yet there is a special type of refinement which preserves 
at least pre-density.

\bdf
\lam{ssqm}
Let $\dP,\dQ\in\ptf$ and 
$D\sq \dP$. 
Say that $\dQ$ \rit{seals\/ $D$ over\/ $\dP$}, 
\index{refinement!seals}%
\index{seals}%
\index{refinement1D@refinement, $\dP\ssa D\dQ$}
\index{arboreal forcing, $\ptf$!refinement1D@refinement, 
$\dP\ssa D\dQ$}%
\index{zzPssq1DQ@$\dP\ssa D\dQ$}%
symbolically $\dP\ssa D\dQ$, if $\dP\ssq\dQ$ holds and  
\imar{ssa D}
every tree $S\in\dQ$ satisfies $S\sqf \bigcup D$.
Then simply $\dP\ssq\dQ$ is equivalent to $\dP\ssa\dP\dQ$. 
\edf

As we'll see now, a sealed set has to be pre-dense both before 
and after the refinement. 
The additional importance of sealing refinements lies in fact that, 
once established, it preserves under further simple 
refinements, that is, $\ssa D$ is transitive in 
a combination with $\ssq$ in the sense of \ref{pqm1} 
of the following lemma:

\ble
\label{pqm}
\ben
\renu
\itlb{pqm0}%
If\/ $\dP\ssa D\dQ$ then\/ $D$ is pre-dense in\/ $\dP\cup\dQ$, 
and if in addition\/ $\dP$ is regular then\/ 
$D$ is pre-dense in\/ $\dP$ as well$;$
\imar{pqm}

\itlb{pqm1}%
if\/ $\dP\ssa D\dQ\ssq\dR$ 
{\rm(note: the second $\ssq$ is not $\ssa D$!)} 
then\/ $\dP\ssa D\dR\;;$ 

\itlb{pqm2}%
if\/ $\sis{\dP_\al}{\al<\la}$  
is a\/ \dd\ssq increasing sequence in\/ \ptf,  
$0<\mu<\la$, and\/ 
$\dP=\bigcup_{\al<\mu}\dP_\al\ssa D\dP_\mu$, then\/
$\dP\ssa D\dQ=\bigcup_{\mu\le\al<\la}\dP_\al\;.$
\een                                             
\ele
\bpf
\ref{pqm0}
To see that $D$ is pre-dense in $\dP\cup\dQ$, 
let $T_0\in\dP\cup\dQ$. 
By \ref{fm}\ref{fm1}, there is a tree $T\in\dQ$, $T\sq T_0$. 
Then $T\sqf \bigcup D$, in particular, there is a tree 
$S\in D$ with $X=[S]\cap[T]\ne\pu$.  
However $X$ is clopen in $[T]$ by \ref{fm}\ref{fm4}.
Therefore 
there is a tree $T'\in\dQ$ 
with $[T']\sq X$, thus $T'\sq S\in D$ and $T'\sq T\sq T_0$. 
We conclude that $T_0$ is compatible with $S\in D$ in $\dP\cup\dQ$.

To see that $D$ is pre-dense in $\dP$ (assuming $\dP$ is regular), 
let $S_0\in\dP$. 
It follows from the above that $S_0$ is compatible with some 
$S\in D$, hence, $S$ and $S_0$ are not \ad. 
It remains to use  Lemma~\ref{regfn}.  

To prove \ref{pqm1} on the top of Lemma~\ref{pqr}\ref{pqr1}, 
let $R\in\dR$. 
Then $R\sqf\bigcup\dQ$, but each $T\in\dQ$ satisfies 
$T\sqf\bigcup D$.
The same for \ref{pqm2}.
\epf


\parf{Refining \muf s} 
\las{refmf}

Let $\jpi,\jqo$ be \muf s. 
Say that $\jqo$ is an \rit{refinement} of $\jpi$, 
\index{multiforcing!refinement, $\jpi\bssq\jqo$}%
\index{zzzpissqqo@$\jpi\bssq\jqo$}%
\kmar{bssq}%
symbolically $\jpi\bssq\jqo$, if $\abc\jpi\sq\abc\jqo$ 
and 
$\jpi(\xi)\ssq\jqo(\xi)$ whenever $\xi\in\abc\jpi$.

\bcor
[of Lemma~\ref{pqr}]
\lam{pqrC}
If\/ $\jpi\bssq\jqo\bssq\jro$ then\/ $\jpi\bssq\jro$. 

If\/ $\jpi\bssq\jqo$ then the set\/ $\mt\jqo$ 
is open dense\snos
{\label{fnmden}%
If $\jP\sq \jR\sq\md{}{}$ then, similarly to 
Footnote~\ref{fnden}, $\jP$ is 
1) 
\rit{dense} in $\jR$ iff 
\index{set!dense (\mut s)}%
\index{dense (\mut s)}%
$\kaz\zr\in \jR\,\sus\zp\in \jP\,(\zp\leq\zr)$, 
2) 
\rit{open dense} in $\jR$ iff in addition 
\index{set!open dense (\mut s)}%
\index{open dense (\mut s)}%
$\kaz\zr\in \jR\,\kaz \zp\in \jP\,(\zp\leq\zr\imp \zp\in\jR)$, 
and 3) 
\rit{pre-dense} in $\jR$ iff the set  
\index{set!pre-dense (\mut s)}%
\index{pre-dense (\mut s)}%
$\jP'=\ens{\zr\in\jR}{\sus\zp\in\jP(\zr\leq \zp)}$ 
is dense in $\jR$.
} 
in\/ $\mt{\jpi\kw\jqo}$. 
\qed
\ecor

Our next goal is to introduce a version of 
Definition~\ref{ssqm} suitable for \muf s; 
we expect an appropriate version of 
Lemma~\ref{pqm} to hold. 

First of all, we accomodate the definition of the 
relation $\sqf$ in Section~\ref{em} for \mut s. 
Namely if $\ju$ is a \mut\ and $\zD$ a collection of \mut s, 
then $\ju\sqf\bigvee \zD$ 
\kmar{sqf}%
will mean that there is a finite 
\index{zzsqf@$\sqf$}%
\index{zzsqfV@$\sqf\bigvee$}%
set $\zD'\sq \zD$ satisfying 
1) $\abc \jv=\abc \ju$ for all $\jv\in \zD'$, and 
2) $[\ju]\sq\bigcup_{\jv\in \zD'}[\jv]$.  

\bdf
\lam{ssl}
Let\/ $\jpi,\jqo$ be \muf s, and\/ $\jpi\bssq\jqo$. 
Say that\/ 
\rit{$\jqo$ seals a set\/ $\zD\sq\mt\jpi$ over\/ $\jpi$}, 
\index{seals}%
symbolically\/ $\jpi\ssb\zD\jqo$,
\kmar{ssb zD}
\index{multiforcing!refinement2D@refinement, $\jpi\ssb\zD\jqo$}%
\index{refinement2D@refinement, $\jpi\ssb\zD\jqo$}%
\index{zzzpissq2Dqo@$\jpi\ssb\zD\jqo$}%
if the following condition holds:
\ben
\fenu
\itlb{ssl*}%
if\/ 
$\zp\in\mt\jpi$, $\ju\in\mt\jqo$, $\abc\ju\sq\abc\jpi$, 
$\abc\ju\cap\abc\zp=\pu$, 
then there is\/ $\zq\in\mt\jpi$  
such that\/ $\zq\leq\zp$, still\/ $\abc{\zq}\cap\abc\ju=\pu$,  
and\/ $\ju\sqf \bigvee \duz \zD\ju\zq$, where\/ 
$$
\duz \zD\ju\zq=\ens{\ju'\in\mt\jpi}
{\abc{\ju'}=\abc\ju\:\text{ and }\:\ju'\cup\zq\in \zD}\,. 
\index{zzDuq@$\duz \zD\ju\zq$}%
\eqno\qed
$$
\een
\eDf

Note that if $\zp,\ju,\zD,\zq$ are as indicated then still 
$\ju\cup\zq\sqf \bigvee\zD'$, where  
$\zD'=\ens{\ju'\cup\zq}{\ju'\in\duz \zD\ju\zq}\sq\zD$. 
Anyway the definition of $\ssb\zD$ in \ref{ssl}
looks somewhat different and more complex 
than the definition of $\ssa D$ in \ref{ssqm}, 
which reflects the fact 
that finite-support products of
forcing notions in \ptf\ behave 
differently (and in more complex way) 
than single arboreal forcings.
Accordingly, the next lemma, similar to Lemma~\ref{pqm}, is 
somewhat less obvious.

\ble
\lam{pqn}
Let\/ $\jpi,\jqo,\jsg$ be \muf s and\/
$\zD\sq\mt\jpi$.
Then$:$
\ben
\renu
\itlb{pqn0}%
if\/ $\jpi\ssb \zD\jqo$ then\/ $\zD$ is 
dense in\/ $\mt\jpi$ and
pre-dense in\/ $\mt{\jpq}\;;$

\vyk{
\itlb{pqn0a}%
if\/ $\jpi\ssb \zD\jqo$ and\/ $\zD\sq\zD'\sq\mt\jpi$   
then\/ $\jpi\ssb{\zD'}\jqo\;;$
}

\itlb{pqn0b}%
if\/ $\jpi$ is regular, 
$\jpi\ssb{\zD_i}\jqo$ for\/ $i=1,\dots,n$, all sets\/ 
$\zD_i\sq\mt\jpi$ 
are open dense in\/ $\mt\jpi$, and\/ $\zD=\bigcap_i\zD_i$,   
then\/ $\jpi\ssb{\zD}\jqo\;;$

\itlb{pqn1}%
if\/ ${\zD}$ is open dense in\/ $\mt\jpi$ 
and\/ $\jpi\ssb {\zD}\jqo\bssq\jsg$ then\/ 
$\jpi\ssb {\zD}\jsg$$;$

\itlb{pqn2}%
if\/ $\sis{\jpi_\al}{\al<\la}$  
is a\/ \dd\bssq increasing sequence in\/ \mfp, 
$0<\mu<\la$, $\jpi=\bkw_{\al<\mu}\jpi_\al$, 
${\zD}$ is open dense in\/ $\mt\jpi$, 
and\/ $\jpi\ssb {\zD}\jpi_\mu$, then\/
$\jpi\ssb {\zD}\jqo=\bkw_{\mu\le\al<\la}\jpi_\al$.
\een                                             
\ele
\bpf
\ref{pqn0}
To check that ${\zD}$ is pre-dense in $\mt\jpq$, 
let $\zr\in \mt\jpq$. 
Due to the product character of $\mt\jpq$, we can assume 
that $\abc\zr\sq\abc\jpi$. 
Let  
$$
X=\ens{\xi\in\abc\zr}{\zc\zr \xi \in\mt\jqo}\,, 
\quad 
Y=\ens{\xi\in\abc\zr}{\zc\zr \xi\in\mt\jpi}\,. 
$$
Then $\zr=\ju\cup \zp$, where 
$\ju=\zr\res X\in\mt\jqo$, $\zp=\zr\res Y\in\mt\jpi$. 
As $\jqo$ seals ${\zD}$, there is a \mut\ $\zq\in\mt\jpi$  
such that $\zq\leq\zp$, $\abc{\zq}\cap\abc\ju=\pu$,  
and $\ju\sqf \bigvee\duz {\zD}\ju\zq$.
By an easy argument, 
there is a \mut\ $\ju'\in \duz {\zD}\ju\zq$ compatible 
with $\ju$ in $\mt\jqo$; let $\jw\in\mt\jqo$, $\jw\leq\ju$, 
$\jw\leq\ju'$, $\abc\jw=\abc{\ju'}=\abc\ju$. 
Then the \mut\ $\zr'=\jw\cup\zq\in\mt\bbpr$ 
satisfies $\zr'\leq\zr$ and $\zr'\leq\ju'\cup\zq\in {\zD}$. 

To check that ${\zD}$ is dense in $\mt\jpi$, 
suppose that $\zp\in\mt\jpi$. 
Let $\ju=\jLa$ (the empty \mut) 
in \ref{ssl*} of Definition~\ref{ssl}, 
so that $\abc\ju=\pu$ and $\duz {\zD}\ju\zq={\zD}$.\vom

\ref{pqn0b} 
Let $\zp\in\mt\jpi$, 
$\ju\in\mt\jqo$, $\abc\ju\sq\abc\jpi$, 
$\abc\ju\cap\abc\zp=\pu$. 
Iterating \ref{ssl*} for $\zD_i$, $i=1,\dots,n$, we find  
a \mut\ $\zq\in\mt\jpi$  
such that\/ $\zq\leq\zp$, $\abc{\zq}\cap\abc\ju=\pu$,  
and\/ $\ju\sqf \bigvee \duz{(\zD_i)}\ju\zq$ for all $i$, where\/ 
$$
\duz{(\zD_i)}\ju\zq=\ens{\ju'\in\mt\jpi}
{\abc{\ju'}=\abc\ju\:\text{ and }\:\ju'\cup\zq\in \zD_i}\,. 
$$
Thus there are finite sets $U_i\sq \duz{(\zD_i)}\ju\zq$ 
such that $[\ju]\sq\bigcup_{\jv\in U_i}[\jv]$ for all $i$. 
Using the regularity assumption and Corollary~\ref{lsad},
we get a finite set $W\sq \mt\jpi$  
such that $\abc{\jw}=\abc\ju$ for all $\jw\in W$, 
$\bigcap_i\bigcup_{\jv\in U_i}[\jv]=\bigcup_{\jw\in W}[\jw]$, 
and if $i=1,\dots,n$ and $\jw\in W$ then   
$[\jw]\sq [\jv]$ for some $\jv\in U_i$ --- hence
$\jw\cup\zq\in \zD_i$.
We conclude that if $\jw\in W$ then $\jw\cup\zq\in \zD$,
hence $\jw\in\duz{\zD}\ju\zq$.
Thus $W\sq\duz{\zD}\ju\zq$. 
However $[\ju]\sq\bigcup_{\jw\in W}[\jw]$ by the choice of $W$.
Thus $\ju\sqf \bigvee \duz{\zD}\ju\zq$.
\vom

\ref{pqn1}
We have $\jpi\bssq\jsg$ by Corollary~\ref{pqrC}. 
To check that $\jsg$ seals ${\zD}$ over $\jpi$,
let $\ju\in\mt\jsg$,  
$\abc\ju\sq\abc\jpi$, 
$\zp\in\mt\jpi$, $\abc\ju\cap\abc\zp=\pu$.
As $\jqo\bssq\jsg$, there is a finite $U\sq\mt\jqo$ such that 
$\abc\jv=\abc\ju$ for all $\jv\in U$, and 
$[\ju]\sq\bigcup_{\jv\in U}[\jv]$.
As $\jpi\ssb\zD\jqo$, 
by iterated application of Definition~\ref{ssl}\ref{ssl*}, 
we get a \mut\ $\zq\in\mt\jpi$  
such that $\zq\leq\zp$, $\abc{\zq}\cap\abc{\ju}=\pu$,  
and if $\jv\in U$ then $\jv\sqf\bigvee\duz {\zD}{\ju}\zq$, 
where\pagebreak[0] 
$$
\duz{\zD}{\ju}\zq=\ens{\jv'\in\mt\jpi}
{\abc{\jv'}=\abc{\jv}=\abc{\ju}\land{\jv'}\cup\zq\in {\zD}}\,.
$$
Note finally that $\ju\sqf\bigvee U$ by construction, hence 
$\ju\sqf\bigvee\duz {\zD}{\ju}\zq$ as well.\vom

\ref{pqn2}
We have to check that $\jqo$ seals ${\zD}$ over $\jpi$.
Let $\ju\in\mt\jqo$, 
$\abc\ju
\sq\abc\jpi$, 
$\zp\in\mt\jpi$, $\abc\ju\cap\abc\zp=\pu$.
There is a finite set $U\sq\mt{\jpi_\mu}$ such that 
$\abc\jv=\abc\ju$ for all $\jv\in U$ and 
$[\ju]\sq\bigcup_{\jv\in U}[\jv]$.
And so on as in the proof of \ref{pqn1}.
\epf

\parf{Generic refinement of a \muf\ by Jensen} 
\las{jex}

Here we introduce a construction, 
due to Jensen in its original form, which brings
refinements 
of forcings and \muf s, of types $\ssa{D}$ and $\ssb{\zD}$. 

\bdf
\lam{dPhi}
Suppose that 
$\jpi=\sis{\dP_\xi}{\xi\in\abs\jpi}$ is a small \muf.\vom

(0)
Let a \rit{\dd\jpi system} be any indexed set of the form  
$\vpi=\sis{\zd\vpi\xi k}{\ang{\xi,k}\in\abs\vpi}$,  
where $\abs\vpi\sq\abs\jpi\ti\om$ is finite and 
$\zd\vpi\xi k=\vpi(\xi,k)\in \clg{\dP_\xi}$ 
for all $\xi,k$. 
(Recall that $\clg{\dP_\xi}$ consists of all 
finite unions of trees in $\dP_\xi$.)
We order the set $\sys\jpi$ of all 
\dd\jpi systems componentwise: $\vpi\leq\psi$ 
($\vpi$ extends $\psi$) 
iff $\abs\psi\sq\abs\vpi$ and 
$\zd\vpi\xi k\sq\zd\psi\xi k$ 
for all $\ang{\xi,k}\in\abs\psi$.
Accordingly the set $\om\ti\sys\jpi$ 
is ordered so that  
$\ang{n,\vpi}\lec\ang{m,\psi}$ 
iff 
$\abs{\psi}\sq\abs{\vpi}$ and 
$\bng{n,\zd\vpi\xi k}\lec\bng{m,\zd{\psi}\xi k}$ 
in $\om\ti\pet$ (Section~\ref{tre}) 
for all $\xi,k$; this implies $m\le n$.\vom

(1) 
Let 
$\cM\in\hc$ be any set.\snos
{\label{hc}%
Recall that $\hc$ = all \rit{hereditarily countable}  
\index{set!hereditarily countable, $\hc$}%
\index{hereditarily countable set, $\hc$}%
\index{zzHC@$\hc$}%
sets, \ie\ those having at most countable transitive 
closures.}The set $\cmp$ of all sets $X\in\hc$, \dd\in definable 
in $\hc$ by formulas with sets in $\cM$ as parameters, is 
still countable. 
Therefore there exists a \dd\lec decreasing sequence 
$\dphi=\sis{\ang{n_j,\vpi_j}}{j<\om}$ 
of pairs $\ang{n_j,\vpi_j}\in\om\ti\sys\jpi$, 
 \dd{\cmp}{\it generic\/} in the sense that it intersects 
every set $D\in\cM\yd D\sq\om\ti\sys\jpi$, 
open dense in $\om\ti\sys\jpi$.\snos
{The density means that for any 
$\ang{m,\psi}\in\om\ti\sys\jpi$ 
there is $\ang{n,\vpi}\in D$ with 
$\ang{n,\vpi}\lec\ang{m,\psi}$. 
The openness means that if $\ang{m,\psi}\in D$ and 
$\ang{n,\vpi}\lec\ang{m,\psi}$ then $\ang{n,\vpi}\in D$.} 
Let us fix any such a \dd{\cmp}generic  sequence 
$\dphi$. 

By definition, each $\vpi_j$ has the form 
$\vpi_j=\sis{\zd{\vpi_j}\xi k}
{\ang{\xi,k}\in\abs{\vpi_j}}$, 
where $\abs{\vpi_j}\sq\abs\jpi\ti\om$ is finite, and each 
tree $\zd{\vpi_j}\xi k$ belongs to $\clg{\dP_\xi}$. 
We have $n_j\to\iy$ by the genericity, so that it 
can be wlog 
assumed that {\ubf $n_0<n_1<n_2<\ldots$ strictly}. \vom

(2) 
Let $\xi\in\abs\jpi$, $k<\om$. 
By the genericity assumption, 
there is a number $j(\xi,k)$ such that 
if $j\ge j(\xi,k)$ then 
$\ang{\xi,k}\in\abs{\vpi_j}$, hence the tree  
$\vpi_j(\xi,k)=\zd{\vpi_j}\xi k\in\clg{\dP_\xi}$
is defined, and we have 
$$
\ldots\lec
\ang{n_{j(\xi,k)+2},\zd{\vpi_{j(\xi,k)+2}}\xi k}
\lec\ang{n_{j(\xi,k)+1},\zd{\vpi_{j(\xi,k)+1}}\xi k}
\lec\ang{n_{j(\xi,k)},\zd{\vpi_{j(\xi,k)}}\xi k}\,, 
$$
with $n_{j(\xi,k)}<n_{j(\xi,k)+1}<n_{j(\xi,k)+2}<\ldots$ 
strictly, by (1) above.\vom

(3)
Then it follows by Lemma~\ref{fus} that each intersection  
$ 
\TS
\qfj{\xi}k
=\bigcap_{j\ge j(\xi,k)}\zd{\vpi_{j}}\xi k
$ 
\index{tree!QFixik@$\qfj{\xi}k$}%
\index{zzQFixik@$\qfj{\xi}k$}%
\kmar{qfj xi k}%
is a tree in $\pes$ (not necessarily in $\dP_\xi$), 
and the relation
$ 
\ang{n_{j},\qfj{\xi}k}
\lec
\ang{n_{j},\zd{\vpi_{j}}\xi k}
$ 
holds for all $j\ge j(\xi,k)$. 
We define 
\kmar{cqf xi}%
\index{zzQFixi@$\cqf\xi$}%
$\cqf\xi=\ens{{\qfj{\xi}k}\ret s}
{k<\om\land s\in{\qfj{\xi}k}}$.\vom

(4) 
We finally let $\jqo=\sis{\cqf\xi}{\xi\in\abs\jpi}$ and
$\jpq=\sis{\dP_\xi\cup\cqf\xi}{\xi\in\abs\jpi}$.\vom

(5)  
Finally if $\jqo=\jqo[\dphi]$ is obtained this way from an 
\dd{\cmp}generic sequence $\dphi$, 
then $\jqo$ is called 
an \dd\cM\rit{generic refinement of $\jpi$}. 
\index{multiforcing!refinement!generic}%
\index{refinement!generic}%
\edf

\ble
[by the countability of $\cmp$]
\lam{gee} 
If\/ 
$\jpi$ is a small \muf\ 
and\/ $\cM\in\hc$ 
then there is an\/ 
\dd\cM generic refinement\/ $\jqo$ of\/ $\jpi$.
\qed
\ele

\bte
\lam{dj}
If\/ $\cM\in\hc$ is transitive, 
$\jpi=\sis{\dP_\xi}{\xi\in\abs\jpi}\in\cM$ is a small \muf, 
and\/ $\jqo=\jqo[\dphi]=\sis{\dQ_\xi}{\xi\in\abs\jpi}$ 
is an\/ \dd\cM generic refinement of\/ $\jpi$, 
then$:$
\ben
\renu
\itlb{dj1}%
$\jqo$ is a small \mdi\ \muf, $\abc\jqo=\abc\jpi$, 
and\/ $\jpi\bssq\jqo\;;$%

\itlb{dj7}%
if pairs\/ $\ang{\xi,k}\ne\ang{\et,\ell}$ 
belong to\/ $\abs\jpi=\abs\jqo$ then\/ 
$[\qfj\xi k]\cap[\qfj\et \ell]=\pu\;;$ 

\itlb{dj8}%
if\/ $\xi\in\abs\jpi$, $S\in\dQ_\xi$ and\/ $T\in\dP_\xi$ 
then\/ $[S]\cap[T]$ is clopen in\/ $[S]$ 
and\/ $T\not\sq S$,  
in particular,\/ $\dQ_\xi\cap \dP_\xi=\pu\;;$  

\itlb{dj9}%
if\/ $\xi\in\abs\jpi$ then 
the set\/ $\dQ_\xi$ is open dense in\/ $\dQ_\xi\cup\dP_\xi\;;$

\itlb{dj3}%
if\/ $\xi\in\abs\jpi$ and a set\/ $D\in\cM$, 
$D\sq\dP_\xi$ is pre-dense in\/ $\dP_\xi$ then\/ 
$\dP_\xi\ssa D\dQ_\xi\;;$


\itlb{dj2}%
if in addition\/ $\jpi=\bkw_{\al<\la}\nor\jpi\al$, 
where\/ $\la<\omi$ and\/ 
$\sis{\nor\jpi\al}{\al<\la}\in\cM$ 
is a\/ \dd\bssq increasing sequence of small \mdi\ \muf s, 
then\/ $\nor\jpi\al\bssq\jqo$ for all\/ $\al<\la$.
\een
\ete  

\bpf
We argue in the notation of Definition~\ref{dPhi}.

\ref{dj7}  
By Corollary~\ref{suz}\ref{suz2}, the set $D$ of all 
pairs $\ang{n,\vpi}\in\om\ti\sys\jpi$, 
where $\vpi$ 
is a pairwise \ad\ system and $\abs\vpi$ 
contains both $\ang{\xi,k}$, $\ang{\et,\ell}$,  
is dense in $\om\ti\sys\jpi$, and obviously 
$D\in\cmp$. 
Thus   
$\ang{n_j,\vpi_j}\in D$ for some $j<\om$. 
Then $\zd{\vpi_j}\xi k\cap\zd{\vpi_j}\et \ell=\pu$ 
since $\vpi_j$ is \ad. 
But $\qfj\xi k\sq \zd{\vpi_j}\xi k$, 
$\qfj\et \ell\sq \zd{\vpi_j}\et \ell$
by construction.%

\ref{dj8}
Let $S=\qfj\xi k$.
To prove the clopenness claim, note that 
the set $D(T)$ of all 
pairs $\ang{n,\vpi}\in\om\ti\sys\jpi$, such that 
$\ang{\xi,k}\in\abc\vpi$ and if $s\in2^n$ 
then either ${\zd\vpi\xi k}\ret s\sq T$ 
or $[\zd\vpi\xi k]\cap [T]=\pu$, 
is dense in $\om\ti\sys\jpi$. 
To prove $T\not\sq S$,  
similarly the set $D'(T)$ of all 
pairs $\ang{n,\vpi}\in\om\ti\sys\jpi$, such that 
$\ang{\xi,k}\in\abc\vpi$ and 
$T\not\sq \zd\vpi\xi k$, is dense.
Note that $D(T),D'(T)\in{\cmp}$ and argue as above.

\ref{dj9}
The openness easily follows from \ref{dj8}. 
To prove the density, let $T\in\dP_\xi$.  
The set $\Da(T)$ of all pairs 
$\ang{n,\vpi}\in\om\ti\sys\jpi$, 
such that $\ang{\xi,k}\in\abs\vpi$ and 
$\zd{\vpi}\xi k=T$ for some $k$, 
belongs to $\cmp$ 
and is dense in $\om\ti\sys\jpi$. 

\ref{dj1} 
By construction, the sets $\jqo(\xi)=\cqf\xi$ are 
special arboreal forcings, and hence 
$\jqo$ is a small \mdi\ \muf, 
and $\abc\jqo=\abc\jpi$. 
To establish $\jpi\bssq\jqo$, 
let $\xi\in\abc{\jpi}$.
We have to prove that $\dP_\xi\ssq\dQ_\xi$. 
Condition \ref{fm1} of Definition~\ref{fm} follows from 
\ref{dj9}, 
condition \ref{fm4} from \ref{dj8}, and \ref{fm3} 
holds since 
$\qfj\xi k\sq \zd{\vpi_j}{\xi}k \in \clg{\dP_\xi}$ 
for some $j$.

\ref{dj3} 
Assume that $\xi\in\abs\jpi$, $k<\om$, $D\in\cmp$ is 
pre-dense in $\dP_\xi$. 
Then the set 
$D'=\ens{T\in\dP_\xi}{\sus S\in D(T\sq S)}$ 
is open dense in $\dP_\xi$, and hence the set 
$\Da\in\cmp$ of all pairs 
$\ang{n,\vpi}\in\om\ti\sys\jpi$, 
such that $\ang{\xi,k}\in\abs\vpi$ and 
$\zd{\vpi}\xi k\ret s\in D'$ 
for all $s\in2^n\cap\zd{\vpi}\xi k$,
is dense in $\om\ti\sys\jpi$ by Lemma~\ref{nq}. 
Thus $\ang{n_j,\vpi_j}\in \Da$ for some $j$, 
and this implies   
$\qfj\xi k\sq \zd{\vpi_j}{\xi}k\sqf\bigcup D$.

\vyk{
and still $D'\in{\cmp}$.
To prove 
$\qfj\xi k\sqf\bigcup D$ for all . 
Then . 
}


\ref{dj2} 
We have to prove that $\nor\jpi\al(\xi)\ssq\jqo(\xi)$
whenever $\xi\in \abc{\nor\jpi\al}$. 
And as $\jpi(\xi)\ssq\jqo(\xi)$ has been checked, 
it suffices to prove that 
$\qfj\xi k\sqf\bigcup \nor\jpi\al(\xi)$. 
However $D=\nor\jpi\al(\xi)$ is pre-dense in 
$\jpi(\xi)=\dP_\xi$ by Lemma~\ref{pqr}\ref{pqr3}, 
and still $D\in{\cmp}$, hence we can refer to 
\ref{dj3}.
\epf

\bcor
\lam{xistt}
In the assumptions of Lemma~\ref{gee}, if\/ 
$\abc\jpi\sq Z\sq\omi$ and\/ $Z$ is countable 
then there is a small \mdi\ \muf\/ $\jqo$ 
such that\/ $\abc\jqo=Z$ 
and\/ $\jpi\bssq\jqo$.
\ecor
\bpf
If $\abc\jpi=Z$ then let $\cM\in\hc$ 
be any countable set containing 
$\jpi$, pick $\jqo$ by Lemma~\ref{gee}, 
and apply Theorem~\ref{dj}. 
If $\abc\jpi\sneq Z$ then we trivially extend the 
construction by $\jqo(\xi)=\dpo$ (see Example~\ref{cloL}) 
for all $\xi\in Z\bez\abc\jpi$.
\epf  

\vyk{
\bcor
\lam{xistC}
Suppose that\/ $\la<\omi$, and\/ 
$\sis{\dP_\al}{\al<\la}$ is an\/ 
\dd\ssq increasing sequence of countable special forcings 
in\/ $\ptf$. 
Then there is a countable special forcing\/ $\dQ\in\ptf$ 
such that\/ $\dP_\al\ssq\dQ$ for each\/ $\al<\la\;.$ 
\ecor
\bpf
If $\al<\la$ then let a \muf\ $\nor\jpi\al$ be defined 
by $\abc{\nor\jpi\al}=\ans{0}$ and 
${\nor\jpi\al}(0)=\dP_\al$. 
By Lemma~\ref{gee} 
and Theorem~\ref{dj} there is a \muf\ $\jqo$ satisfying  
$\abc\jqo=\ans{0}$ and   
${\nor\jpi\al}\bssq\jqo$, $\kaz\al$. 
Let $\dQ=\jqo(0)$.
\epf
}

\parf{Generic refinement: sealing dense sets} 
\las{pres}

This Section proves a special consequence of 
\dd{\cmp}genericity of \muf\ refinements, 
the relation $\ssb{\zD}$ of Definition~\ref{ssl} between
a \muf\ and its refinement, via a dense set $\zD$.

\bte
\lam{uu4}
Under the assumptions of Theorem~\ref{dj}, if\/ 
${\zD}\in\cmp$, ${\zD}\sq\mt\jpi$, and\/ 
${\zD}$ is open dense in $\mt\jpi$,
then $\jpi\ssb {\zD}\jqo$. 
\ete
\bpf
We suppose that $\jqo=\jqo[\dphi]$ is obtained from an
decreasing 
\dd{\cmp}generic sequence $\dphi$ of pairs 
$\ang{n_j,\vpi_j}\in\om\ti\sys\jpi$,
as in Definition~\ref{dPhi}(1), 
and argue in the notation of \ref{dPhi}.
Suppose that $\zp\in\mt\jpi$, $\ju\in\mt\jqo$, 
$\abc\ju\cap\abc\zp=\pu$, as in \ref{ssl*} of 
Definition~\ref{ssl}; the extra condition 
$\abc\ju\sq\abc\jpi$ holds automatically as we have 
$\abc\jqo=\abc\jpi$. 
We have to find a \mut\ $\zq$ which witnesses 
\ref{ssl}\ref{ssl*} for $\ju$.

Each term $\zc\ju\xi$ of $\ju$ 
($\xi\in\abs\ju$) 
is equal to some $\qfj\xi{,k_\xi}\ret{t_{\xi}}$, 
where $t_{\xi}\in\qfj\xi{,k_\xi}$. 
We can wlog assume that simply $t_\xi=\La$, 
so that $\zc\ju\xi=\qfj\xi{,k_\xi}$, $\kaz \xi$.  

\bdf
\lam{uu4d}
If $n<\om$ then 
let $\syn{\jpi}n$ contain all systems $\vpi\in\sys\jpi$ 
such that  
$\ang{\xi,k_\xi}\in\abs\vpi$ for all $\xi\in\abs\ju$, 
and $\zd\vpi\xi{k}\ret t\in\dP_\xi=\jpi(\xi)$ 
(not just $\in\clg{\dP_\xi}$!)
for all $\ang{\xi,k}\in\abs\vpi$ and 
$t\in2^n\cap\zd\vpi\xi{k}$.
If $\vpi\in\syn{\jpi}{n}$ then 
let 
$\bsn{n}{\vpi}
$ 
contain all multistrings 
$\bs=\sis{s_\xi}{\xi\in\abs\ju}$ such that 
$s_\xi\in2^n\cap \zd\vpi\xi{,k_\xi}$, $\kaz\xi\in\abs\ju$. 
If $\bs=\sis{s_\xi}{\xi\in\abs\ju}\in\bsn{n}{\vpi}$ 
then define 
$\jvn{\bs}{\vpi}\in\mt\jpi$ by 
$\abs{\jvn{\bs}{\vpi}}=\abs\ju$ 
and $\zc{\jvn{\bs}{\vpi}}\xi=\zd\vpi\xi{,k_{\xi}}\ret{s_\xi}$ 
for all $\xi\in\abs\ju$.
\edf

\ble
\lam{syn}
Let\/ $n<\om$ and\/ $\vpi\in\sys\jpi$. 
There exists a system\/ $\psi\in\syn\jpi n$ satisfying\/ 
$\ang{n,\psi}\lec\ang{n,\vpi}$. 
\ele
\bpf
Add each absent 
$\ang{\xi,k_{\xi}}\nin\abs\vpi$ 
to $\abs\psi$ and define  
$\zd\psi\xi{,k_{\xi}}\in{\dP_\xi}$ arbitrarily.
If $\ang{\xi,k}\in\abs{\psi}$ 
and $t\in2^n\cap\zd\psi\xi{k}$, but  
$\zd\psi\xi{k}\ret t\in\clg{\dP_\xi}\bez\dP_\xi$,
then shrink $\zd\psi\xi{k}$ to a tree in $\dP_\xi$ 
by Lemma~\ref{suz}\ref{suz0}, and do this 
for all triples $\xi,k,t$ as indicated. 
\epF{Lemma}


\ble
\lam{uu4a}
If\/ $\zr\in\mt\jpi$, $\abs\zr\cap\abs\ju=\pu$, then 
the set\/ $\Da_{\zr}\in\cM$ of all pairs\/   
$\ang{n,\vpi}\in \om\ti\sys\jpi$, 
such that\/ $\vpi\in\syn\jpi n$  and 
there is\/ $\zq\in \mt\jpi$ 
satisfying\/ $\zq\leq\zr$, $\abs \ju\cap\abs \zq=\pu$, 
and\/ 
$(1)$ 
if\/ $\bs\in\bsn{n}{\vpi}$ then\/ 
$\jvn{\bs}{\vpi}\cup\zq\in \zD$, ---  
is dense in\/ $\om\ti\sys\jpi$.
\ele 
\bpf[Lemma]
Let $\ang{n,\psi}\in\om\ti\sys\jpi$. 
We'll find 
a pair $\ang{n,\vpi}\in \Da_{\zr}$ (same $n$!) 
with 
$\ang{n,\vpi}\lec\ang{n,\psi}$. 
We wlog assume that $\psi\in\syn\jpi n$, by Lemma~\ref{syn}.

Let $\bs=\sis{s_\xi}{\xi\in\abs\ju}\in\bsn{n}{\psi}$. 
Consider the \mut\ $\jvn{\bs}{\psi}\in\mt\jpi$.
As $\zD$ is dense, 
there are \mut s $\zr',\jv\in\mt\jpi$ 
such that $\abs{\jv}=\abs\ju$, $\jv\leq\jvn{\bs}{\psi}$, 
$\abs{\zr'}\cap\abs\ju=\pu$, $\zr'\leq\zr$, and 
$\jv\cup\zr'\in \zD$.
Define a system  
$\psi'\in \sys\jpi$ with $\abs{\psi'}=\abs\psi$, 
that extends $\psi$ by shrinking each tree 
$\zd\psi\xi{,k_{\xi}}\ret{s_\xi}$ 
to $\zc\jv\xi$, so that 
$\zd{\psi'}\xi{,k_{\xi}}\ret{s_\xi}=\zc\jv\xi$, 
but 
$\zd{\psi'}\xi{,k_{\xi}}\ret{t}=
\zd{\psi}\xi{,k_{\xi}}\ret{t}$ 
for all $t\in2^n\cap\zd\psi\xi{,k_{\xi}}$, $t\ne s_\xi$, 
and $\zd{\psi'}\et{k}=\zd{\psi}\et{k}$ whenever 
$\ang{\et,k}\in\abs\psi$ does not have 
the form $\ang{\xi,k_\xi}$, where $\xi\in\abs\ju$.
We have $\ang{n,\psi'}\lec\ang{n,\psi}$ by construction, 
therefore $\bsn{n}{\psi'}=\bsn{n}{\psi}$.

This construction can be iterated, so that all strings 
$\bs\in\bsn{n}{\psi}$ are considered one by one. 
This results in a system $\vpi\in\sys\jpi$, such that 
$\abs\vpi=\abs\psi$ and 
$\ang{n,\vpi}\lec\ang{n,\psi}$ --- and then 
$\bsn{n}{\vpi}=\bsn{n}{\psi}$, 
and a \mut\ $\zq\in\mt\jpi$ with $\zq\leq\zr$ and 
still $\abs\zq\cap\abs\ju=\pu$, such that if 
$\bs
\in\bsn{n}{\psi}$ 
then the \mut\ $\jvn{\bs}{\psi}$, 
satisfies $\jvn{\bs}{\psi}\cup\zq\in \zD$.
Then $\zq$ witnesses that $\ang{n,\vpi}\in\Da_r$.
\epF{Lemma}  

By the lemma, we have $\ang{n_j,\vpi_j}\in\Da_{\zp}$ 
for some $j$.
Let this be witnessed by a \mut\ $\zq\in\mt\jpi$,  
so that $\zq\leq\zp$, $\abs \ju\cap\abs \zq=\pu$, 
and (1) of Lemma~\ref{uu4a} holds 
for $n=n_j$, $\vpi=\vpi_j$.
We easily conclude that 
$[\ju]\sq\bigcup_{\bs\in\bsn{n}{\vpi_j}}[\jvn{\bs}{{\vpi_j}}]$. 
Yet $\jvn{\bs}{{\vpi_j}}\in\duz \zD\ju\zq$, 
$\kaz\bs$, by (1).
\epF{Theorem}

\bcor
\lam{uu4c}
Under the assumptions of Theorem~\ref{dj}, if\/ 
a set\/ $\zD\in\cM,$ $\zD\sq\mt\jpi$  
is pre-dense in\/ $\mt\jpi$, 
then it remains pre-dense in\/ $\mt\jpq$. 
\ecor
\bpf
Assume wlog that $\zD$ is open dense in $\mt\jpi$. 
(If not then consider 
$\zD'=\ens{\zp\in \mt\jpi}{\sus\zq\in\zD\,(\zp\leq\zq)}$.) 
Note that $\jpi\ssb \zD\jqo$ by Theorem~\ref{uu4}, 
and use Lemma~\ref{pqn}\ref{pqn0}.
\epf

\parf{Real names and direct forcing} 
\las{rn}

Our next goal is to 
introduce a suitable notation related to names of reals 
in $\dn$ in the context of forcing notions of the form  
$\mt\jpi$.

\bdf
\lam{rk'}
A \rit{real name} 
\index{real name}%
is any set $\rc\sq\md\ti(\om\ti2)$ such that the 
sets $\kkc ni=\ens{\zp\in\md}{\ang{\zp,n,i}\in\rc}$ 
\imar{kkc ni}%
\index{zzKcni@$\kkc ni$}%
satisfy the following: 
if $n<\om$ and $\zp\in \kkc n0$, 
$\zq\in \kkc n1$, then 
$\zp,\zq$ are \sad.\snos 
{Recall that the condition of 
somewhere almost disjointness \sad\ 
(Definition~\ref{dsad}) 
is equivalent to the incompatibility of $\zp,\zq$ 
in $\md$ and in any set of the form $\mt\jpi$, 
where $\jpi$ is a regular \muf, by corollary~\ref{lsad}.} 
Let $\kic n=\kkc n0\cup\kkc n1\sq\mt\jpi$    
\imar{kic n}%
\index{zzKcn@$\kic n$}%

A \qn\ $\rc$ 
\index{real name!small}%
is \rit{small} if each $\kic n$ is at most countable ---  
then the set  
$\abs\rc=\bigcup_{n}\bigcup_{\zp\in\kic n}\abc\zp$, 
\index{zzcII@$\abs\rc$}%
and $\rc$ itself, are countable, too.

Let $\jpi$ be a \muf. 
A \qn\ $\rc$ is 
\rit{\pol\jpi} 
\index{real name!picomplete@$\pi$-complete}%
if every set 
$\kipc n\jpi=\ens{\zp\in\mt\jpi}
{\sus\zq\in \kic n\,(\zp\leq\zq)}$ 
\imar{kipc n jpi}%
\index{zzKcnpi@$\kipc n\jpi$}%
(the \dd\jpi cone of $\kic n$) 
is pre-dense in $\mt\jpi$.
In this case, if a set (a filter) $G\sq\mt\jpi$ is  
\dd{\mt\jpi}generic over the family of 
all sets $\kic n$, then we define a real  
$\rc[G]\in\dn$ so that  $\rc[G](n)=i$ 
\index{real name!evaluation, $\rc[G]$}%
\index{zzcG@$\rc[G]$}%
iff $G\cap \kc ni\ne\pu$.

We do not require in this case that  
$\rc\sq\mt\jpi\ti(\om\ti2)$, or equivalently, 
$\kic n\sq\mt\jpi$ for all $n$, 
but if this inclusion  
indeed holds then this will be explicitly mentioned.
\edf

\vyk{
\bdf
\lam{rk'}
A \rit{real name} 
\index{real name}%
is any set $\rc\sq\md\ti(\om\ti2)$ such that the 
sets $\kkc ni=\ens{\zp\in\md}{\ang{\zp,n,i}\in\rc}$ 
\imar{kkc ni}%
\index{zzKcni@$\kkc ni$}%
satisfy the following: 
\ben
\fenu
\itlb{rk1}%
if $n<\om$ and $\zp\in \kkc n0$, 
$\zq\in \kkc n1$, then 
$\zp,\zq$ are \sad\ 
(somewhere almost disjoint, Definition~\ref{dsad}),  
hence $[\zp]\cap[\zq]=\pu$, so $\zp,\zq$ are incompatible.
\een
Let $\jpi$ be a \muf. 
A \qn\ $\rc$ is said to be a
\rit{\rn\jpi} 
\index{real name!pirealname@\rn\jpi}%
if, in addition to \ref{rk1}, 
$\rc\sq\mt\jpi\ti(\om\ti2)$ and the following 
condition holds:
\ben
\fenu
\atc
\itlb{rk2}%
all sets $\kic n=\kkc n0\cup\kkc n1\sq\mt\jpi$ are 
\imar{kic n}%
\index{zzKcn@$\kic n$}%
pre-dense in $\mt\jpi$.
\vyk{
\rit{pre-dense for\/ $\mt\jpi$}, in the sense that 
\index{set!pre-dense for}%
\index{pre-dense for}%
the set 
$\kipc n\jpi=\ens{\zp\in\mt\jpi}
{\sus\zq\in \kic n\,(\zp\leq\zq)}$ 
\imar{kipc n jpi}%
is dense in $\mt\jpi$. 
}
\een
A \qn\ $\rc$ 
\index{real name!small}%
is \rit{small} if each $\kic n$ is at most countable ---  
then the set  
$\abs\rc=\bigcup_{n}\bigcup_{\zp\in\kic n}\abc\zp$, 
\index{zzcII@$\abs\rc$}%
and $\rc$ itself, are countable, too.
\edf

\vyk{
Generally speaking, we do not assume that 
$\kic n\sq\mt\jpi$. 
However if, in addition to \ref{rk1}, \ref{rk2}  
above, $\kic n\sq\mt\jpi$ holds for all $n$, 
then say 
that $\rc$ is a \rit{true \rn\jpi}.%
\index{real name!true@true \rn\jpi} 
\index{real name!pirealname@\rn\jpi!true}%
Then each set $\kic n=\kkc n0\cup\kkc n1$  
is a pre-dense subset of $\mt\jpi$. 
}

\bdf
\lam{preta}
Let $\rc$ be a \qn\ and
$G\sq\md$ a pairwise compatible set.
Define the
\rit{evaluation} $\rc[G]$, a partial function $\om\to2$, 
so that
\index{real name!evaluation, $\rc[G]$}%
$$
\rc[G](n)=i
\quad\text{iff}\quad 
\index{zzcG@$\rc[G]$}%
\sus \zp\in G\:\sus\zq\in \kkc ni\:(\zp\le\zq)\,.
\eqno\qed
$$
%
\eDf

\bre
\lam{ggem}
Let $\jpi$ be a \muf, $\rc$ be a \rn\jpi, and   
a set $G\sq\mt\jpi$ be \dd{\mt\jpi}generic  over 
the collection of all  sets $\kic n$ as in
\ref{rk2}.
If $n<\om$ then $\kic n$ is pre-dense, thus 
$G\cap(\kic n)\ne\pu$, and hence the value 
$\rc[G](n)=0,1$ is defined --- therefore $\rc[G]\in\dn$.
\ere

\bre
\lam{ai}
If $\jpi$ is \rit{\mre} 
then the notions of being \sad\ 
and being incompatible in $\mt\jpi$ are equivalent by 
Corollary~\ref{lsad}, so that 
a \rn\jpi\ in our sense 
is the same as a \dd{\mt\jpi}name 
for an element of 
$\dn$ in the general theory of forcing.
\ere
}

\vyk{
\bpri
\lam{proj1}
If $\xi<\omi$ then let $\rpi_{\xi}$
\index{real name!x.xi@$\rpi_{\xi}$}%
\index{zzx.xi@$\rpi_{\xi}$}%
be a \qn\ such that if $i=0,1$ then 
the set $\kk ni{\rpi_{\xi}}$ 
consists of a lone \mut\ $\zr=\zr^{\xi}_{ni}$ with 
$\abc{\zr}=\ans{\xi}$ and 
$\zc\zr \xi  =\ens{t\in\bse}{\lh t\le n\lor t(n)=i}$, 
while if $i\ge2$ then $\kk ni{\rpi_{\xi}}=\pu$.
Note that each tree of the form 
$\ens{t\in\bse}{\lh t\le n\lor t(n)=i}$ 
belongs to the set $\dpop$, Example~\ref{cloL}.
\kmar{dpop}

If $\jpi\in\md$ and $\xi\in\abc\jpi$ 
then $\rpi_{\xi}$ is a \rn\jpi\ of the real 
$x_{\xi}=x_{\xi}[G]\in\dn$, 
the $\xi$th term of a \dd{\mt\jpi}generic sequence 
$\sis{x_{\xi}[G]}{\xi\in\abc\jpi}$. 
That is, if $G\sq\mt\jpi$ is generic then 
$x_{\xi}[G]$ introduced by \ref{adds} 
coincides with $\rpi_{\xi}[G]$
defined by \ref{preta}.
\epri
}

Assume that 
$\rc$ is a real name, in the sense of \ref{rk'}. 
Say that a multitree $\zp$: 
\bit
\item
\rit{directly forces\/ $\rc(n)=i$}, 
where $n<\om$ and $i=0,1$, iff there is a 
\mut\ $\zq\in\kkc ni$ 
such that $\zp\leq \zq$; 

\item
\rit{directly forces\/ $s\su\rc$},  
where $s\in\bse,$ iff for all $n<\lh s$, $\zp$ 
directly forces $\rc(n)=i$, where $i=s(n)$;

\item
\rit{directly forces\/ $\rc\nin[T]$},  
where $T\in\pet$, iff there is a string $s\in\bse\bez T$ 
such that $\zp$ directly forces $s\su \rc$. 
\eit
The definition of direct forcing is not 
explicitly associated with 
any concrete forcing notion, but in fact 
it is compatible with any \muf.

\ble
\lam{avo}
Let\/ $\jpi$ be a \muf, $\rc$ a\/ \rn\jpi, $\zp\in\mt\jpi$. 
If\/ $n<\om$ then there exists\/ $i=0,1$ and a \mut\/ 
$\zq\in\mt\jpi\yd \zq\leq\zp$, 
which directly forces\/ $\rc(n)=i$.  
If\/ $T\in\pet$ then there exists\/ $s\in T$ 
and a \mut\/  
$\zq\in\mt\jpi$, $\zq\leq\zp$, which\/  
directly forces\/ $\rc\nin[T\ret s]$. 
\ele

\bpf
To prove the first claim 
use the density of sets $\kic n$ by 
Definition~\ref{rk'} above.
To prove the second claim, 
pick $n$ such that $T\cap2^n$ contains at least two strings.
By the first claim, 
there is a \mut\ $\zq\in\mt\jpi\yd \zq\leq\zp$, 
and a string $t\in T\cap2^n$ such that 
$\zq$ directly forces $t\su \rc$. 
Now take any $s\in T\cap2^n$, $s\ne t$. 
\epf

\parf{Sealing real names and avoiding refinements} 
\las{lorn}

The next definition extends Definition~\ref{ssl}  
to real names. 

\bdf
\lam{ssl+}
Assume that\/ $\jpi,\jqo$ are \muf s, 
$\rc$ is a real name, and\/ $\jpi\bssq\jqo$. 
Say that\/ 
\rit{$\jqo$ seals\/ $\rc$ over\/ $\jpi$}, 
\index{seals}%
symbolically\/ $\jpi\ssb\rc\jqo$, 
\kmar{ssb rc}
\index{multiforcing!refinement2c@refinement, 
$\jpi\ssb\rc\jqo$}%
\index{refinement2c@refinement, $\jpi\ssb\rc\jqo$}%
\index{zzzpissq2cqo@$\jpi\ssb\rc\jqo$}%
if $\jqo$ seals, over $\jpi$, each set 
$\kipc n\jpi=\ens{\zp\in\mt\jpi}
{\sus\zq\in \kic n\,(\zp\leq\zq)}$, 
in the sense of Definition~\ref{ssl}.
\edf

\bcor
\lam{uu4c2}
Under the assumptions of Theorem~\ref{dj}, if\/ 
${\rc}\in\cmp$ and\/ $\rc$ is a\/ \rn\jpi\ 
then $\jpi\ssb {\rc}\jqo$. 
\ecor
\bpf
Each set $\kipc n\jpi$ belongs to $\cmp$
(as so do $\rc$ and $\jpi$) and is 
open dense in $\mt\jpi$, 
so it remains to apply 
Theorem~\ref{uu4}.
\epf

\ble
\lam{pqs}
Let\/ $\jpi,\jqo,\jsg$ be \muf s and\/ $\rc$ be a 
\qn.
Then
\ben
\renu
\itlb{pqs1}%
if\/ $\jpi\ssb \rc\jqo$ then\/ $\rc$ is a\/
\pol\jpi\ and   a\/
\rn{(\jpq)}$;$

\itlb{pqs2}%
if\/ $\jpi\ssb {\rc}\jqo\bssq\jsg$ then\/ 
$\jpi\ssb {\rc}\jsg\;;$ 

\itlb{pqs3}%
if\/ $\sis{\jpi_\al}{\al<\la}$  
is a\/ \dd\bssq increasing sequence in\/ \mfp, 
$0<\mu<\la$, $\jpi=\bkw_{\al<\mu}\jpi_\al$, 
and\/ $\jpi\ssb {\rc}\jpi_\mu$, then\/
$\jpi\ssb {\rc}\jqo=\bkw_{\mu\le\al<\la}\jpi_\al\;.$ 
\een
\ele
\bpf
\ref{pqs1} 
By definition, we have $\jpi\ssb{\kipc n\jpi}\jqo$ 
for each $n$, 
therefore $\kipc n\jpi$ is dense in $\mt\jpi$  
(then obviously open dense) 
and pre-dense in $\mt\jpq$ by Lemma~\ref{pqn}\ref{pqn0}.
It follows that $\kipc n{(\jpq)}$ is dense in $\mt\jpq$.

To check \ref{pqs2}, \ref{pqs3} apply 
\ref{pqn1}, \ref{pqn2} of 
Lemma~\ref{pqn}.
\epf


If $\jpi$ is a \muf\  
then $\mt\jpi$ adds a collection of 
\rit{principal generic reals} $x_{\xi}=x_{\xi}[G]\in\dn$, 
$\xi\in\abc\jpi$, 
where each  
$x_{\xi}$ is \dd{\jpi(\xi)}generic 
over the ground set universe, see Remark~\ref{adds}. 
Obviously many more reals are added, and given a 
\rn\jpi\ $\rc$, one can elaborate different requirements
for a condition $\zp\in\mt\jpi$ to force that $\rc$ is a name 
of a real of the form $x_{\xi k}$ or to force the opposite. 
The next definition provides such a condition related 
to the ``opposite'' direction.

\bdf
\lam{npri}
Let $\jpi$ be a \muf, $\xi\in\abc\jpi$. 
A \qn\  $\rc$ is 
\rit{non-principal over\/ $\jpi$ at\/ $\xi$}, 
if the following set 
is open dense in\/ $\mt\jpi$:
$$
\ddj\xi \rc\jpi
\kmar{ddj xi rc jpi}
=\ens{\zp\in\mt\jpi}{\xi\in\abs\zp\land
\zp\,\text{ directly forces }\,
\rc\nin[\zc\zp\xi]}\,.\eqno\qed
$$
\eDf

We'll show below (Theorem~\ref{npn}\ref{npn1}) that the   
non-principality implies $\rc$ being {\ubf not} a name 
of the real $x_{\xi}[\uG]$. 
And further, the avoidance condition in the next definition 
will be shown to imply $\rc$ being a name 
of a non-generic real.

\bdf
\lam{avod}
Let $\jpi,\jqo$ be \muf s, $\jpi\bssq\jqo$, 
$\xi\in\abc\jpi$; 
\rit{$\jqo$ avoids a\/ \qn\ $\rc$ over\/ $\jpi$ 
at\/ $\xi$}, 
in symbol\/ $\jpi\ssd\rc\xi\jqo$, 
\kmar{ssd rc xi }
\index{avoids}%
\index{zzzpissqcqo@$\jpi\ssd\rc\xi \jqo$}%
\index{multiforcing!refinementc@refinement, 
$\jpi\ssd\rc\xi \jqo$}%
\index{refinementc@refinement, $\jpi\ssd\rc\xi \jqo$}%
if for each 
$Q\in\jqo(\xi)$, $\jqo$ seals the set
$$
\dqu\rc Q\jpi
\kmar{dqu rc Q jpi}
=
\ens{\zr\in\mt\jpi} 
{\xi\in\abs\zr\land
\zr\,\text{\rm\ directly forces }\, \rc\nin[Q]}\,, 
$$ 
over $\jpi$ in the sense of Definition~\ref{ssl} --- 
that is formally $\jpi\ssb{\dqu\rc Q\jpi}\jqo$.
\edf


\ble
\lam{pqo}
Assume that\/ $\jpi,\jqo,\jsg$ are \muf s, 
$\xi\in\abc\jpi$, and\/ 
$\rc$ is a\/ \rn\jpi.
Then$:$
\ben
\renu
\itlb{pqo0}%
if\/ $\jpi\ssd\rc\xi \jqo$ and\/ $Q\in\jqo(\xi)$ then\/ 
the set\/ $\dqu\rc Q\jpi$ is open dense in\/ $\mt\jpi$ 
and pre-dense in\/ $\mt\jpq\;;$ 

\itlb{pqo1}%
if\/ $\jpi\ssd\rc\xi \jqo\bssq\jsg$ then\/ 
$\jpi\ssd\rc\xi \jsg\;;$ 

\itlb{pqo2}%
if\/ $\sis{\jpi_\al}{\al<\la}$  
is a\/ \dd\bssq increasing sequence in\/ \mfp, 
$0<\mu<\la$, $\jpi=\bkw_{\al<\mu}\jpi_\al$, 
and\/ $\jpi\ssd\rc\xi \jpi_\mu$, then\/
$\jpi\ssd\rc\xi \jqo=\bkw_{\mu\le\al<\la}\jpi_\al$. 
\een                                             
\ele
\bpf
\ref{pqo0} 
Apply Lemma~\ref{pqn}\ref{pqn0}.
To prove \ref{pqo1} 
let 
$S\in\jsg(\xi)$. 
Then, as $\jqo\bssq\jsg$, there is a finite set 
$\ans{Q_1,\dots,Q_m}\sq \jqo(\xi)$ such that 
$S\sq Q_1\cup\dots\cup Q_m$.
We have 
$\jpi\ssb{\dqu\rc {Q_i}\jpi}\jqo$ for all $i$ 
as $\jpi\ssd\rc\xi \jqo$, thus  
$\jpi\ssb{\dqu\rc {Q_i}\jpi}\jsg$, $\kaz i$, 
by Lemma~\ref{pqn}\ref{pqn1}.
Note that $\bigcap_{i}\dqu\rc {Q_i}\jpi\sq \dqu\rc {S}\jpi$ 
since $S\sq \bigcup_iQ_i$. 
We conclude that $\jpi\ssb{\dqu\rc{S}\jpi}\jsg$ 
by Lemma~\ref{pqn}\ref{pqn0b}.
Therefore $\jpi\ssd\rc\xi \jsg$, as required. 

To prove \ref{pqo2} make use of Lemma~\ref{pqn}\ref{pqn2} 
the same way.
\kmar{rassmotret'}
\epf

\parf{Generic refinement avoids non-principal names} 
\las{emf2}

The following theorem says that generic refinements as 
in Section~\ref{jex} avoid nonprincipal names. 
It resembles Theorem~\ref{uu4} to some extent, yet the 
latter is not directly applicable here as both 
the \mut\ $Q$ and  
the set $\dqu\rc Q\jpi$ depend on $\jqo$, and hence the 
sets $\dqu\rc Q\jpi$ do not necessarily belong to $\cmp$. 
However the proof will be based on rather similar
arguments. 

\bte
\lam{K}
Under the assumptions of Theorem~\ref{dj}, if\/ 
$\et\in\abc\jpi\sq\cM$ and\/
$\rc\in\cM$ is a\/ \rn\jpi\ non-principal 
over\/ $\jpi$ at\/ $\et$ then\/ $\jpi\ssd\rc\et \jqo$.  
\ete   
\bpf
Assume that $\jqo=\jqo[\dphi]$ is obtained from an 
\dd{\cmp}generic sequence $\dphi$ in $\om\ti\sys\jpi$,
as in Definition~\ref{dPhi}. 
We stick to the notation of \ref{dPhi}.

Let 
$Q\in\jqo(\et)=\djqo\et$; 
we have to prove that $\jqo$ seals the set
$\dqu\rc Q\jpi$ over\/ $\jpi$. 
By construction $Q=\qfj\et K\ret s$ 
for some $K<\om$ and $s\in\qfj\et K$; 
it can be assumed that simply $Q=\qfj\et K$. 
Following the proof of Theorem~\ref{uu4}, 
we suppose that $\zp\in\mt\jpi$, $\ju\in\mt\jqo$, 
$\abc\ju\cap\abc\zp=\pu$,  
and 
$\zc\ju\xi=\qfj\xi{,k_{\xi}}$, 
for each $\xi\in\abs\ju$. 
We have to find a \mut\ $\zq$ which witnesses 
\ref{ssl}\ref{ssl*} for $\ju,\,\zp,\,\zD=\dqu\rc Q\jpi$.
Note that $\et$ may or may not belong to the set 
$\abs\ju$, and even if $\et\in\abs{\ju}$, 
so $k_\et$ is defined,  then $K$ 
may or may not be equal to $k_\et$.
In the remainder of the proof, we use 
{\ubf the notation of Definition~\ref{uu4d}}, 
in particular, $\syn\jpi n$, 
$\bsn n\vpi$, $\jvn \bs\vpi$.

Assume that $\zr\in\mt\jpi$, $\abs\zr\cap\abs\ju=\pu$.
Consider the set $\Da_{\zr}\in\cM$ of all pairs   
$\ang{n,\vpi}\in \om\ti\sys\jpi$, 
such that $\vpi\in\syn\jpi n$ (see Def.~\ref{uu4d}), 
$\ang{\et,K}\in \abc\vpi$, 
and there is 
a \mut\ $\zq\in \mt\jpi$ satisfying $\zq\leq\zr$, 
still $\abs \ju\cap\abs \zq=\pu$, and 
\ben
\nenup
\itlb{ku41}%
if $\bs\in\bsn n\vpi$ and $t\in\zd\vpi{\et}K\cap2^n$ then 
$\jvn\bs\vpi\cup\zq$ 
directly forces $\rc\nin[\zd\vpi{\et}K\ret t]$.
\een
Condition \ref{ku41} is similar to (1)  
of Lemma~\ref{uu4a}, of course. 
Note that direct forcing of $\rc\nin[Q]$ 
cannot be used in \ref{ku41} 
because $Q$ is not necessarily an element of $\cM$, 
but $\rc\nin[\zd\vpi{\et}K]$ will be an effective 
replacement.

\ble
\lam{Ka}
If\/ $\zr\in\mt\jpi$, $\abs\zr\cap\abs\ju=\pu$, 
then\/ 
$\Da_{\zr}$ is dense in\/ $\om\ti\sys\jpi$.
\ele
\bpf
We follow the proof of  Lemma~\ref{uu4a}.
Let $\ang{n,\psi}\in\om\ti\sys\jpi$. 
We wlog assume that $\psi\in\syn\jpi n$ 
(see Lemma~\ref{uu4a}), 
so $\ang{\xi,k_{\xi}}\in\abs\psi$ for all $\xi\in\abs\ju$ 
and    
$\zd\psi\xi{k}\ret t\in\dP_\xi$ 
for all $\ang{\xi,k}\in\abs{\psi}$ 
and $t\in2^n\cap\zd\psi\xi{k}$,
and $\ang{\et,K}\in\abs\psi$ as well.

We have to define a system $\vpi\in\syn\jpi n$ 
such that $\ang{n,\vpi}\lec\ang{n,\pi}$ and 
$\vpi\in\Da_r$. 
As in the proof of Lemma~\ref{uu4a}, it suffices to 
fulfill \ref{ku41} for {\ubf one particular pair} of 
$\bs=\sis{s_\xi}{\xi\in\abs\ju}\in
\bsn{n}{\psi}$ 
and $t\in\zd\psi{\et}K\cap2^n$; the final goal is then 
achieved by simple iteration through all such pairs. 
We have two cases.\vom

{\ubf Case 1:} 
$\et\in\abs\ju$, $K=k_\et$, $t=s_\et$.  
Consider the \mut\ $\jvn{\bs}{\psi}\in\mt\jpi$.
The set $\ddj \et\rc\jpi$, as in Definition~\ref{npri}, 
is dense by the non-principality of $\rc$. 
It follows that there are \mut s $\zq,\jv\in\mt\jpi$ 
such that $\abs{\jv}=\abs\ju$, $\jv\leq\jvn{\bs}{\psi}$, 
$\abs{\zq}\cap\abs\ju=\pu$, $\zq\leq\zr$, and 
$\jv\cup\zq\in \ddj \et\rc\jpi$. 
Therefore $\jv\cup\zq$ 
directly forces $\rc\nin[\zc{\zq}{\et}]$.
Define a system  
$\vpi\in \sys\jpi$ with $\abs{\vpi}=\abs\psi$, 
from $\psi$ by: 
\ben
\aenu
\itlb{ka1}%
shrinking each tree 
$\zd\psi\xi{,k_{\xi}}\ret{s_\xi}$ 
($\xi\in\abs\ju$)
to $\zc\jv\xi$, so that 
$\zd{\vpi}\xi{,k_{\xi}}\ret{s_\xi}=\zc\jv\xi$, 

\itlb{ka2}%
in particular, shrinking 
$\zd\psi\et{K}\ret{t}$ 
to $\zc\jv\et$, so that 
$\zd{\vpi}\et{K}\ret{t}=\zc\jv\et$, 
\een
and no other changes. 
We have $\ang{n,\vpi}\lec\ang{n,\psi}$, 
$\jvn{\bs}{\vpi}=\jv$, and 
$\zd{\vpi}\et{K}\ret{t}=\zc\jv\et$ 
by construction.
In particular, $\jvn{\bs}{\vpi}\cup\zq$ 
directly forces $\rc\nin[\zd{\vpi}\et{K}\ret{t}]$, 
thus \ref{ku41} holds.\vom

{\ubf Case 2:} 
not Case 1. 
By Lemma~\ref{avo}, there exist \mut s 
$\zq,\jv\in\mt\jpi$ and a tree $T\in\dP_\et$ 
such that $T\sq \zd\psi\et{K}\ret{t}$, 
$\abs{\jv}=\abs\ju$, $\jv\leq\jvn{\bs}{\psi}$, 
$\abs{\zq}\cap\abs\ju=\pu$, $\zq\leq\zr$, and 
$\jv\cup\zq$ 
directly forces $\rc\nin[T]$.
Define a system  
$\vpi\in \sys\jpi$ with $\abs{\vpi}=\abs\psi$, 
that extends $\psi$ by \ref{ka1} above and: 
\ben
\aenu
\atc\atc
\itlb{ka3}%
shrinking $\zd\psi\et{K}\ret{t}$ 
to $T$, so that 
$\zd{\vpi}\et{K}\ret{t}=T$, 
\een
and no other changes. 
Note that \ref{ka1} and \ref{ka3} do not contradict 
each other since $\ang{\et,T,t}\ne\ang{\xi,k_\xi,s_\xi}$ 
for all $\xi\in\ju$ by the Case 2 hypothesis.
We have $\ang{n,\vpi}\lec\ang{n,\psi}$, 
$\jvn{\bs}{\vpi}=\jv$, and 
$\zd{\vpi}\et{K}\ret{t}=\zc\jv\et$ 
by construction.
In particular, $\jvn{\bs}{\vpi}\cup\zq$ 
directly forces $\rc\nin[\zd{\vpi}\et{K}\ret{t}]$, 
thus \ref{ku41} holds.   
\epF{Lemma}

Come back to the theorem.
As $\Da_{\zp}\in{\cmp}$, we have 
$\ang{n_j,\vpi_j}\in\Da_{\zp}$ for some $j$
by the lemma. 
Let this be witnessed by a   
\mut\ $\zq\in\mt\jpi$, so that 
$\zq\leq\zp$, $\abs \ju\cap\abs \zq=\pu$, and 
\ref{ku41} holds for $n=n_j$, $\vpi=\vpi_j$. 
In particular, as 
$\zd{\vpi_j}{\et}K=\bigcup_{t\in\zd{\vpi_j}{\et}K\cap2^n}
\zd{\vpi_j}{\et}K\ret t$, 
the \mut\ $\jvn\bs{\vpi_j}\cup\zq$ 
directly forces $\rc\nin[\zd{\vpi_j}{\et}K]$ 
whenever $\bs\in\bsn n{\vpi_j}$, hence 
directly forces $\rc\nin[Q]$ as well, because 
$Q=\qfj\et K\sq \zd{\vpi_j}{\et}K$ by construction.
Thus if $\bs\in\bsn n{\vpi_j}$ then 
$\jvn\bs{\vpi_j}\cup\zq\in \dqu\rc Q\jpi$, 
and hence 
$\jvn\bs{\vpi_j}\in\duz{\dqu\rc Q{\vpi_j}}\ju\zq$.
On the other hand, 
$[\ju]\sq\bigcup_{\bs\in\bsn{n}{\vpi_j}}[\jvn{\bs}{{\vpi_j}}]$, 
so that  
$\ju\sqf\bigvee\duz{\dqu\rc Q\jpi}\ju\zq$, 
as required. 
\epf

\parf{Consequences for generic extensions} 
\las{crn}

We first prove a lemma on adequately representation 
of reals in 
\dd{\mt\jpi}generic extensions    
by real names. 
Then Theorem~\ref{npn} will show corollaries  
for non-principal names.
 
\ble
\lam{r2n}
Suppose that\/ $\jpi$ is a \mre\ \muf\ and 
$G\sq\mt\jpi$ is 
generic over the ground set universe\/ $\rV$.    

If\/ $x\in\rV[G]\cap\dn$  
then 
there is a\/ \rn\jpi\ $\rc\in\rV$, 
$\rc\sq\mt\jpi\ti\om\ti2$,  
such that\/ $x=\rc[G]$. 

If\/ $\mt\jpi$ is a CCC forcing\snos
{The CCC property means that   
\index{ccc@CCC}%
\index{forcing!ccc@CCC}%
every antichain\/ $A\sq\mt\jpi$ is at most countable.}
in\/ $\rV$, and\/ $\rc\in\rV$, 
$\rc\sq\mt\jpi\ti\om\ti2$ is a\/ \rn\jpi,
then there is a\/ {\ubf small} 
\rn\jpi\ $\rd\in\rV$, 
$\rd\sq\mt\jpi\ti\om\ti2$, such that\/  
$\mt\jpi$ forces\/
$\rc[\uG]=\rd[\uG]$ over\/ $\rV$.
\ele
\bpf
The first claim is an instance of a general
forcing theorem.  
To prove the second one, 
extend each set $\kic n\sq\mt\jpi$  
to an open dense set 
$\kipc n\jpi=\ens{\zp\in\mt\jpi}
{\sus\zq\in\kkc n{}\,(\zp\leq\zq)}$,
choose maximal antichains $A_n\sq \kipc n\jpi$
in those sets
--- which are countable by CCC, and then let
$A_{ni}=
\ens{\zp\in A_n}{\sus\zq\in\kkc n{i}\,(\zp\leq\zq)}$ 
and 
$\rd=\ens{\ang{\zp,n,i}}{\zp\in A_{ni}}$.
\epf

\bte
\lam{npn}
Let\/ $\jpi$ be a \mre\ \muf\ and\/ $\xi\in\abc\jpi$. 
Then
\ben
\renu
\itlb{npn1}%
if\/ $\mt\jpi$ is CCC, 
a set\/ $G\sq\mt\jpi$ is generic over the ground set 
universe\/ $\rV$,   
and\/ $x\in\rV[G]\cap\dn,$
then\/ $x\ne x_{\xi}[G]$ if and only if 
there is a small\/ \rn\jpi\
$\rc\sq\mt\jpi\ti(\om\ti2)$, 
non-principal over\/ $\jpi$ at\/ $\xi$  
and such that\/ $x=\rc[G]$$;$ 

\itlb{npn2}%
if\/ $\rc\sq\mt\jpi\ti(\om\ti2)$ is a \/ \rn\jpi, 
$\jqo$ is a \muf, $\jpi\ssd\rc\xi\jqo$, 
and a set\/ $G\sq\mt{\jpi\kw\jqo}$ is generic 
over\/ $\rV$ then\/ 
$\rc[G]\nin
\bigcup_{Q\in\jqo(\xi)}[Q]$.
\een
\ete
\bpf
\ref{npn1}
Let $x\ne x_{\xi}[G]$.
By a known forcing theorem, there is a 
\rn\jpi\ $\rc$ such that $x=\rc[G]$ and 
$\mt\jpi$ forces that $\rc\ne x_{\xi}[\uG]$, and, 
by Lemma~\ref{r2n}, $\rc$ is small 
since $\mt\jpi$ is CCC.
It remains to show that 
$\rc$ is a non-principal name over $\jpi$ at $\xi$,  
that is, the set
$$
\ddj\xi\rc\jpi
=\ens{\zp\in\mt\jpi}
{\xi\in\abs\zp\land\zp\,\text{ directly forces }\,
\rc\nin[\zc\zp\xi ]}\,.
$$
is open dense in\/ $\mt\jpi$. 
The openness is clear, let us prove the density. 
Consider any $\zq\in\mt\jpi$. 
Then $\zq$ \dd{\mt\jpi}forces $\rc\ne x_{\xi }[\uG]$ by the 
choice of $\rc$, hence we can assume that, for some $n$, 
$\rc(n)\ne x_{\xi}[\uG](n)$ 
is \dd{\mt\jpi}forced by $\zq$.
Then by Lemma~\ref{avo} there is a \mut\ 
$\zp\in\mt\jpi$, $\zp\leq\zq$, and $s\in \om^{n+1},$ 
such that $\zp$ directly forces $s\sq \rc$. 
Now it suffices to show that $s\nin \zc\zp\xi$.
Suppose otherwise:  $s\in \zc\zp\xi$. 
Then the tree $T=\req{\zc\zp\xi}s$ still belongs to $\mt\jpi$. 
Therefore the \mut\ $\zr$ defined by 
$\zc\zr\xi=T$ and 
$\zc\zr{\xi'}=\zc\zp{\xi'}$ for each  
$\xi'\ne\xi$, belongs to $\mt\jpi$ and 
satisfies $\zr\leq\zp\leq\zq$.
However $\zr$ directly forces both 
$\rc(n)$ and $x_{\xi}[\uG](n)$ 
to be equal to one and the same value $\ell=s(n)$, which 
contradicts to the choice of $n$.

To prove the converse let $\rc\sq\mt\jpi\ti(\om\ti2)$ 
be a \rn\jpi\ non-principal over $\jpi$ at $\xi$, 
and $x=\rc[G]$.
Assume to the contrary that 
$x=x_{\xi}[G]$. 
There is a \mut\ $\zq\in G$ which \dd{\mt\jpi}forces 
$\rc=x_{\xi}[\uG]$.
As $\rc$ is non-principal, 
there is a \mut\ $\zp\in G\cap \ddj\xi \rc\jpi$, 
$\zp\leq\zq$. 
Thus $\zp$ directly forces $\rc\nin[\zc\zp\xi]$, and hence 
\dd{\mt\jpi}forces the same statement.
Yet $\zp$ \dd{\mt\jpi}forces 
$x_{\xi}[\uG]\in[\zc\zp\xi]$, 
of course, and this is a contradiction.\vom

\ref{npn2}
Suppose towards the contrary that 
$Q\in\jqo(\xi)$ and 
$\rc[G]\in[Q]$.
By definition, $\jqo$ seals, over\/ $\jpi$, the set
$$
\dqu\rc Q\jpi
=
\ens{\zr\in\mt\jpi} 
{\xi\in\abs\zr\land
\zr\,\text{\rm\ directly forces }\, \rc\nin[Q]}\,. 
$$
Therefore $\dqu\rc Q\jpi$ is pre-dense in 
$\mt{\jpi\kw\jqo}$ by Lemma~\ref{pqn}, and hence  
$G\cap \dqu\rc Q\jpi\ne \pu$. 
In other words, there is a \mut\ $\zr\in\mt\jpi$  
which directly forces $\rc\nin[Q]$. 
It easily follows that $\rc[G]\nin[Q]$, 
which is a contradiction.
\epf

\parf{Combining refinement types} 
\las{comb}  

Here we summarize the properties 
of generic refinements considered above. 
The next definition combines the 
refinement types 
$\ssa{D}\yi\ssb {\zD}\yi\ssb {\rc}\yi\ssd\rc\xi$.

\index{zzzpissqMqo@$\jpi\ssm\cM\jqo$}%
\index{multiforcing!refinementM@refinement, $\jpi\ssm\cM\jqo$}%
\index{refinementM@refinement, $\jpi\ssm \cM\jqo$}%

\bdf
\lam{extM}
Suppose that $\jpi\bssq\jqo$ are \muf s 
and $\cM\in\hc$ is any set. 
Let $\jpi\ssm\cM\jqo$ 
\kmar{ssm mm}
mean that the four following requirements hold:
\ben
\nenu
\itlb{extM1}%
if $\xi\in\abc\jpi$, $D\in\cM$, $D\sq\jpi(\xi)$, 
$D$ is pre-dense in 
$\jpi(\xi)$, then 
$\jpi(\xi)\ssa D\jqo(\xi)$;

\itlb{extM2}%
if ${\zD}\in\cM$, ${\zD}\sq\mt\jpi$, 
${\zD}$ is open dense in $\mt\jpi$,
then $\jpi\ssb {\zD}\jqo$; 

\itlb{extM2c}%
if $\rc\in\cM$ is a \rn\jpi\  
then $\jpi\ssb {\rc}\jqo$; 

\itlb{extM3}%
if $\xi\in\abc\jpi$ and $\rc\in\cM$ is a 
\rn\jpi, \rit{non-principal over\/ $\jpi$ at\/ $\xi$}, 
then $\jpi\ssd\rc\xi\jqo$.\qed 
\een
\eDf

\bcor
[of lemmas \ref{pqm}, \ref{pqn}, \ref{pqs}, \ref{pqo}]
\lam{mpq}
Let\/ $\jpi,\jqo,\jsg$ be \muf s and\/ 
$\cM$ be a countable set.
Then$:$
\ben
\renu
\itlb{mpq1}%
if\/ $\jpi\ssm\cM\jqo\bssq\jsg$ then\/ $\jpi\ssm\cM\jsg\;;$ 

\itlb{mpq2}%
if\/ $\sis{\jpi_\al}{\al<\la}$  
is a\/ \dd\bssq increasing sequence in\/ \mfp, 
$0<\mu<\la$, $\jpi=\bkw_{\al<\mu}\jpi_\al$, 
and\/ $\jpi\ssm{\cM}\jpi_\mu$, then\/
$\jpi\ssm{\cM}\jqo=\bkw_{\mu\le\al<\la}\jpi_\al$.\qed 
\een                                             
\ecor

\bcor
\lam{geec} 
If\/ $\jpi$ is a small \muf, $\cM\in\hc$, and\/ 
$\jqo$ is an\/ \dd\cM generic refinement of\/ $\jpi$ 
{\rm(exists by Lemma~\ref{gee}!)}, 
then\/ 
$\jpi\ssm\cM\jqo$.
\ecor
\bpf
We have $\jpi\ssm\cM\jqo$ by a combination of 
\ref{dj}\ref{dj3}, \ref{uu4}, \ref{uu4c}, and 
\ref{K}.
\epf

\parf{Increasing sequences of \muf s} 
\las{incS}

Recall that $\mfp$ is the collection of all \muf s 
(Section~\ref{mul}). 
\imar{mfp}
Let 
\imar{mf}%
$$
\bay{rclcccc}
\mf&=&\ens{\jpi\in\mfp} 
{\jpi\ \text{is a \rit{small} \muf}};\\[1ex] 
\spmf&=&\ens{\jpi\in\mfp} 
{\jpi\ \text{is a small and \rit{\oso} \muf}}. 
\eay
$$
\imar{spmf}%
\index{multiforcing!sMF@$\mf$}%
\index{zzsMF@$\mf$}%
\index{multiforcing!spMF@$\spmf$}%
\index{zzspMF@$\spmf$}%
Thus a \muf\ $\jpi\in\mfp$ belongs to $\mf$ if 
$\abc\jpi\sq\omi$ 
is (at most) countable and if $\xi\in \abc\jpi$ then 
$\jpi(\xi)$ is a countable forcing in $\ptf$, 
and $\jpi\in\spmf$ requires that in addition each 
$\jpi(\xi)$ is \oso\ 
(Definition~\ref{ptfn}).

\bdf
\lam{sdf}
If $\ka\le\omi$ then let $\vmf_\ka$
be the set of all 
\index{zzMF-@$\vmf$}%
\index{zzMF-w@$\vmi$}%
\index{length!$\len\vjpi$}%
\index{zzlenpi@$\len\vjpi$}%
\dd\bssq increasing sequences 
$\vjpi=\sis{\nor\jpi\al}{\al<\ka}$ of 
\imar{vmf vmi}%
\muf s $\nor\jpi\al\in\spmf$, 
\rit{domain-continuous}  
in the sense that if $\la<\ka$ is a limit ordinal then 
$\abs{\nor\jpi\la}=\bigcup_{\al<\la}\abs{\nor\jpi\al}$.
Let $\vmf=\bigcup_{\ka<\omi}\vmf_\ka$.

We order $\vmf\cup\vmi$ by the usual 
relations $\sq$ and $\su$ of extension of sequences. 
Thus 
$\vjpi\su\vjqo$ iff $\ka=\len\vjpi<\la=\len\vjqo$ and 
$\nor\jpi\al=\nor\jqo\al$ for all $\al<\ka$.
In this case, if $\cM$ is any set, and $\nor\jqo\ka$ 
(the first term of $\vjqo$ absent in $\vjpi$) 
satisfies $\jpi\ssm\cM\nor\jqo\ka$, where 
$\jpi=\bkw_{\al<\ka}\nor\jpi\al$, 
then we write $\vjpi\su_\cM\vjqo$. 
\index{multiforcing!extension, $\vjpi\su\vjqo$}%
\index{zzzpisuqo-@$\vjpi\su\vjqo$}%
\index{multiforcing!Mextension@\dd\cM extension, 
$\vjpi\su_\cM\vjqo$}%
\index{zzzpisuMqo-@$\vjpi\su_\cM\vjqo$}%
\index{zzsuM@$\su_\cM$}%

If $\vjpi\in\vmf_\ka$ then let 
$\mt\vjpi= \mt\jpi$, where 
\index{multitree!$\mt\vjpi$}%
\index{zzMTpi-@$\mt\vjpi$}%
$\jpi=\bkw\vjpi=\bkw_{\al<\ka}\nor\jpi\al$ 
(componentwise union). 
Accordingly, a \rit{\rn\vjpi} means 
\index{real name!pi-realname@\rn\vjpi}%
a \rn\jpi. 
\edf

\ble
\lam{142}
If\/ $\vjpi,\vjqo\in\vmf$,
$\rc\in\mm$ is a\/ \rn\vjpi, and\/  
$\vjpi\su_{\ans\rc}\vjqo$,
then\/ $\rc$ is a\/ \rn\vjqo. 
\ele
\bpf
Let $\ka=\len\vjpi<\la=\len\vjqo$ and 
$\jpi=\bkw\vjpi=\bkw_{\al<\ka}\nor\jpi\al$. 
Then by definition $\jpi\ssm{\ans\rc}\nor\jqo\ka$, 
hence $\jpi\ssb\rc\nor\jqo\ka$ because 
$\rc$ is a \rn\jpi. 
However  
$\jpi\ssb\rc\jqo=\bkw_{\ka\le\al<\la}\nor\jqo\al$
by Lemma~\ref{pqs}\ref{pqs3}. 
Therefore $\rc$ is a \pol{(\jpi\kw\jqo)} name
by Lemma~\ref{pqs}\ref{pqs1}. 
However, 
$\jpi\kw\jqo=\bkw_{\al<\la}\nor\jqo\al=\bkw\vjqo$.
\epf

\bdf
\lam{zflm} 
Let $\zflm$ be the subtheory of $\ZFC$ 
\index{zzZFLM@$\zflm$}%
\index{theory!zzZFLM@$\zflm$}%
including all axioms except for the power set axiom, 
plus the axiom of constructibility $\rV=\rL$, 
and plus the axiom saying that 
$\cP(\om)$ exists. 
(Then $\omi$, $\hc$, and generally sets related to the 
continuum, like $\dn\yd\pet$, exist, too.)
The axiom of choice is included in $\zflm$ in the form of 
the wellorderability principle. 

If $x\in\hc$ 
($\hc$= hereditarily countable sets, Footnote~\ref{hc}) 
then let
\index{zzLx@$\fl(x)$}%
\index{model!zzLx@$\fl(x)$}%
$\fl(x)$ be the least \rit{countable} transitive model 
(CTM) of $\zflm$ containing $x$ and satisfying 
\index{model!CTM, countable transitive model}%
\index{CTM, countable transitive model}%
$x\in(\hc)^{\fl(x)}$. 
It necessarily has the form $\fl(x)=\rL_\mu$ for 
some $\mu=\mu_x<\omi$.

An ordinal $\xi<\ka$ is a \rit{crucial ordinal} 
\index{crucial ordinal}%
\index{ordinal!crucial ordinal}%
of a sequence 
$\vjpi=\sis{\nor\jpi\al}{\al<\ka}\in\vmf_\ka$ if 
$(\bkw_{\al<\xi}{\nor\jpi\al})
\ssm{\fl(\vjpi\res\xi)}\jpi_\xi$ 
holds.
This is equivalent to 
$\vjpi\res\xi\su_{\fl{(\vjpi\res\xi)}}\vjpi$.
\edf 

\ble
\lam{xisc}
Suppose that\/ $\ka\le\omi$ and\/  
$\vjpi=\sis{\nor\jpi\al}{\al<\ka}\in\vmf_\ka$. 
Then$:$
\ben
\renu
\itlb{xisc0}%
$\jpi=\bkw\vjpi=\bkw_{\al<\ka}\nor\jpi\al$ is 
a \mre\ \muf$;$ 

\itlb{xisc1}%
if\/ $\ka<\la\le\omi$ and\/ 
$\cM\in\hc$ 
then 
there is a \muq\/ $\vjqo\in\vmf$ satisfying\/ 
$\len\vjqo=\la$ and\/ $\vjpi\su_\cM\vjqo\;;$

\vyk{
\itlb{xisc1+}%
if\/ $\sis{s_\al}{\al<\la}$ is
a\/ \dd\su increasing sequence of countable sets\/
$s_\al\sq\omi$,
$s_\al=\abc{\nor\jpi\al}$ for all\/ $\al<\ka$,
and\/ $s_\ga=\bigcup_{\al<\ga}s_\al$ for all
limit\/ $\ga<\la$, then
there is a \muq\/ $\vjqo\in\vmf$ satisfying\/ 
$\len\vjqo=\la$,
$\abc{\nor\jqo\al}=s_\al$ for all\/ $\al<\la$,
and\/ $\vjpi\su_\cM\vjqo\;;$

\itlb{xisc1*}%
if\/ $\vjpi,\vjro,\vjqo\in\vmf$
and\/ $\vjpi\su_\cM\vjro\sq\vjqo$ then\/ 
$\vjpi\su_\cM\vjqo\;;$ 
}

\itlb{xisc2}%
if\/ $\xi<\ka$ is a crucial ordinal of\/ $\vjpi$, 
$\jpi_{<\xi}=\bkw_{\al<\xi}\nor\jpi\al$, 
$\jpi_{\ge\xi}=\bkw_{\xi\le\ba<\ka}\nor\jpi\ba$,
then\/ 
$\jpi_{<\xi}\ssm{\fl{(\vjpi\res\xi)}}\jpi_{\ge\xi}$  
and 
$\jpi_{<\xi}
\ssm{\fl{(\vjpi\res\xi)}}\nor\jpi\ba$ 
for\/ $\xi\le\ba<\ka$, 
hence
\ben
\aenri
\itlb{xisc2x}\msur%
$\mt{\jpi_{\ge\xi}}$ is open dense in\/ $\mt\vjpi\,,$

\vyk{
\itlb{xisc2a}%
if\/ $\xi<\ka$, $D\in\mm$, $D\sq\jpi(\xi)$, 
$D$ is pre-dense in\/ $\jpi(\xi)$, then\/ 
$D$ remains pre-dense in\/ $\jpi(\xi)\cup\jqo(\xi)\,,$ 
}

\itlb{xisc2b}%
if ${\zD}\in\fl{(\vjpi\res\xi)}$, 
${\zD}\sq\mt{\vjpi\res\xi}$, 
${\zD}$ is open dense in\/ $\mt{\vjpi\res\xi}$,
then ${\zD}$ is pre-dense in 
$\mt{\jpi_{<\xi}\kw\jpi_{\ge\xi}}=\mt\vjpi\,.$ 
%
\een
\een
\ele

\bpf
\ref{xisc0} 
Make use of Lemma~\ref{pqr}\ref{pqr3}.

\ref{xisc1} 
We define terms $\nor\jqo\al$ of the \muq\ $\vjqo$ required
by induction.

Naturally put $\nor\jqo\al=\nor\jpi\al$ for each $\al<\ka$.
To define the crucial term $\nor\jqo\ka$, 
we wlog assume that $\cM$ contains $\vjpi$ and satisfies 
$\ka\sq\cM$ (otherwise take a bigger set). 
By Lemma~\ref{gee}, there is an \dd\cM generic 
refinement $\jpi'$ of $\jpi=\bkw_{\al<\ka}\nor\jpi\al$. 
By Theorem~\ref{dj}, $\jpi'$ is a small \mdi\ \muf,  
$\jpi\bssq\jpi'$, and $\nor\jpi\al\bssq\jpi'$ 
for all $\al<\ka$.
In addition $\jpi\ssm\cM\jpi'$ by Corollary~\ref{geec}.
We let $\nor\jqo\ka=\jpi'$. 
The extended \muq\  
$\vjqo_+=\sis{\nor\jqo\al}{\al<\ka+1}$ 
belongs to $\vmf_{\ka+1}$ and satisfies 
$\vjpi\su_\cM\vjqo_+$.

The following steps are pretty similar, except that 
we can take $\cM=\pu$.

To prove the main claim of \ref{xisc2} make use of
Corollary \ref{mpq}. 


To prove \ref{xisc2x} apply Corollary~\ref{pqrC}. 

\ref{xisc2b}
As $\jpi_{<\xi}\ssm{\fl{(\vjpi\res\xi)}}\jpi_{\ge\xi}$ and 
$\zD\in\fl{(\vjpi\res\xi)}$, 
we have $\jpi_{<\xi}\ssb \zD\jpi_{\ge\xi}$, 
therefore $\zD$ is pre-dense in $\mt\vjpi$ by 
Lemma~\ref{pqn}\ref{pqn0}.
%
\epf

\vyk{
Our plan regarding the forcing notion for Theorem~\ref{Tsep} 
will be to define a certain \muq\ $\vjPi$ in $\vmi$ and 
the ensuing \muf\ $\jPi=\bkw\vjPi$ with 
remarkable properties related to definability and its own 
genericity of some sort. 
}

\parf{The key sequence} 
\las{sek4}  

In this section we define the forcing notion to prove 
Theorem~\ref{Tsep}. 
It will have the form $\mt\jPi$, for a certain \muf\ 
$\jPi$ with $\abc\jPi=\omi$. 
The \muf\ $\jPi$ will be equal to the componentwise 
union of terms of a certain sequence $\vjPi\in\vmi$.
The construction of this sequence 
in $\rL$, the constructble universe, will employ 
some ideas related to diamond-style constructions, 
as well as to some sort of 
\rit{definable genericity}. 
The following definition introduces another 
important notion 
involved in the construction. 

\bdf
\lam{blo}
A \muq\ $\vjpi\in\vmf$ 
\rit{\se s\/} a set $W$ if either  
\index{block}%
\index{block!positive}%
\index{block!negative}%
$\vjpi\in W$ (\rit{positive} \sen) 
or there is no \muq\  
$\vjqo\in W$ 
extending $\vjpi$ (\rit{negative} \sen).   
\edf

Recall that $\hc$ = all  
hereditarily countable sets, Footnote~\ref{hc}.

\bdf
\lam{HCdef}
We use standard notation 
$\is\hc n$, $\ip\hc n$, $\id\hc n$ 
(slanted $\is{}{},\ip{}{},\id{}{}$) 
\index{zzzsig@$\is\hc n$, $\ip\hc n$, $\id\hc n$}%
\index{definability classes!
lightface $\is\hc n$, $\ip\hc n$, $\id\hc n$}%
for classes of \rit{lightface} definability in $\hc$  
(no parameters allowed), and 
$\hs n$, $\hp n$, $\hd n$ 
\index{definability classes!
boldface $\fs\hc n$, $\fp\hc n$, $\fd\hc n$}%
\index{zzzsigb@$\hs n$, $\hp n$, $\hd n$}%
for \rit{boldface} definability in $\hc$ 
(parameters in $\hc$ allowed).
It is well-known that if $n\ge1$ and $X\sq\dn$ then 
\bce
$X\in\is\hc n\leqv X\in\is1{n+1}$\,, 
\quad and\quad 
$X\in\hs n\leqv X\in\fs1{n+1}$,
\ece 
and the same for 
$\ip{}{}$, $\fp{}{}$, $\id{}{}$, $\fd{}{}$.  
\edf



\bte
[in $\rL$]
\lam{kyt}
Let\/ $\nn\ge3$.
There exists a \muq\/ 
$\vjPi=\sis{\jpn\al}{\al<\omi}\in\vmi$ satisfying 
the following requirements$:$ \ 
\ben
\renu
\itlb{kyt1} 
the \muq\/ 
$\vjPi$ belongs to the definability class\/ $\id\hc{\nn-2}\;;$ 

\itlb{kyt0} 
$\abs{\bkw\vjPi}=\omi\;;$

\itlb{kyt2} 
if\/ $\nn\ge4$ and\/ $W\sq\vmf$ is a 
boldface\/ $\hs{\nn-3}$ set  
then there is an ordinal\/ $\ga<\omi$ such that the 
\muq\/ 
$\vjPi\res\ga$ \se s\/ $W\;;$ 

\itlb{kyt3}
there is a club\snos
{\label{club}%
Closed unbounded set = club.} 
\index{set!closed unbounded, club}%
\index{club (closed unbounded)}
such that every\/ $\ga\in\dC$ 
is a crucial ordinal for\/ $\vjPi$.
\een 
\ete
\bpf
We argue under $\rV=\rL$. 
Let $\ufo_\nn(p,x)$ be a canonical universal 
\index{universal formula, $\ufo_m(p,x)$}%
$\is{}{\nn-3}$ formula, so that the family of all 
boldface $\hs{\nn-3}$ sets $X\sq\hc$ 
is equal to  the family of all sets 
of the form 
$\ups_\nn(p)=\ens{x\in\hc}{\hc\models \ufo_\nn(p,x)}$, 
\index{sets $\ups_\nn(p)$}%
\index{zzzYpsnz@$\ups_\nn(p)$}%
$p\in\hc$.

For $\al<\omi$, 
define a sequence $\vjpi[\al]\in\vmf$ by induction 
as follows. 

We let $\vjpi[0]=\pu$, the empty sequence.\vhm  

{\ubf Step $\al\to\al+1$.} 
Suppose that 
$\vjpi[\al]\in\vmf$ is defined, $\ka=\dom{\vjpi[\al]}$, 
$\cM=
\fl{(\vjpi[\al])}$,  
and $p_\al$ is the $\al$-th element of $\hc=\rL_{\omi}$ 
in the sense of the Goedel wellordering $\lel$. 
\index{zzzz<L@$\lel$, Goedel wellordering}%
\index{Goedel wellordering $\lel$}%
By Lemma~\ref{xisc}\ref{xisc1}, there is a 
sequence $\vjta\in\vmf_{\ka+1}$ satisfying 
$\vjpi[\al]\su_\cM\vjta$. 
By Corollary~\ref{xistt}, there is a 
sequence $\vjqo\in\vmf_{\ka+2}$ satisfying 
$\vjta\su\vjqo$ and $\al\in\abs{\vjqo(\ka+1)}$. 
Finally if $\nn\ge4$ then 
there is a sequence $\vjpi\in\vmf$ satisfying   
$\vjqo\su\vjpi$ and \se ing the set $\ups_\nn(p_\al)$.  
Let $\vjpi[\al+1]$ be the $\lel$-least of such 
sequences $\vjpi$.\vom
 
\vyk{
To get a sequence $\vjqo\in W$ extending a 
given $\vjpi\in\vmf_\ka$, $\ka<\omi$, we 
first extend $\vjpi$ by Lemma~\ref{xisc}\ref{xisc1} 
to a sequence $\vjqo\in\vmf_{\ka+2}$.
Its last term $\vjqo(\ka+1)=\nor\jqo{\ka+1}$, 
a \muf,  
satisfies $\nor\jqo{\ka}\bssq\nor\jqo{\ka+1}$.  
If still $\xi\nin\abs{\nor\jqo{\ka+1}}$ then add 
$\xi$ to $\abs{\nor\jqo{\ka+1}}$ and put 
$\nor\jqo{\ka+1}(\xi)=\dpo$ 
(the Cohen forcing, see Definition~\ref{cloL}). 
The final sequence $\vjqo$ belongs to $W$ 
and $\vjpi\su\vjqo$.
}%

{\ubf Limit step.} 
If $\la<\omi$ is limit then we naturally 
define $\vjpi[\la]=\bigcup_{\al<\la}\vjpi[\al]$.\vom  

We have $\al<\ba\imp \vjpi[\al]\su\vjpi[\ba]$ 
by construction, hence 
$\vjPi=\bigcup_\al\vjpi[\al]\in\vmf_{\omi}$.  
To prove \ref{kyt1}, 
note first of all that the relation 
$R(\vjpi,\vjqo,\cM):=\text{``$\vjpi\su_\cM\vjqo$''}$ 
is absolute for all 
transitive models of $\zflm$, hence $R$ is $\id\hc{1}$. 
Easily the assignment $\vjpi\mto\fl{(\vjpi)}$ 
is $\id\hc{1}$ as well. 
Finally ``to \se\ $\ups_\nn(p)$'' is a $\id\hc{\nn-2}$
relation. 
Using these facts, it's a routine estimation to 
verify  \ref{kyt1}.

To check \ref{kyt0}, 
note that $\al\in\abs{\bkw\vjpi[\al+1]}$
by construction.

To check \ref{kyt2} ($\nn\ge4$), note that any   
boldface\/ $\hs{\nn-3}$ set $W\sq\vmf$   
is equal to $\ups_\nn(p_\al)$ for some $\al<\omi$, 
so $\ga=\dom{\vjpi[\al+1]}$ is as required.  

\ref{kyt3} 
The set $\dC=\ens{\dom{\vjpi[\al]}}{\al<\omi}$ 
is a club by the limit step of the 
construction. 
Moreover if $\ga=\dom{\vjpi[\al]}\in \dC$ then 
$\vjPi\res\ga=\vjpi[\al]$, and $\ga$ is crucial for 
$\vjPi$ by construction.
\epf

\bgg
[in $\rL$]
\lam{keys}
From now on we fix {\ubf a number $\nn\ge3$}  
\index{zzn@$\nn$}%
as in Theorem~\ref{Tsep}.
We also  
{\ubf fix a \muq}
$\vjPi=\sis{\jpn\al}{\al<\omi}\in\vmi$ satisfying 
\index{sequence!key \muq\ $\vjPi$}%
\index{zzzPivec@$\vjPi$}%
\index{zzzPixi@$\jPi(\xi)$}%
\ref{kyt1} -- \ref{kyt3}  
of Theorem~\ref{kyt} for this $\nn$. 
We call this fixed $\vjPi$ 
\rit{the key \muq}. 
\egg

\ble
\lam{18a}
If\/ $\nn\ge4$ and\/ $W\sq\vmf$ is a\/ $\hs{\nn-3}$ set  
dense in\/ $\vmf$ then there is an ordinal\/ 
$\ga<\omi$ such that\/ $\vjPi\res\ga\in W$. 
\ele
\bpf
By \ref{keys}, $\vjPi$ satisfies 
\ref{kyt2} of Theorem~\ref{kyt}, hence 
there is an ordinal $\ga<\omi$ such that 
$\vjPi\res\ga$ \se s $W$. 
The negative \sen\ is impossible by the density 
of $W$, hence in fact $\vjPi\res\ga\in W$. 
\epf

\parf{Key forcing notion} 
\las{bpf}

We continue to argue in $\rL$, and we'll  make use of 
the key \muq\ $\vjPi=\sis{\jPi_\al}{\al<\omi}$  
introduced by \ref{keys}. 

\bdf
[in $\rL$]
\lam{bang}
Define the \muf s 
$$
\bay{lclcrlllll}
\kmar{jPi}%
\jPi
&=&\bkw_{\al<\omi}\jPi_\al
&\in&\mfp,& \\[1ex] 
\index{key elements!Pi@$\jPi$}%
\index{multiforcing!Pi@$\jPi$}%
\index{zzzPi@$\jPi$}%
\index{key elements!Pi<ga@$\pilg\ga$}%
\index{multiforcing!Pi<ga@$\pilg\ga$}%
\index{zzzPi<ga@$\pilg\ga$}%
\kmar{pilg ga}%
\pilg\ga
&=&\bkw_{\al<\ga}\jPi_\al&\in& \mf , & 
\text{for each $\ga<\omi$},\\[1ex]
\index{key elements!Pi>ga@$\pigg\ga$}%
\index{multiforcing!Pi>ga@$\pigg\ga$}%
\index{zzzPi>ga@$\pigg\ga$}%
\kmar{pigg ga}%
\pigg\ga
&=&\bkw_{\ga\le\al<\omi}\jPi_\al&\in& \mfp, & 
\text{for each $\ga<\omi$}.
\eay%
$$%
We further define $\fP=\mt{\jPi}=\mt{\vjPi}$, 
\kmar{fP}%
\index{key elements!Pd@$\fP$}%
\index{multiforcing!Pd@$\fP$}%
\index{zzPd@$\fP$}%
and, for all $\ga<\omi$,
$$
\fpl\ga=\mt{\pilg\ga}=\mt{\vjPi\res\ga}\,,\quad 
\kmar{fpl ga fpg ga}%
\index{key elements!Pd<ga@$\fpl\ga$}%
\index{multiforcing!Pd<ga@$\fpl\ga$}%
\index{zzPd<ga@$\fpl\ga$}%
\fpg\ga=\mt{\pigg\ga}=\mt{\vjPi\res(\omi\bez\ga)}\,.
\eqno\qed
$$
\eDf
\index{key elements!Pd>ga@$\fpg\ga$}%
\index{multiforcing!Pd>ga@$\fpg\ga$}%
\index{zzPd>ga@$\fpg\ga$}%
The set $\fP=\mt{\jPi}$ will be our 
{\ubf key forcing notion}. 
\index{forcing!key forcing}%
\index{key forcing}%

\bcor
[in $\rL$, by \ref{kyt}\ref{kyt0}]
\lam{doml}
$\jPi$ is a \mre\ \muf\ and\/ 
$\abc\jPi=\omi$, thus\/ 
$\fP=\prod_{\xi<\omi}\jPi(\xi)$ 
(with finite support).\qed  
\ecor

\vyk{
\bpf
To prove $\abc\jPi=\omi$, note that
the set $W$ of all \muq s $\vjpi\in\vmf$ satisfying 
$\xi\in\abc{\bkw\vjpi}$ is $\hs1$ 
(with $\xi$ as a parameter of definition). 
In addition $W$ is dense in $\vmf$.\snos 
{To get a sequence $\vjqo\in W$ extending a 
given $\vjpi\in\vmf_\ka$, $\ka<\omi$, we 
first extend $\vjpi$ by Lemma~\ref{xisc}\ref{xisc1} 
to a sequence $\vjqo\in\vmf_{\ka+2}$.
Its last term $\vjqo(\ka+1)=\nor\jqo{\ka+1}$, 
a \muf,  
satisfies $\nor\jqo{\ka}\bssq\nor\jqo{\ka+1}$.  
If still $\xi\nin\abs{\nor\jqo{\ka+1}}$ then add 
$\xi$ to $\abs{\nor\jqo{\ka+1}}$ and put 
$\nor\jqo{\ka+1}(\xi)=\dpo$ 
(the Cohen forcing, see Definition~\ref{cloL}). 
The final sequence $\vjqo$ belongs to $W$ 
and $\vjpi\su\vjqo$.} 
Therefore by Lemma~\ref{18a}  
there is an ordinal $\ga<\omi$ such that 
$\vjPi\res\ga\in W$, as required.
\epf
}

If $\xi<\omi$ then, following the corollary,
let $\al(\xi)<\omi$
be the least ordinal $\al$ satisfying  
$\xi\in\abc{\jpn{\al}}$. 
Thus a forcing $\jpn{\al}(\xi)\in\ptf$
is defined whenever $\al$ satisfies 
$\al(\xi)\le\al<\omi$, and    
$\sis{\jpn{\al}(\xi)}{\al(\xi)\le\al<\omi}$ is a 
\dd\bssq increasing sequence of special forcings 
in\/ $\ptf$.
Note that
$\jPi(\xi)= \bigcup_{\al(\xi)\le\al<\omi} \jpn{\al}(\xi)$
by construction.

\bcor
[in $\rL$]
\lam{doml2}
The sequence of ordinals\/ 
$\sis{\al(\xi)}{\xi<\omi}$ and the sequence of
forcings\/ 
$\sis{\jpn{\al}(\xi)}{\xi<\omi,\,\al(\xi)\le\al<\omi}$
are\/ $\id\hc{\nn-2}$. 
\ecor
\bpf
By construction the following double equivalence holds:
$$
\bay{cclc}
\al<\al(\xi)
&\eqv&
\sus\jpi(\jpi=\jpn\al\land \xi\in\dom\jpi)
&\eqv 
\\[1ex]
&\eqv&
\kaz\jpi(\jpi=\jpn\al\imp \xi\in\dom\jpi)
&. \eay
$$
However $\jpi=\jpn\al$ is a $\id\hc{\nn-2}$ relation by 
Theorem~\ref{kyt}\ref{kyt1}.
It follows that so is the sequence $\sis{\al(\xi)}{\xi<\omi}$.
The second claim is similar.
\epf

\bcor
[in $\rL$, of Lemma~\ref{pqr}\ref{pqr3}]
\lam{xiden}
If\/ $\xi<\omi$ and\/ $\al(\xi)\le\al<\omi$ 
then the set\/ $\jpn\al(\xi)$ 
is pre-dense in\/ $\jPi(\xi)$ and in\/ $\jPi$.\qed  
\ecor

In spite of Corollary~\ref{doml}, the sets 
$\abc{\pilg\ga}$ can be 
quite arbitrary (countable) subsets of $\omi$. 
However we get the next corollary:

\bcor
[in $\rL$, of Corollary~\ref{doml}]
\lam{domc}
$\dcp=\ens{\ga<\omi}{\abc{\pilg\ga}=\ga}$ 
\index{key elements!C'@$\dcp$}%
\index{set!C'@$\dcp$}%
\index{zzCd'@$\dcp$}%
\imar{dcp}%
is a club in\/ $\omi$.\qed
\ecor

To prove the CCC property, 
we'll need the following result.

\ble
[in $\rL$]
\lam{lclub}
If\/ $X\sq\hc=\rL_{\omi}$ then the set\/ $\skri O_X$ 
of all ordinals\/ $\ga<\omi$, such that\/ 
$\stk{\rL_\ga}{X\cap\rL_\ga}$ is an elementary submodel 
of\/ $\stk{\rL_{\omi}}{X}$ and\/ 
$X\cap\rL_\ga\in\fl(\vjPi\res\ga)$, 
is stationary, hence unbounded in\/ $\omi$.

More generally, if\/ $X_n\sq\hc$ for all\/ $n$ 
then the set\/ $\skri O$ of all 
ordinals\/ $\ga<\omi$, such that\/ 
$\stk{\rL_\ga}{\sis{X_n\cap\rL_\ga}{n<\om}}$ 
is an elementary submodel of\/  
$\stk{\rL_{\omi}}{\sis{X_n}{n<\om}}$ 
and\/ $\sis{X_n\cap\rL_\ga}{n<\om}\in\fl(\vjPi\res\ga)$, 
is stationary, hence unbounded in\/ $\omi$.
\ele
\bpf
Let $C\sq\omi$ be a club. 
Let $M$ be a countable elementary submodel of 
$\rL_{\om_2}$ 
containing $C\yi\omi\yi X\yi\vjPi$, 
and such that $M\cap\rL_{\omi}$ is transitive. 
Let $\phi:M\onto\rL_\la$ be the Mostowski collapse, 
and $\ga=\phi(\omi)$. 
Then 
$$
\ga<\la<\omi,\quad 
\phi(X)=X\cap\rL_\ga,\quad 
\phi(C)=C\cap\ga,\quad 
\phi(\vjPi)=\vjPi\res\ga
$$
by the choice of $M$. 
It follows that 
$\stk{\rL_\ga}{X\cap\rL_\ga,C\cap\ga,\vjPi\res\ga}$ 
is an elementary submodel of 
$\stk{\rL_{\omi}}{X,C,\vjPi}$, so $\ga\in\skri O_X$.
Moreover, $\ga$ is uncountable in $\rL_\la$, hence 
$\rL_\la\sq\fl(\vjPi\res\ga)$. 
(See Definition \ref{zflm} on models 
$\fl(\vjpi)\models\zflm$.)
We conclude that   
$X\cap\rL_\ga\in\fl(\vjPi\res\ga)$ since
$X\cap\rL_\ga\in\rL_\la$ by construction.
On the other hand, $C\cap\ga$ is unbounded in $\ga$ 
by the elementarity, therefore $\ga\in C$, as required. 

The second, more general claim does not differ much.
\epf

\bcor
[in $\rL$]
\lam{ccc}
The forcing\/ $\fP$ satisfies CCC.
Therefore\/ 
\dd\fP generic extensions of\/ $\rL$ preserve 
cardinals.
\ecor
\bpf
Suppose that $A\sq\fP=\mt{\vjPi}$ is a maximal antichain. 
By \ref{keys} and Theorem~\ref{kyt}\ref{kyt3}, 
there is a club $\dC\sq\omi$ such that every 
$\ga\in\dC$ is a crucial ordinal for $\vjPi$.
By Lemma~\ref{lclub}, 
there is an ordinal $\ga\in\dC$ such that 
$A'=A\cap\fpl\ga$ is a maximal antichain in 
$\fpl\ga=\mt{\vjPi\res\ga}$ and 
$A'\in\fl(\vjPi\res\ga)$. 
It follows that the set 
$D(A')=\ens{\zp\in\fpl\ga}{\sus\zq\in A\,(p\leq q)}
\in\fl(\vjPi\res\ga)$ 
is open dense in $\fpl\ga$. 

Yet $\ga$ is a crucial ordinal for $\vjPi$, 
therefore by Lemma~\ref{xisc}\ref{xisc2b} 
both the set $D(A')$, 
and hence $A'$ itself as well, remain  
pre-dense in the whole set $\fP=\mt\vjPi$. 
We conclude that $A=A'$ is countable. 
\epf

\bcor
[in $\rL$]
\lam{Cjden}
If a set\/ $D\sq\fP$ is pre-dense in\/ $\fP$ 
then there is an ordinal\/ 
$\ga<\omi$ such that\/ $D\cap\fpl\ga$ is already pre-dense 
in\/ $\fP$.
\ecor
\bpf
We can assume that in fact $D$ is dense.
Let $A\sq D$ be a maximal antichain in $D$; then $A$ is 
a maximal antichain in $\fP$ because of the density of $D$.
Then $A\sq \fpl\ga$ for some $\ga<\omi$ by Lemma~\ref{ccc}.
But $A$ is pre-dense in\/ $\fP$.
\epf

\vyk{
We finish this Section with the following corollary 
that transforms the genericity
of the sequence $\vjPi$ in $\vmf$ by \ref{kyt}\ref{kyt2} 
into a similar genericity of the set $\fP=\mt\vjPi\sq\md$
inside the set $\md$ of all \mut s.
Diven a set $X\sq\md$, we let 
$X^-=\ens{\zp\in\md}{\neg\:\sus\zq\in X\,(\zq\leq\zp)}$, 
the set of all \mut s not extendable to \mut s in X. 
Then $X\cup X^-$ is obviously dense in $\md$. 

\bcor
[in $\rL$]
\lam{kyt2c}
Let\/ $X\sq\md$ be\/ $\hs{\nn-3}$, $\nn\ge4$. 
Then\/ $(X\cup X^-)\cap\fP$ is dense in\/ $\fP$. 
Thus if\/ $X$ is already dense in\/ $\md$ 
then\/ $X\cap\fP$ is dense in\/ $\fP$. 
\ecor
\bpf
Assume that $\zp_0\in\fP$; we have to find a \mut\ 
$\zq\in\fP\cap(X\cup X^-)$, $\zp\leq\zp_0$.
By construction there is $\ga_0<\omi$ such that 
$\zp_0\in\fP_{<\ga_0}=\mt{\vjPi\res\ga_0}$. 
The set $W$ of all sequences $\jpi\in\vmf$, satisfying: 
$\sus \zq\in\mt{\jpi}(\zq\leq\zp_0\land\zq\in X)$,
is $\hs{\nn-3}$ along with $X$. 
Recall that $\vjPi$ satisfies \ref{kyt2} of 
Theorem~\ref{kyt} by \ref{keys}.
Therefore there is an ordinal $\ga<\omi$ 
such that $\vjPi\res\ga$ \se s $W$. 
We can \noo\ assume that $\ga_0\le\ga$. 

If $\vjPi\res\ga\in W$, and this is witnessed by 
$\zq\in\mt{\vjPi\res\ga}$, then $\zq$ is as required.

Thus we can assume that $\vjPi\res\ga$ \se s $W$ 
{\ubf negatively}, 
so there is no sequence $\vjqo\in W$ extending 
$\vjPi\res\ga$. 
Recall that 
$\zp_0\in\mt{\vjPi\res\ga_0}\sq\mt{\vjPi\res\ga}$.
As $\vjPi$ is \dd\ssq increasing, there is a \mut\ 
$\zq\in\mt{\jPi(\ga)}$, $\zq\leq\zp_0$. 
We claim that $\zq\in X^-$.

Assume to the contrary that there is a \mut\ 
$\zr\in X$, $\zr\leq\zq$. 
Then $\zr\leq\zp_0$.
Corollary~\ref{pqrC} (the last claim) gives a 
small special \muf\ $\jsg$ such that $\zr\in\mt\jsg$ 
and any \muf\ $\jpi$ satisfying $\jpi\ssq\jPi(\ga)$ 
also satisfies $\jpi\ssq\jsg$.  
It follows that the extension $\jqo$ of $\vjPi\res\ga$, 
of length $\dom\jqo=\ga+1$, by the last term 
$\jqo(\ga)=\jsg$, is still a sequence in $\vmf$. 
Note that $\zr\in X\cap\mt\jqo$ by construction, hence 
$\jqo\in W$.
But this contradicts to the negative \sen\ assumption. 

Thus $\zq\in X^-$, as required. 
($\zq\leq\zp_0$ and $\zq\in\fP$ hold by construction.)
\epf
}

\parf{Basic generic extension} 
\las{kge}

Recall that the key \muq\ 
$\vjPi=\sis{\jPi_\al}{\al<\omi}$ 
of small \oso\ \muf s $\jPi_\al$ 
is defined in $\rL$ by \ref{keys},   
the componentwise union 
$\jPi=\bkw_{\al<\omi}\jPi_\al$ is a \muf, 
$\abc\jPi=\omi$ in $\rL$, 
and $\fP=\mt\vjPi=\mt\jPi\in\rL$ is our key 
forcing notion, equal to the finite-support product 
$\prod_{\xi<\omi}\jPi(\xi)$ of arboreal  
forcings $\jPi(\xi)$ in $\rL$. 
See Section~\ref{bpf}, where some properties of $\fP$ are 
established, including CCC and definability of the 
factors $\jPi(\xi)$ in $\rL$.
Our goal will be to show that certain submodels of 
\dd\fP generic models prove Theorem~\ref{Tsep}.

\bre
\lam{omil}
From now on, we'll typically 
argue in $\rL$ and in \dd\omil preserving 
generic extensions $\rL$
(this includes, \eg, \dd\fP generic extensions 
by Corollary~\ref{ccc}). 
Thus it will always be the case that 
$\omil=\omi$.
This allows us to still think that $\abc\jPi=\omi$ 
(rather than $\omil$).
\ere


\bdf
\lam{bstr}
Let a set $G\sq\fP$ be generic over the 
constructible set universe $\rL$.
If $\xi<\omi$ 
then following Remark~\ref{adds}, we  
\bit
\item[$-$]
define 
$G(\xi)=\ens{\zc \zp\xi}{\zp\in G\land\xi\in\abc\zp}
\sq  \jPi(\xi)$;

\item[$-$]
let $x_{\xi}=x_{\xi}[G]\in\dn$ be the only real in 
$\bigcap_{T\in G(\xi)}[T]$.

\item[$-$]
\index{zzXXG@$\bX=\bX[G]$}%
let 
$\bX=\bX[G]=\sis{x_{\xi}[G]}{\xi<\omi}
=\ens{\ang{\xi,x_{\xi}[G]}}
{\xi<\omi}$..
\eit
Thus $\fP$ adjoins an array  
$\bX[G]$ 
of reals to $\rL$, 
where each   
\index{zzxxxikG@$x_{\xi}=x_{\xi}[G]$}%
$x_{\xi}=x_{\xi}[G]\in\dn\cap\rL[G]$ 
is a \dd{\jPi(\xi)}generic real over $\rL$, 
and $\rL[G]=\rL[\bX[G]]$.

If $\Da\sq\omi$ then let 
$
\fP\res\Da=\mt{\jPi\res\Da}=\ens{\zp\in \fP}{\abs\zp\sq\Da}
$.
\edf

The next lemma makes use of the 
product structure of  $\fP$.

\ble
\lam{prods}
Suppose that\/ $\Da\in\rL$, $\Da\sq\omi$. 
Then\/ $\fP=\mt\jPi$ is equal to the 
product\/ $(\fP\res\Da)\ti(\fP\res\Da')$, where\/ 
$\Da'=\omi\bez\Da$. 
If\/ $G\sq\fP$ is generic over $\rL$, then the set\/ 
$ 
G\res\Da=\ens{\zp\in G}{\abs\zp\sq\Da}
$  
is \dd{(\fP\res\Da)}generic over\/ $\rL$. 

If\/ $\xi<\omi$, $\xi\nin\Da$, then\/ 
$x_{\xi}[G]\nin \rL[G\res\Da]$.\qed
\ele

\parf{Definability of generic reals} 
\las{dgr}

Recall that the factors $\jPi(\xi)$ of the forcing 
notion $\fP=\mt\jPi=\prod_{\xi<\omi}\jPi(\xi)$ are defined by 
$\jPi(\xi)=\bigcup_{\al(\xi)\le\al<\omi}\jpn{\al}(\xi)$, 
where $\al(\xi)<\omi$, 
and the sets $\jpn{\al}(\xi)$ are 
countable sets of perfect trees, whose  definability in 
$\rL$ is determined by Corollary~\ref{doml2}.
We'll freely use the notation introduced by 
Definition~\ref{bstr}.

\bte
\lam{xik}
Assume that a set\/ $G\sq\fP$ is\/ 
\dd\fP generic over\/ $\rL$, 
$\xi<\omi$, and\/ $x\in\rL[G]\cap\dn.$ 
The following are equivalent$:$
\ben
\nenu
\itlb{xik1}
$x=x_{\xi}[G]\;;$
\hfill 
$(2)$
$x$ is\/ \dd{\jPi(\xi)}generic over\/ $\rL\,;
\qquad\qquad$

\atc
\itlb{xik3}
$x\in\bigcap_{\al(\xi)\le\al<\omi}
\bigcup_{T\in\jPi_\al(\xi)}[T]$.
\een
\ete
\bpf
$\ref{xik1}\imp(2)$ is a routine (see Remark~\ref{adds}). 
To check $(2)\imp\ref{xik3}$ recall that each set 
$\jPi_\al(\xi)$ is pre-dense in $\jPi(\xi)$ by 
Lemma~\ref{pqr}\ref{pqr3}.
It remains to establish $\ref{xik3}\imp\ref{xik1}$.
Suppose 
that $x\in\rL[G]\cap\dn$ but \ref{xik1} fails, 
that is, $x\ne x_{\xi}[G]$. 
By Theorem~\ref{npn}\ref{npn1} there is a 
small ($\fP=\mt\jPi$ is CCC by \ref{ccc}) 
\rn\jPi\ $\rc\in\rL$, 
such that $\rc\sq\fP\ti\om\ti2$,  $x=\rc[G]$, 
and $\rc$ is non-principal over $\jPi$ at $\xi$, 
meaning that the set
$$
\ddj\xi\rc\jPi
=\ens{\zp\in\fP=\mt\jPi}
{\xi\in\abs\zp\land\zp\,\text{ directly forces }\,
\rc\nin[\zc\zp\xi]}\,.
$$
is open dense in\/ $\fP=\mt\jPi$.
By the smallness of $\rc$, 
there is an ordinal $\ga<\omi$ 
such that $\rc$ is a \rn{\pilg\ga}, and 
we can assume, by Corollary~\ref{Cjden},  that 
$\ddi\xi k\rc\jPi\cap\fpl{\ga}$ is 
pre-dense in $\fP$, therefore, open dense in 
$\fpl{\ga}$ --- and then $\rc$ is 
non-principal over $\pilg\ga$ at $\xi$. 
We can further assume that 
$\rc\in\fl{(\vjPi\res\ga)}$. 
And finally, we can assume that $\ga$ belongs to the 
set $\dC$ of Theorem~\ref{kyt}\ref{kyt3}, 
in other words, $\ga$ is crucial for $\vjpi$, 
that is, $\pilg \ga \ssm{\fl(\vjPi\res\ga)} \jPi_\ga$.
It follows that 
$\pilg \ga \ssm{\fl(\vjPi\res\ga)} \pigg\ga$
by Lemma~\ref{xisc}\ref{xisc2}.
Then $\pilg \ga \ssd\rc\xi \pigg\ga$ holds
as well by \ref{extM}\ref{extM3}, 
since $\rc\in\fl(\vjPi\res\ga)$ and 
because of the non-principality of $\rc$.
Now Theorem~\ref{npn}\ref{npn2} with 
$\jpi=\pilg\ga$ and $\jqo=\pigg\ga$ 
(note that  $\jpi\kw\jqo=\jPi$) 
implies 
$x=\rc[G]\nin\bigcup_{Q\in\pigg\ga(\xi)}[Q]$, 
in particular, 
$x\nin\bigcup_{Q\in\jPi_\ga(\xi)}[Q]$. 
In other words, \ref{xik3} fails as well.
\epf  

\bcor
\lam{mod2}
Assume that\/ $G\sq\fP$ is\/ \dd\fP generic 
over\/ $\rL$, and\/ $M$ is a generic extension 
of\/ $\rL$ satisfying\/ $\dn\cap M\sq\rL[G]$. 
Then\/ $\gfu[G]\cap M$
is a set of definability class\/ $\ip\hc{\nn-2}$ in\/ $M$. 
\ecor
\bpf
By the theorem, it holds in $M$ 
that $\ang{\xi,x}\in \gfu[G]$ iff 
$$
\kaz \al<\omi\:\sus T\in\jPi_\al(\xi)
\big(\al(\xi)\le\al \imp x\in[T]\big),
$$
which can be re-written as 
$$
\kaz \al<\omi\:\kaz \mu<\omi\:\kaz Y\:
\sus T\in Y
\big(\mu=\al(\xi)\land Y=\jPi_\al(\xi)\land\mu\le\al
\imp x\in[T]\big).
$$
Here the equality $\mu=\al(\xi)$ 
is $\id\hc{\nn-2}$ by Corollary~\ref{doml2}, and so is 
the equality $Y=\jPi_\al(\xi)$ by 
Corollary~\ref{doml2}. 
It follows that the whole relation is $\ip\hc{\nn-2}$, 
since the quantifier $\sus T\in Y$ is bounded.
\epf


\bcor
\lam{mod3}
If\/ $G\sq\fP$ is\/ 
\dd\fP generic over\/ $\rL$  
then 
it holds in\/ $\rL[G]$  that there is 
a ``good''\/ $\id1{\nn}$ wellordering of\/ 
$\dn$ of length\/ $\omi$. 
\ecor
\bpf
If $\ga<\omi$ then let 
$\bX_{\ga}=\sis{x_{\xi}[G]}{\xi<\ga}$. 
The equality $Y=\bX_{\ga}$ is a $\ip\hc{\nn-2}$ relation 
in $\rL[G]$ (with $\ga,Y$ as arguments)
by Corollary~\ref{mod2}. 
If $x\in\dn\cap\rL[G]$ then let $\ga(x)$ be 
the least $\ga<\omi$ such that $x\in\rL[\bX_{\ga}]$, 
and $\nu(x)<\omi$ be the index of $x$ in the 
canonical wellordering of $\dn$ in $\rL[\bX_{\ga}]$.
We wellorder $\dn\cap\rL[G]$ according to 
the lexicographical ordering of the triples 
$\ang{\tmax\ans{\ga(x),\nu(x)},\ga(x),\nu(x)}$. 
This is $\id\hc{\nn-1}$ by the above, hence $\id1{\nn}$.
The ``goodness'' (that is, the set of all coded proper  
initial segments has to be $\is1\nn$) 
can be easily verified. 
\epf

\parf{The non-separation model} 
\las{ns}

The model for Theorem~\ref{Tsep}  
will be defined on the base of a 
\dd\fP generic extension $\rL[G]$ of $\rL$. 
More exactly, it will have the form  $\rL[G\res\Da]$, 
where $\Da\sq\omil$ will itself be a generic set 
over $\rL[G]$.

Let $\fQ={\idt}^{\omil}\cap\rL$ with 
\rit{countable} support; 
a typical element of $\fQ$ is a partial map 
$q\in\rL$ from $\omil$ to the 3-element set 
$\idt$, with 
a domain $\dom q\sq\omil$ countable in $\rL$, 
that is, just bounded in $\omil$. 
(The choice of the 3-element set $\idt$ is explained 
by later considerations, see Definition \ref{ggh}.)
We order $\fQ$ opposite to extension, that is, let 
$q\leq q'$ ($q$ is stronger) iff $q'\sq q$. 
Thus $\fQ\in\rL$, and, inside $\rL$, $\fQ$ is   
equal to the product $\idt^{\omi}$ with countable support.
Accordingly a \dd\fQ generic object is a full 
\dd\fQ generic map 
$H:\omil\to\idt$. 

Recall that $\fP$ is a CCC forcing 
in $\rL$ by Corollary~\ref{ccc}.

\ble
\lam{prek}
$\fP$ remains CCC in 
any\/ \dd\fQ generic extension\/ $\rL[H]$ of\/ $\rL$,
therefore\/ $\fP\ti\fQ$ preserves cardinals over\/ $\rL$.
\ele

\bpf
Suppose towards the contrary that some $q'\in\fQ$ forces 
that $C$ is an uncountable 
antichain in $\fP$, $C$ being a \dd\fQ name. 
Note that, in $\rL$, $\fQ$ is \rit{countably complete\/}:  
if $q_0\geq q_1\geq q_2\geq \dots$ is a sequence 
in $\fQ$ then there is a condition 
$q=\bigcup_kq_k\in\fQ$; $q\leq q_k\yd\kaz k$.
Therefore, \rit{arguing in\/ $\rL$}, we can define by induction 
a decreasing sequence $\sis{q_\xi}{\xi<\omi}$ in $\fQ$ and 
a sequence of pairwise incompatible conditions $p_\xi\in\fP$, 
such that $q_0\leq q'$ and each $q_\xi$ forces that 
$p_\xi\in C$. 
But then $A=\ens{p_\xi}{\xi<\omi}\in\rL$ 
is an uncountable antichain in $\fP$, a contradiction.
\epf

\ble
\lam{prer}
Assume that a set\/ $G\ti H$ is\/ 
$\fP\ti\fQ$-generic over\/ $\rL$. 
Then 
\ben
\renu
\itlb{prer1}%
$\dn\cap\rL[G,H]\sq\rL[G]$, 
hence\/ $\omil=\omi^{\rL[G]}=\omi^{\rL[G,H]}\;;$ 

\itlb{prer2}%
if\/ $\Da\in\rL$, $\Da\sq\omil$ 
then\/ $\rL[G\res\Da,H]\cap\dn\sq\rL[G\res\Da]\;;$

\itlb{prer3}%
if\/ $\Da\in\rL[H]$, $\Da\sq\omil$, and\/ 
$\xi<\omil$ then\/  
$x_{\xi}[G]\in\rL[G\res\Da]$ iff\/ $\xi\in\Da$. 
\een
\ele

\bpf
Note that $\fQ$ may not be countably complete 
in $\rL[G]$ any more, 
so that the most elementary way to prove \ref{prer1}  
does not work.
However consider $\rL[G,H]$ as a \dd\fP generic extension 
$\rL[H][G]$ of $\rL[H]$. 
Let $x\in\dn\cap\rL[H][G]$. 
As $\fP=\mt\jPi$ is CCC in $\rL[H]$ by Lemma~\ref{prek}, 
there exists a small \rn\jPi\ $\rc\in\rL[H]$, such 
that $\rc\sq\fP\ti\om\ti2$ and $x=\rc[G]$. 
Because of the smallness, $\rc$ is effectively 
coded by a real, hence $\rc\in\rL$   
because $\rL[H]$ has just the same reals as $\rL$. 
Thus $\rc\in\rL$ and $x=\rc[G]\in\rL[G]$. 

The proof of \ref{prer2} is similar.

\ref{prer3}
In the nontrivial direction, suppose that $\xi\nin\Da$. 
Consider the set $\Da'=\omil\bez\ans{\xi}\in\rL$.
As obviously $G\res\Da\in\rL[G\res\Da',H]$, any real in 
$\rL[G\res\Da]$ belongs to $\rL[G\res\Da']$ by \ref{prer2}. 
But $x_{\xi}[G]\nin \rL[G\res\Da']$ by Lemma~\ref{prods}.
\epf

Recall that if $\nu\in\Ord$ then the ordinal 
product $2\nu$ is considered as the ordered sum of $\nu$ 
copies of $2=\ans{0,1}$. 
(Contrary to $\nu2=\nu+\nu$.)
Thus if $\nu=\la+m$, where $\la$ is a limit ordinal or $0$ 
and $m<\om$, then $2\nu=\la+2m$ and $2\nu+1=\la+2m+1$, 
and $\ang{\nu,i}\mto2\nu+i$ is a bijection of 
$\omi\ti2$ onto $\omi$.

\bdf
\lam{ggh}
If $H:\omil\to\idt$ then define sets 
$$
\bay{c}
\ah H=\ens{2\nu}{H(2\nu)=1}
\,,\quad   
\bh H=\ens{2\nu}{H(2\nu)=2} 
\,,\quad  
\ch H=\ens{2\nu}{H(2\nu)=3}\,,
\\[1.5ex]
\dh H=\ens{2\nu+1}{H(2\nu+1)=1}
\,,\quad   
\eh H=\ens{2\nu+1}{H(2\nu+1)=2} 
\,,
\eay
$$
and 
$\fh H=\ens{2\nu+1}{H(2\nu+1)=3}$, 
and further
\vspace*{2ex} 

$ 
\bay[t]{cclccccccc}
\Da_H
&=&
\ens{4\nu}{2\nu\in \ah H\cup \ch H}
\cup
\ens{4\nu+1}{2\nu\in \bh H\cup \ch H}
\cup\phantom{x}\\[1ex]
&&\hspace*{7ex}\phantom{x}\cup 
\ens{4\nu+2}{2\nu+1\in \dh H}
\cup
\ens{4\nu+3}{2\nu+1\in \eh H}\,.
\eay
$ 
\hfill
$
\bay[t]{ccccccc}
\,\\[1ex]
\qed
\eay
$\vhm 
\eDf

Note that $\mor H\sq\rL[G]$ is not necessarily 
true since the set $\Da_H$ 
does not necessarily belong to $\rL[G]$,  
but we have $\mor H\sq\rL[G][H]$, of course.

\parf{Non-separation theorem: the $\hc$ version} 
\las{nspf}

Now we prove the following result, the 
$\hc$-definability version of Theorem~\ref{Tsep}.

\bte
\lam{Tsep'}
Let a set\/ $G\sq\fP$ be\/ \dd\fP generic over\/ $\rL$ and\/ 
$H:\omil\to\idt$ be a map\/ \dd\fQ generic over\/ $\rL[G]$. 
Then it is true in\/ $\mor H$ that\/ 
\ben
\renu
\itlb{Tsep'1}
$\ah H\yd\bh H$ are disjoint\/ 
$\ip{\hc}{\nn-1}$ 
sets, not separable by disjoint\/ 
$\hs{\nn-1}$ sets$;$

\itlb{Tsep'2}
$\dh H\yd\eh H$ are disjoint\/ 
$\is{\hc}{\nn-1}$ 
sets, not separable by disjoint\/ 
$\hp{\nn-1}$ sets.
\een
\ete

The proof of Theorem~\ref{Tsep'} below in this Section 
includes a reference to the following result, which 
will have its own lengthy proof in the remainder.

\bte
[will be proved in Section \ref{ee}]
\lam{shon}
Assume that\/ $X\in\rL$, $X\sq\omil$ is unbounded in\/ 
$\omil$, and a set\/ $G\sq\fP$ is\/ 
\dd\fP generic over\/ $\rL$. 
Then\/ $\rL[G\res X]\cap\dn$ is an elementary submodel of\/ 
$\rL[G]\cap\dn$ \poo\ all\/ $\is1{\nn-1}$ formulas with real 
parameters in\/ $\rL[G\res X]$. 
\ete  

\bcor
\lam{shoC}
Under the assumptions of Theorem~\ref{Tsep'}, 
$\hc^{\mor H}$ 
is an elementary submodel of\/ 
$\hc^{\rL[G]}$ \poo\ all\/ $\is{}{\nn-2}$ formulas.
\ecor

Note that $\hc^{\mor H}\sq\hc^{\rL[G]}$
by Lemma~\ref{prer}, while $\mor H\not\sq\rL[G]$.

\bpf
[Corollary]
We have 
$\omil=\omi^{\rL[G]}=\omi^{\rL[G\res\Da_H]}$ and 
$\Da_H\cap\la\in\rL$ for any $\la<\omil$   
by Lemma~\ref{prer}. 
It remains to cite theorem~\ref{shon}, having in mind 
that $\is{\hc}{\nn-2}$-definability 
corresponds to $\is1{\nn-1}$-definability.\vom 

\epF{Corollary~\ref{shoC} from Theorem~\ref{shon}}

\bpf[Theorem~\ref{Tsep'}]
\ref{Tsep'1}
To check that, say, $\ah H$ is $\ip{\hc}{\nn-1}$ in $\mor H$, 
it suffices to prove that the equality
$$
\ah H=\ens{2\nu<\omi}
{\neg\:\sus x\,(\ang{4\nu+1,x}\in\gfu)}
$$
holds in $\mor H$, where 
$\gfu=\gfu[G]\cap\mor H$ is a $\ip\hc{\nn-2}$ set 
in $\mor H$ by  
Corollary~\ref{mod2}.
(For $\bh H$ it would be $\ang{4\nu,x}\in\gfu$ 
in the displayed formula.)

First suppose that $\nu<\omil$, $\xi=4\nu+1$, 
$x\in\mor H\cap\dn$, and $\ang{\xi,x}\in\gfu$; 
prove that $2\nu\nin \ah H$. 
Now, by definition $x=x_{\xi}[G]$, and $\xi\in\Da_H$ 
by Lemma~\ref{prer}\ref{prer3}. 
%
But then $2\nu\in \bh H\cup \ch H$, so $2\nu\nin \ah H$, 
as required. 

To prove the converse, let $2\nu\nin \ah H$, so that 
$2\nu\in \bh H\cup \ch H$. 
Then $\xi=4\nu+1\in\Da_H$, and hence $x=x_{\xi}\in\mor H$ 
and $\ang{\xi,x}=\ang{4\nu+1,x}\in \gfu$, as required.

To prove the non-separability claim, 
suppose towards the contrary that, in $\mor H$, 
the sets $\ah H,\bh H$ are separated by disjoint 
$\hs{\nn-1}$ sets 
$A,B\sq\omi=\omil$.
The sets $A,B$ are defined, in the set 
$\hc^{\mor H}$, by $\is{}{\nn-1}$ 
formulas, resp., $\vpi(a,\xi)\yi\psi(a,\xi)$, with a real 
parameter $a\in\mor H\cap\dn;$ hence, $a\in\rL[G]$ 
by Lemma~\ref{prer}. 
Let $\la<\omil$ be a limit ordinal such that 
$a\in\rL[G\res \Da_{H\la}]$, 
where $\Da_{H\la}=\Da_H\cap\la\in\rL$. 

If $K:\omil\to\idt$ 
(for instance, $K=H$), then let 
$$
\ahh{K}=\ens{\xi<\omil}
{\vpi(a,\xi)^{\hgh GK}}\,,\quad
\bhh{K}=\ens{\xi<\omil}
{\psi(a,\xi)^{\hgh GK}}\,.
\eqno(\ast)
$$
Then by definition $\ah H\sq A=\ahh H$, $\bh H\sq B=\bhh H$, 
and $\ahh H\cap\bhh H=\pu$.
Fix a condition $q_0\in\fQ$ compatible with $H$ 
(here meaning that simply $q_0\su H$), 
which forces the mentioned properties of $A,B$, 
so that, 
\ben
\fenu
\atc
\itla{fdag} 
if $K:\omil\to\idt$ is a map \dd\fQ generic over 
$\rL[G]$ and compatible with $q_0$, then $\ah K\sq \ahh K$,  
$\bh K\sq \bhh K$, and $\ahh K\cap\bhh K=\pu$.
\een
We may assume that $\dom{q_0}\sq\la$, 
otherwise just increase $\la$.

Let $\nu_0$ be any ordinal, $\la\le\nu_0<\omi$. 
Consider the maps $H_1\yi H_2\yi H_{3}:\omil\to\idt$, 
generic over $\rL[G]$, compatible with $q_0$, 
and satisfying ${H_i}(2\nu_0)=i$, $i=1,2,3$, and 
${H_1}(\al)={H_2}(\al)={H_{3}}(\al)$ 
for all $\al\ne 2\nu_0$.
Then $\Da_{H_{3}}=\Da_{H_1}\cup\ans{4\nu_0+1}$ 
by Definition~\ref{ggh}, hence,  
$\mor{H_1}\sq\mor{H_{3}}$. 
It follows by Corollary~\ref{shoC} 
that $\ahh{H_1}\sq\ahh{H_{3}}$. 
Therefore $\ah {H_1}\sq\ahh{H_1}\sq\ahh{H_{3}}$ 
by \ref{fdag}. 
We conclude that $2\nu_0\in \ahh{H_{3}}$, just because 
$2\nu_0\in \ah {H_1}$ by the choice of $H_1$.

And we have $2\nu_0\in \bhh{H_{3}}$ by a similar argument 
(with $H_2$). 
Thus $\ahh{H_{3}}\cap\bhh{H_{3}}\ne\pu$, 
contrary to \ref{fdag}.
The contradiction ends the proof of \ref{Tsep'1}.

The proof of \ref{Tsep'}\ref{Tsep'2} is pretty similar.\vom

\epF{Theorem \ref{Tsep'} modulo Theorem \ref{shon}}

\parf
[The main theorem modulo the elementary equivalence theorem] 
{The main theorem modulo theorem \ref{shon}} 
\las{meet}

\bpf[Theorem \ref{Tsep} modulo theorem \ref{shon}]
\ref{Tsep1}
We argue under the assumptions of Theorem~\ref{Tsep'}. 
To define a non-separable pair of $\ip1{\nn}$ sets 
in $\mor H$, 
let $\wo\sq\dn$ be the $\ip11$ set of codes of 
countable ordinals, and for $w\in\wo$ 
let $\abs w<\omi$ be the ordinal coded by $w$.
As $\omil=\omi$ by Corollary~\ref{ccc}, for any $\xi<\omi$ 
there is a code $w\in\wo\cap\rL$ with $\abs w=\xi$. 
Let $w_\xi$ be the \dd\lel least of those, and  
$X=\ens{w_\xi}{\xi\in \ah H}$, 
$Y=\ens{w_\xi}{\xi\in \bh H}$.

The sets $X,Y\sq\wo\cap\rL$ are $\ip{\hc}{\nn-1}$ in 
$\rL[G\res\Da_H]$ 
together with $\ah H$ and $\bh H$, and hence $\ip1{\nn}$, 
and $X\cap Y=\pu$. 
Suppose to the contrary that, in $\rL[G\res\Da_H]$,  
$X',Y'\sq\dn$ are 
disjoint sets in $\fs1{\nn}$, hence  
$\hs{\nn-1}$, 
such that $X\sq X'$ and $Y\sq Y'$. 
Then, in $\rL[G\res\Da_H]$,  
$$
A=\ens{\xi<\omil}{w_\xi\in X'}
\quad\text{and}\quad
B=\ens{\xi<\omil}{w_\xi\in Y'}
$$
are disjoint sets in $\hs{\nn-1}$, 
and we have $\ah H\sq A$ and $\bh H\sq B$ by construction, 
contrary to Theorem~\ref{Tsep'}.
The contradiction ends the proof of \ref{Tsep1}.
The proof of \ref{Tsep}\ref{Tsep2} is pretty similar.\vom

\epF{Theorem \ref{Tsep} modulo Theorem \ref{shon}}

\vyk{
\parf{The failure of $\fs1{\nn}$ separation}
\las{fail}

%
A version of the model of Section~\ref{ns}, 
where the $\fs1{\nn}$ separation fails for very 
similar reasons, can be can easily manufactured. 
Namely, coming back to Definition~\ref{ggh}, 
we make use of the sets 
$ 
\Da'_H=\ens{2\nu}{\nu\in \ah H}\cup
\ens{2\nu+1}{\nu\in \bh H}\, 
$ 
instead of $\Da_H$.
Then, similarly to Theorem~\ref{Tsep'}, it is true in 
the model  
$\morp H$ 
that $\ah H$ and $\bh H$ are disjoint $\is{\hc}{\nn-1}$ sets
not separable by disjoint\/ $\hp{\nn-1}$ sets.\snos
{In fact it is possible to prove, following 
\cite{kl28}, that a certain 
counterexample to the $\fs1{\nn}$ Separation 
in $\rL$ survives 
in the extension say of the type considered in 
Section~\ref{ns}.}

Moreover, it is possible to maintain both constructions 
together in the same model, so that Separation 
fails in the model for both $\fs1{\nn}$ and $\fp1{\nn}$. 

\qeD{Theorem \ref{Tsep} modulo Theorem \ref{shon}}
}


\parf{Auxiliary forcing relation}
\las{auxA}

Here we begin a lengthy proof of Theorem \ref{shon}. 
It involves an auxiliary forcing relation, not 
explicitly connected with any particular forcing notion,
in particular, with the key forcing $\fP$. 

\bgg
\lam{bass4}
We'll assume that $\nn\ge4$, since if $\nn=3$ then 
Theorem \ref{shon} holds by the Shoenfield
absoluteness.
\egg

{\ubf We argue in $\rL$.}
%
Consider 2nd order arithmetic language, with
variables $k,l,m,n,\dots$ of type $0$ over $\omega$ and 
variables $a,b,x,y,\dots$ of type $1$ over $\dn$,  
whose atomic formulas are those of the form $x(k)=n$.
Let $\cL$ be the extension of this language,
\index{formula!Lformula@\dd\cL formula}%
\index{zzL@$\cL$, language}%
which allows to substitute
variables of type $0$ with natural numbers
and variables of type $1$ with 
{\ubf small real names} (Definition~\ref{rk'}) 
$\rc\in\rL$.

\vyk{
By \rit{\dd\cL formulas} we understand formulas of this
\index{formula!Lformula@\dd\cL formula}%
\index{zzLformula@\dd\cL formula}%
extended language.
}

We define natural classes $\ls1n$, $\lp1n$ ($n\ge1$)
\index{formula!LSP@$\ls1n$, $\lp1n$, $\lsp11$}%
\index{zzLSP@$\ls1n$, $\lp1n$, $\lsp11$}%
of \dd\cL formulas.
Let $\lsp11$ be the closure of $\ls11\cup\lp11$ under 
${\neg}\yi{\land}\yi{\lor}$ and quantifiers over $\om$.
If $\vpi$ is a formula in $\ls1n$ (resp., $\lp1n$), then 
let $\vpi^-$ be  
\index{formula!phi-@$\vpi^-$}%
\index{zzzphi-@$\vpi^-$}%
the result of canonical transformation of 
$\neg\:\vpi$ 
to the $\lp1n$ (resp., $\ls1n$) form.

\vyk{
In the remainder, we'll need a somewhat more flexible 
relations between names and \muf s than in 
Definition~\ref{rk'}

\bdf
\lam{rk}
Let $\jpi$ be a \muf\ and $\rc$ is a \qn. 
Define a \qn\ 
\index{zzcpi@$\re\rc\jpi$}%
$\re\rc\jpi=\ens{\ang{\zp,n,i}}
{\zp\in\mt\jpi\land \sus\zq\,
(\zp\leq\zq\land\ang{\zq,n,i}\in\rc)}$, 
so that 
\imar{kk ni {re rc jpi}}%
\index{zzKcpini@$\kk ni{\re\rc\jpi}$}%
$$
\kk ni{\re\rc\jpi}=
\ens{\zp}
{\zp\in\mt\jpi\land \sus\zq\,
(\zp\leq\zq\land \zq\in\kk ni\rc)}
\,,\qquad n<\om,\;i=\ans{0,1}\,.
$$
Say that $\rc$ is a \rit{\srn\jpi} iff $\re\rc\jpi$
\index{real subname!pirealsubname@\srn\jpi}%
is a \rn\jpi. 
In this case, all sets 
\index{zzKcpin@$\kk n{}{\re\rc\jpi}$}%
$\kk n{}{\re\rc\jpi}=\kk n{0}{\re\rc\jpi}\cup
\kk n{1}{\re\rc\jpi}$ 
are pre-dense in $\mt\jpi$ by Definition \ref{rk'},
and then, obviously, open dense,

If $G\sq\mt\jpi$ is  
\dd{\mt\jpi}generic 
over the family of all  sets 
$\kk n{}{\re\rc\jpi}$, 
then say that $G$ is \rit{\dd\rc generic}.  
\index{set!cgeneric@\dd\rc generic}%
\index{generic@\dd\rc generic}%
In this case we can define the evaluation
\index{real name!evaluation, $\rc[G]$}%
\index{real subname!evaluation, $\rc[G]$}%
$\rc[G]=(\re\rc\jpi)[G]\in\dn$
($(\re\rc\jpi)[G]$ is defined by Definition~\ref{rk'}).

Accordingly if $\vpi$ is a \dd\cL formula and 
$G\sq\mt\jpi$ 
is \dd\rc generic for any name $\rc$ occurring 
in $\vpi$ then 
say that that $G$ is \rit{\dd\vpi generic} 
\index{set!figeneric@\dd\vpi generic}%
\index{generic@\dd\vpi generic}
and let $\vpi[G]$
be the result of substitution of $\rc[G]$ for any
name $\rc$ in $\vpi$.
Thus $\vpi[G]$ is a 2nd order arithmetic
formula, which may include numbers and
elements of $\dn$ as parameters.

If $\vjpi\in\vmf$ and $\jpi=\bkw\vjpi$ then 
\index{zzcpi@$\re\rc\vjpi$}%
$\re\rc\vjpi$, $\kk n{}{\re\rc\vjpi}$, 
\srn\vjpi, \etc\ will mean 
$\re\rc\jpi$, $\kk n{}{\re\rc\jpi}$, \srn\jpi, \etc 
\edf

\bcor
\lam{223}
If\/ $\vjpi,\vjqo\in\vmf$, $\mm$ is any set, 
$\vjpi\su_\mm\vjqo$, $\rc\in\mm$ is a\/ \srn\vjpi, 
then\/ $\rc$ is a\/ \srn\vjqo\ as well. 
\ecor
\bpf
Apply Lemma \ref{142} for $\re\rc\vjpi$. 
\epf
}

\vyk{
\bdf
\lam{zfm} 
Let $\zflm$ be the theory containing all axioms of 
$\zfc^-$ (minus the Power Set axiom) plus the axiom 
of constructibility $\rV=\rL$. 
\index{model!$\mm(x)$}%
\index{zzMx@$\mm(x)$}%
\index{model!CTM, countable transitive model}%
\index{CTM, countable transitive model}%
If $x\in\hc$ 
($\hc$= all hereditarily countable sets) then let
$\fl(x)$ be the least \rit{countable} transitive model 
CTM of $\zflm$ containing $x$. 
\edf 
}  

Now we define 
a relation $\zp\foe\vjpi\vpi$ between \mut s $\zp$, 
\muq s $\vjpi\in\vmf$,
and closed \dd\cL formulas $\vpi$ in 
$\lsp11$ or $\ls1n\cup\lp1n$, $n\ge2$, 
which will suitably approximate the true \dd\fP forcing 
relation.  
The definition goes on by 
induction on the complexity of $\vpi$.

\ben
\cenu
\itlb{fo2}
Let $\vjpi\in\vmf$, $\zp\in\md$    
(not necessarily $\zp\in\mt\vjpi$), and 
$\vpi$ is a closed $\lsp11$ formula. 
We define 
$\zp\foe\vjpi\vpi$ iff   
\index{forcing!forc@$\fof$}%
\index{zzforc@$\fof$}%
there is a CTM $\mm\mo\zflm$ 
(recall Definition \ref{zflm} on $\zflm$),
an ordinal $\vt<\dom\vjpi$, 
and a \mut\ $\zpo\in\mt{\vjpi\res\vt}$,   
such that

(1)
$\zp\leq \zpo$ 
(meaning: $\zpo$ is weaker),

(2)
$\mm$ contains $\vjpi\res\vt$ 
(then contains $\mt{\vjpi\res\vt}$ and $\zpo$ as well),

(3)
every name $\rc$ in $\vpi$ belongs to
$\mm$ and is \pol{\vjpi\res\vt},

(4)
$\vjpi\res\vt\su_\mm\vjpi$ --- therefore
$\vjpi\res\vt\su_{\ans\rc}\vjpi$
for any name $\rc$ in $\vpi$, \ and 

(5)
$\zpo$ \dd{\mt{\vjpi\res\vt}}forces $\vpi[\uG]$ 
over $\mm$ in the usual sense.%
\vyk{
\snos
{\label{snosfo2}%
As \ref{fo1} is assumed, every  
parameter $\rc$ in $\vjpi$ is a \srn\vjpi, 
and hence all sets 
$\kk n{}{\re\rc\vjpi}$ are open dense in $\mt\vjpi$. 
All these sets  
$\kk n{}{\re\rc\vjpi}$ belong to $\mm$ 
because $\rc,\vjpi\in\mm$. 
Therefore every set $G\sq\mt\vjpi$, 
generic over $\mm$, is \dd\vpi generic as well, 
so that the valuation $\vpi[G]$ is well-defined, 
and the \dd{\mt{\vjpi\res\vt}}forcing of $\vpi[\uG]$ 
in \ref{fo2} is well-defined as well.}
}%
\snos 
{Item \ref{fo2} not only requires $\vpi[\uG]$ 
to be forced, but also  
suitably seals this status by 
$\vjpi\res\vt\su_\mm\vjpi$. 
This will help us to prove the consistency of 
$\fof$ in Lemma~\ref{wf1}.}

\itlb{fo3}
If $\vpi(x)$ is a $\lp1n$ formula, $n\ge1$, then 
we define $\zp\foe\vjpi\sus x\,\vpi(x)$ iff there is a 
small \qn\ $\rc$ 
such that $\zp\foe\vjpi\vpi(\rc)$.

\itlb{fo4}
If $\vpi$ is a closed $\lp1n$ formula, $n\ge2$, 
then we define $\zp\foe\vjpi\vpi$ iff  
there is no \muq\   
$\vjta\in\vmf$ and \mut\ $\zp'\in\mt{\vjta}$ 
such that $\vjpi\sq\vjta$,
$\zp'\leq\zp$, 
and $\zp'\foe{\vjta}\vpi^-$.
\een

\bre
\lam{22c}
The condition
``$\zpo$ \dd{\mt{\vjpi\res\vt}}forces $\vpi[\uG]$ 
over $\mm$'' 
in \ref{fo2} does not depend on the 
choice of a CTM $\mm$ containing $\vjpi\res\vt$ 
and $\vpi$, 
since if $\vpi$ is $\lsp11$ then all transitive 
models agree on the formula $\vpi[G]$ 
by the Mostowski absoluteness theorem  
\cite[Theorem 25.4]{jechmill}.
\ere

\ble
\lam{mont}
Assume that \muq s\/ $\vjpi\sq\vjqo$ belong to\/ $\vmf$, 
$\zq,\zp\in\md$, 
$\zq\leq\zp$, 
$\vpi$ is an\/ \dd\cL formula.
Then\/ $\zp\foe\vjpi\vpi$ implies\/ $\zq\foe\vjqo\vpi$.
\ele
\bpf
If $\vpi$ is a $\lsp11$ formula,  
$\zp\foe\vjpi\vpi$, and this is witnessed 
by $\mm$, $\vt$, $\zpo$ as in \ref{fo2},
then the exactly same $\mm$, $\vt$, $\zpo$ 
witness $\zq\foe\vjqo\vpi$. 
\vyk{
For instance, to check \ref{fo1} in the part saying   
that every parameter $\rc$ in $\vpi$ is a \srn\vjqo, 
note that 
$\vjpi\res\vt\su_\mm\vjqo$ and $\rc\in\mm$ 
hold by definition, 
hence it remains to refer to Corollary~\ref{223}.
}

The induction step $\sus$, as in \ref{fo3}, is 
pretty elementary.

Now the induction step $\kaz$, as in \ref{fo4}.
Let $\vpi$ be a closed \dd{\lp1n}formula, $n\ge2$, and
$\zp\foe\vjpi\vpi$. 
Assume that $\zq\foe\vjqo\vpi$ fails.
Then by \ref{fo4} there exist:
a \muq\ $\vjqo'\in\vmf$ and \mut\ $\zq'\in\mt{\vjqo'}$ 
such that $\vjqo\sq\vjqo'$, $\zq'\leq\zq$,
and $\zq'\foe{\vjqo'}\vpi^-$.
But then $\vjpi\sq\vjqo'$ and $\zq'\leq\zp$,
hence $\zp\foe\vjpi\vpi$ fails by \ref{fo4}.
\epf

\bdf
\lam{Fd}
If $K$ is one of the classes
$\lsp11$, $\ls1n$, $\lp1n$ ($n\ge2$), 
then let $\fos K$ consist
of all triples $\ang{\vjpi,\zp,\vpi}$ such that 
$\zp\foe\vjpi\vpi$.
\edf

Then $\fos K$ is a subset of $\hc$.

\ble
[definability, in $\rL$]
\lam{deff}
$\fos{\lsp11}\in\id\hc1$.
If\/ $n\ge2$ then\/ 
$\fos{\ls1n}$ belongs to\/ $\is\hc{n-1}$ 
and\/ $\fos{\lp1{n}}$ belongs to\/ $\ip\hc{n-1}$. 
\ele
\bpf
Relations like $\vjpi\in\vmf$,  
``being a formula in $\lsp11$, $\ls1n$, $\lp1n$'', 
$\zp\in \mt{\vjro}$, forcing over a CTM, etc.\  
are definable in $\hc$ by bounded formulas, hence 
$\id\hc1$. 
On the top of this, the model $\mm$ can be 
tied by both $\sus$ and $\kaz$ in \ref{fo2}, 
see Remark~\ref{22c}.
This wraps up the $\id\hc1$ estimation for $\lsp11$. 

The inductive step by \ref{fo3} is quite simple.  

Now the step by \ref{fo4}.
Assume that 
$n\ge2$, and it is already established that 
$\fos{\ls1n}\in\is\hc{n-1}$. 
Then $\ang{\vjpi,\zp,\vpi}\in\fos{\lp1n}$
iff $\vjpi\in\vmf$, $\zp\in\md$, $\vpi$ 
is a closed $\lp1n$ formula, 
and, by \ref{fo4}, there exist no triple 
$\ang{\vjta,\zp',\psi}\in\fos{\ls1n}$
such that $\vjta\in\vmf$, 
$\vjpi\sq\vjta$, $\zp'\in\mt{\vjta}$, $\zp'\leq\zp$, 
and $\psi$ is $\vpi^-$.
We easily get the required 
estimation $\ip\hc{n-1}$ of $\fos{\lp1n}$. 
\epf



\vyk{
The following two lemmas are mainly related to 
the relation $\fof$  with respect to 
formulas in the class $\lsp11$.
}

\ble
[in $\rL$]
\lam{lsp}
Let\/ $\vjpi\in\vmf$, 
$\zp\in\mt\vjpi$, 
$\vpi$ is a formula in\/ $\lsp11$.  
%
\ben
\renu
\itlb{lsp1}
If\/ $\vjpi\sq\vjqo\in\vmf\cup\vmi$,
$\gN\mo\zflm$ is a TM containing\/ $\vjqo$ and\/ $\vpi$,
and\/ $\zp\foe\vjpi\vpi$, then\/ $\zp$
\dd{\mt\vjqo}forces\/ $\vpi[\uG]$ 
over\/ $\gN$ in the usual sense$.$

\itlb{lsp2}
If a TM\/ $\gN\mo\zflm$ contains\/ $\vjpi$, 
each name\/ $\rc$ in\/ $\vpi$
belongs to\/ $\gN$ and is\/ \pol\vjpi,
and\/ $\zp$ \dd{\mt\vjqo}forces\/ $\vpi[\uG]$ 
over\/ $\gN$, then there exists\/
$\vjqo\in\vmf$ such that\/
$\vjpi\su_\gN\vjqo$ and\/ $\zp\foe\vjqo\vpi$.
\een
\ele
\bpf
\ref{lsp1}
By definition there is an ordinal $\vt<\dom\vjpi$, 
a \mut\ $\zpo\in\mt{\vjpi\res\vt}$, 
and a CTM $\mm\mo\zfm$ containing $\vjpi\res\vt$ and 
such that
$\zp\leq\zpo$, 
every name $\rc$ in $\vpi$ belongs to $\mm$
and is \pol{\vjpi\res\vt}, 
$\vjpi\res\vt\su_\mm\vjpi$, 
and $\zpo$ \dd{\mt{\vjpi\res\vt}}forces $\vpi[\uG]$ 
over $\mm$. 
We can \noo\ assume that $\mm\sq\gN$. 
(Otherwise  $\gN\sq\mm$, and we replace $\gN$ by $\mm$.) 

Now suppose that $G\sq\mt\vjqo$ is a set 
\dd{\mt\vjqo}generic 
over $\gN$ and $\zp\in G$ --- then $\zpo\in G$, too. 
We have to prove that $\vpi[G]$ is true in $\gN[G]$. 

We claim that the set
$G'=G\cap\mt{\vjpi\res\vt}$ is 
\dd{\mt{\vjpi\res\vt}}generic over $\mm$.
Indeed, let a set $\zD\in\mm$, $\zD\sq\mt{\vjpi\res\vt}$, 
be open dense in\/ $\mt{\vjpi\res\vt}$.
Then, as ${\vjpi\res\vt}\su_\mm\vjqo$, 
$\zD$ is pre-dense in $\mt\vjqo$ 
by \ref{xisc}\ref{xisc2b}, 
and hence $G\cap\zD\ne\pu$ by the choice of $G$. 
It follows that  $G'\cap\zD\ne\pu$.

Now if $\rc$ is a name in $\vpi$ then 
$\rc\in\mm$ and $\rc$ is \pol{\vjpi\res\vt}.
It follows by the above that $\rc[G']\in\dn$
is defined.
Therefore $\rc[G]=\rc[G']$, because $G'\sq G$.
Thus $\vpi[G]$ coincides with $\vpi[G']$.
Note also that $\zpo\in G'$.
We conclude that $\vpi[G']$ holds in $\mm[G']$ as 
$\zpo$ forces $\vpi[\uG]$ over $\mm$. 
The same formula $\vpi[G]$ is holds 
$\gN[G]$ by the Mostowski absoluteness.

\ref{lsp2}
Lemma \ref{xisc}\ref{xisc1} yields $\vjqo\in\vmf$
such that $\vjpi\su_\gN\vjqo$.
\epf  

\ble
[in $\rL$]
\lam{wf1}
Let\/ $\vjpi\in\vmf$, $\zp\in\mt\vjpi$, 
$\vpi$ be a formula in\/ $\lsp11$ or\/ $\ls1n$, 
$n\ge2$. 
Then\/ 
$\zp\foe\vjpi\vpi$ and\/ $\zp\foe\vjpi\vpi^-$ 
cannot hold together. 
\vyk{
 and\/ $(2)$ there is a \muq\ 
$\vjqo\in\vmf$ and\/ $\zq\in\mt\vjqo$ such that\/ 
$\vjpi\su\vjqo$ and\/ $\zq\foe\vjqo\vpi$ or\/ 
$\zq\foe\vjqo\vpi^-$.
}
\ele
\bpf 
Let $\vpi\in\lsp11$. 
If 
both $\zp\foe\vjpi\vpi$ and $\zp\foe\vjpi\vpi^-$
then, by Lemma \ref{lsp}, 
$\zp$ \dd{\mt\vjpi}forces both $\vpi[\uG]$ 
and $\vpi^-[\uG]$ over a large enough CTM $\mm$, 
a contradiction.
If $\vpi\in\ls1n$ then the result 
follows by \ref{fo4}.
\vyk{
(2) 
If $\vpi\in\ls1n$, $n\ge2$, then the result 
follows by \ref{fo4}, so assume 
that $\vpi\in\lsp11$.
Let $\mm\mo\zflm$ be a CTM containing $\vpi$ and $\vjpi$.
By Lemma~\ref{xisc}\ref{xisc1},   
there is a \muq\ 
$\vjqo\in\vmf$ with\/ $\vjpi\su_{\mm}\vjqo$. 
The forcing $\mt{\vjpi}\in\mm$ contains $\zp$. 
There is a stronger condition 
$\zq\in\mt{\vjpi}$, $\zq\leq\zp$, 
which \dd{\mt{\vjpi}}forces either $\vpi$ or $\vpi^-$. 
Accordingly either $\zp\foe\vjqo\vpi$ or 
$\zp\foe\vjqo\vpi^-$, by \ref{fo2}.
}%
\epf

\parf{Tail invariance}
\las{tail}

Invariance theorems are very typical for 
all kinds of forcing. 
We prove two major invariance theorems on the 
auxiliary forcing. 
The first one shows tail invariance, 
while the other one 
(Section \ref{aut}) 
explores the permutational invariance.    

If $\vjpi=\sis{\jpi_\al}{\al<\la}\in\vmf$ and   
$\ga<\la=\dom\vjpi$ then let the 
\dd\ga\rit{tail} 
\index{tail@\dd\ga tail}%
$\vjpi\reg\ga$ be the restriction 
$\vjpi\res{[\ga,\la)}$ to the ordinal semiinterval 
$[\ga,\la)=\ens{\al}{\ga\le\al<\la}$.
Then the set 
$\mt{\vjpi\reg\ga}=\bkw_{\ga\le\al<\la}\vjpi(\al)$ 
is open dense in $\mt\vjpi$ 
by Lemma~\ref{xisc}\ref{xisc2x}.
Therefore it can be expected that if 
$\vjqo$ is another \muq\ of the same length 
$\la=\dom\vjqo$, and $\vjqo\reg\ga=\vjpi\reg\ga$, 
then the relation 
$\foe\vjpi$ coincides with $\foe\vjqo$.
And indeed this turns out to be the case (almost).

\bte
\lam{tat}
Assume that\/ $\vjpi,\vjqo$ are \muq s in\/ $\vmf$, 
$\ga<\la=\dom\vjpi=\dom\vjqo$, 
$\vjqo\reg\ga=\vjpi\reg\ga$, 
$\zp\in\md$, $n\ge2$, and\/ $\vpi$ is a formula 
in\/ $\lp1n\cup\ls1{n+1}$.
Then\/ 
$\zp\foe\vjpi\vpi$ iff\/ $\zp\foe{\vjqo}\vpi$. 
\ete

\bpf
{\ubf Part 1:} the $\lp12$ case.  
Let $\psi(x)$ be a $\ls11$ formula. 
Suppose 
that 
$\zp\foe{\vjqo}\kaz x\,\psi(x)$ {\ubf fails},
so there is  
$\vjqo'\in\vmf$ and a \mut\ $\zq\in\mt{\vjqo'}$ 
such that $\vjqo\sq\vjqo'$, $\zq\leq\zp$, and 
$\zq\foe{\vjqo'}\sus x\,\psi^-(x)$.
We can assume that 
$\zq\in\mt{\vjqo'\reg\ga}$.
By definition there is a small \qn\ $\rc$ 
such that $\zq\foe{\vjqo'}\psi^-(\rc)$.
\vyk{
There is a \mut\ 
$\zr\in\mt{\vjqo'\reg\ga}=\mt{\vjpi'\reg\ga}$, 
$\zr\leq\zq$.
Then still $\zr\leq\zp$ and $\zr\foe{\vjqo'}\psi^-(\rc)$.
}%

Let $\la'=\dom\vjqo'$. 
Define a \muq\ $\vjpi'$ so that 
$\dom{\vjpi'}=\la'=\dom{\vjqo'}$, $\vjpi\sq\vjpi'$, 
and $\vjpi'\reg\la=\vjqo'\reg\la$. 
Then $\vjpi'\reg\ga=\vjqo'\reg\ga$, 
hence $\zq\in\mt{\vjpi'\reg\ga}\sq\mt{\vjpi'}$. 

Consider any CTM $\gN\mo\zflm$ containing $\psi$,  
$\rc$, $\vjpi'$, $\vjqo'$. 
Then $\zq$ \dd{\mt{\vjqo'}}forces $\psi^-(\rc)[\uG]$ 
over $\gN$ by Lemma~\ref{lsp}. 
However the forcing 
notions $\mt{\vjpi'}$, $\mt{\vjqo'}$ contain one and 
the same dense set 
$\mt{\vjpi'\reg\ga}=\mt{\vjqo'\reg\ga}$.
Therefore $\zq$ also \dd{\mt{\vjpi'}}forces 
$\psi^-(\rc)[\uG]$ over $\gN$. 
Then by definition $\zq\foe{\vjpi'}\psi^-(\rc)$ and 
$\zq\foe{\vjpi'}\sus x\,\psi^-(x)$, 
hence 
$\zp\foe\vjpi\kaz x\,\psi(x)$ 
{\ubf fails}, as required.\vom 

{\ubf Part 2:} 
the step $\lp1n\to\ls1{n+1}$, $n\ge2$.  
Let $\vpi(x)$ be a formula in $\lp1n$. 
Assume that $\zp\foe\vjpi\sus x\,\vpi(x)$. 
By definition (see \ref{fo3} in Section~\ref{auxA}), 
there is a small \qn\ $\rc$ such that 
$\zp\foe\vjpi\vpi(\rc)$. 
Then we have $\zp\foe{\vjqo}\vpi(\rc)$ 
by the inductive hypothesis, thus  
$\zp\foe{\vjqo}\sus x\,\psi(x)$.\vom  

{\ubf Part 3:} 
the step $\ls1n\to\lp1{n}$, $n\ge3$. 
Let $\vpi$ be a $\lp1{n}$ formula, 
and $\zp\foe{\vjqo}\vpi$ {\ubf fails}.
Then by \ref{fo4} of Section~\ref{auxA}, 
there is a \muq\ 
$\vjqo'\in\vmf$ and a \mut\ $\zp'\in\mt{\vjqo'}$ 
such that $\vjqo\sq\vjqo'$, $\zp'\leq\zp$, and 
$\zp'\foe{\vjqo'}\vpi^-$.
By Lemma~\ref{xisc}\ref{xisc2x},
there is a \mut\ $\zr\in\mt{\vjqo'\reg\ga}$, $\zr\leq\zp'$.
Then $\zr\leq\zp$ and $\zr\foe{\vjqo'}\vpi^-$.
Define a \muq\ $\vjpi'\in\vmf$ by 
$\dom{\vjpi'}=\la'=\dom{\vjqo'}$, $\vjpi\sq\vjpi'$, 
and $\vjpi'\reg\la=\vjqo'\reg\la$. 
Then 
$\zr\in\mt{\vjpi'\reg\ga}$, 
$\zr\leq\zp$, and also $\zr\foe{\vjpi'}\vpi^-$   
by the inductive hypothesis.  
We conclude that  
$\zp\foe\vjpi\vpi$
{\ubf fails} as well. 
\epf

\parf{Permutations}
\las{aut}

Still {\ubf arguing in $\rL$}, 
we let $\aut$ be the set of all 
bijections $\hh:\omi\onto\omi$, such that 
$\hh=\hh\obr$ and  
the \rit{non-identity domain} 
\index{non-identity domain@$\abk\hh$}%
\index{zzNID@$\abk\hh$}%
$\abk\hh=\ens{\xi}{\hh(\xi)\ne\xi}$ 
is at most countable. 
Elements of $\aut$ will be called \rit{permutations}.
\index{permutation!$\aut$}%
\index{zzperm@$\aut$}%
\vyk{
If $m<\om$ then let $\aut_m$ consist of those 
permutations $\hh\in\aut$ satisfying 
$\abk\hh\sq\omi\ti(\om\bez m)$, 
in other words,  
$\hh(\xi,k)=\ang{\xi,k}$ 
for all $\xi<\omi$, $k<m$.
}

Let $\hh\in\aut$. 
\index{permutation!action}%
We extend the action of $\hh$  as follows.
\bit
\item 
if $\zp$ is a \mut\ then $\hh\zp$ is 
a \mut, 
$\abs{\hh\zp}=\ima\hh\zp=
\ens{\hh(\xi)}{\xi\in\abs\zp}$, 
and $(\hh\zp)(\hh(\xi))=\zp(\xi)$
whenever $\xi\in\abs\zp$, 
in other words, $\hh\zp$ coincides with the 
superposition $\zp\circ{(\hh\obr)}$;

\item 
if $\zpi\in\md$ is a \muf\ then 
$\hh\ap\zpi=\zpi\circ{(\hh\obr)}$ is 
a \muf, 
$\abs{\hh\ap\zpi}=\ima\hh\zpi$ 
and $(\hh\ap\zpi)(\hh(\xi))=\zpi(\xi)$
whenever $\xi\in\abs\zpi$;  


\item
if $\rc\sq\md\ti(\om\ti\om)$ is a \qn, 
then put 
$\hh\rc=\ens{\ang{\hh\zp, n,i}}{\ang{\zp, n,i}\in\rc}$, 
thus easily $\hh\rc$ is a \qn\ as well;

\item
if $\vjpi=\sis{\jpi_\al}{\al<\ka}\in\vmf$, then   
$\hh\vjpi=\sis{\hh\ap\jpi_\al}{\al<\ka}$, 
still a \muq\ in $\vmf$; 

\item
if $\vpi:=\vpi(\rc_1,\dots,\rc_n)$ is a \dd\cL formula 
(with all names explicitly indicated), then   
$\hh\vpi$ is $\vpi(\hh\rc_1,\dots,\hh\rc_n)$.
\eit

Many notions and relations defined 
above are clearly \dd\aut invariant, \eg, 
$\zp\in\mt\jpi$ iff $\hh\zp\in\mt{\hh\ap\jpi}$, 
$\jpi\ssb{}\jqo$ iff ${\hh\ap\jpi}\ssb{}{\hh\ap\jqo}$, 
{\sl et cetera}. 
The invariance also takes place with respect to 
the relation $\fof$.

\bte
\lam{auT}
Assume that\/ $\vjpi\in\vmf$, $\zp\in\mt\vjpi$, 
$\hh\in\aut$, 
$n\ge2$, and\/ 
$\vpi$ belongs to\/ $\lp1n\cup\ls1{n+1}$.  
Then\/ $\zp\foe\vjpi\vpi$ iff\/ 
$(\hh\zp)\foe{\hh\vjpi}(\hh\vpi)$. 
\ete

\bpf
Let $\vjqo=\hh\vjpi$, $\zq=\hh\zp$. 

{\ubf Part 1:} the $\lp12$ case. 
Assume that $\vpi(x)$ is a $\ls11$ formula, 
$\psi(x):=\hh\vpi(x)$, 
and $\zq\foe{\vjqo}\kaz x\,\psi(x)$ {\ubf fails}.
Then by definition
(\ref{fo2}, \ref{fo3} of Section~\ref{auxA})
there is a \muq\ $\vjqo'\in\vmf$,
a \mut\ $\zq'\in\mt{\vjqo'}$, and
a small \qn\ $\rd$, 
such that
$\vjqo\su\vjqo'$,
$\zq'\leq\zq$, and 
$\zq'\foe{\vjqo'}\sus x\,\psi^-(\rd)$.
The \muq\ $\vjpi'=\hh\obr\vjqo'$ then satisfies 
$\vjpi\su\vjqo$, the \mut\ $\zp'=\hh\obr\zq'$ 
belongs to $\mt{\vjpi'}$, $\zp'\leq\zp$, 
 and $\rc=\hh\obr{\rd}$ 
is a small \qn. 
However we cannot now claim that 
$\zp'\foe{\vjpi'}\vpi^-(\rc)$, since the existence 
of $\gM$, $\vt$ as in \ref{fo2} in Section~\ref{auxA} 
is not necessarily preserved by the action of $\hh\obr$ 
or $\hh$.

To circumwent this difficulty, let $\gM\mo\zflm$ be 
a CTM containing $\vjpi'\yi\vjqo'\yi\hh\yi\rc\yi\rd$ and 
(all names in) $\vpi,\psi$.
Then $\zq'$ \dd{\mt{\vjqo'}}forces 
$\psi^-(\rd)[\uG]$ over $\gM$
by Lemma \ref{lsp}\ref{lsp1}.
Then $\zp'$ \dd{\mt{\vjpi'}}forces 
$\vpi^-(\rc)[\uG]$ over $\gM$, by the 
standard theorems of forcing.
Lemma\ref{lsp}\ref{lsp2} yields a 
\muq\ $\vjta\in\vmf$ with $\vjpi'\su \vjta$,
such that $\zp'\foe{\vjta}\vpi^-(\rc)$, hence 
$\zp'\foe{\vjta}\sus x\,\vpi^-(x)$
by \ref{fo3}.
However $\vjpi\su\vjpi'\su\vjta$ and $\zp'\leq\zp$,
therefore, $\zp\foe{\vjpi}\kaz x\,\vpi(x)$
{\ubf fails} by \ref{fo4}, as required.

{\ubf Part 2:} 
the step $\lp1n\to\ls1{n+1}$, $n\ge2$. 
Let $\vpi(x)$ be a formula in $\lp1n$ and  
$\psi(x):=\hh\vpi(x)$. 
Assume that $\zp\foe\vjpi\sus x\,\vpi(x)$. 
By definition (\ref{fo3} in Section~\ref{auxA}), 
there is a small \qn\ $\rc$ such that 
$\zp\foe\vjpi\vpi(\rc)$. 
Then we have $\zq\foe{\vjqo}\psi(\rd)$ 
by inductive assumption, where 
$\rd=\hh\rc$ is a small \qn\ itself. 
Thus  
$\zq\foe{\vjqo}\sus x\,\psi(x)$. 

{\ubf Part 3:} 
the step $\ls1n\to\lp1{n}$, $n\ge3$. 
Let 
$\vpi$ be a formula in $\lp1n$, and
$\zq\foe{\vjqo}\psi$ {\ubf fails}, 
where $\zq=\hh\zp$, $\vjqo=\hh\vjpi$, 
and $\psi$ is $\hh\vpi$, as above. 
By \ref{fo4},   
there is a \muq\ 
$\vjqo'\in\vmf$ and a \mut\  $\zq'\in\mt{\vjqo'}$ 
such that $\vjqo\sq\vjqo'$,
$\zq'\leq\zq$, and $\zq'\foe{\vjqo'}\psi^-$.
%
Now let $\zp'=\hh\obr\zq'$ and $\vjpi'=\hh\obr\vjqo'$, 
so that $\zp'\leq\zp$ and $\vjpi\sq\vjpi'$. 
We have $\zp'\foe{\vjpi'}\vpi^-$ 
by inductive assumption. 
We conclude that $\zp\foe\vjpi\vpi$ {\ubf fails}, 
as required.  
\epf

\parf{Forcing inside the key \muq}
\las{fkm}

The following Theorem~\ref{fofot} will show that 
the forcing relation $\foe\vjpi$, considered 
with countable initial segments $\vjpi =\vjPi\res\al$ 
of the key sequence $\vjPi$, coincides 
with the true $\fP$-forcing relation up to
level $\nn-1$.\vom 

{\ubf We argue in $\rL$.} 
Recall that the key \muq\ 
$\vjPi=\sis{\jpn\al}{\al<\omi}\in\vmi$, satisfying 
\ref{kyt1}, \ref{kyt0}, \ref{kyt2}, \ref{kyt3} of 
Theorem~\ref{kyt}, was introduced by 
\ref{keys}, and 
$\fP=\mt\vjPi$ is our forcing notion. 
In addition, $\nn\ge4$ by \ref{bass4}.

\bdf
\lam{foe}
We write 
$\zp\foe\al\vpi$ instead of $\zp\foe{\vjPi\res\al}\vpi$, 
for the sake of brevity. 
Let $\zp\fof\vpi$ mean: $\zp\foe\al\vpi$ for some $\al<\omi$.
\edf

\ble
[in $\rL$]
\lam{pi}
Assume that\/ $\zp\in\fP$, $\al<\omi$, and\/ 
$\zp\foe\al\vpi$. 
Then$:$
\ben
\renu
\itlb{pi1}
if\/ $\al\le\ba<\omi$, $\zq\in\fpl\ba=\mt{\vjPi\res\ba}$, 
and\/ $\zq\leq\zp$, 
then\/ $\zq\foe\ba\vpi\;;$ 

\itlb{pi2}
if\/ $\zq\in\fP$, $\zq\leq\zp$, then\/ $\zq\foe\ba\vpi$ 
for some\/ $\ba\,;$ $\al\le\ba<\omi\;;$ 

\itlb{pi3}
if\/ $\zq\in\fP$  
and\/ $\zq\fof\vpi^-$ then\/ $\zp,\zq$ 
are\/ \sad$;$

\itlb{pi4}
therefore, 1st, if\/ $\zp,\zq\in\fP$, $\zq\leq\zp$, 
and\/ $\zp\fof\vpi$ then\/ $\zq\fof\vpi$, and 2nd,\/  
$\zp\fof\vpi$, $\zp\fof\vpi^-$ cannot hold together.
\een
\ele
\bpf
To prove \ref{pi1} apply Lemma~\ref{mont}.
To prove \ref{pi2} pick $\ba$ such that
$\al<\ba<\omi$ 
and $\zq\in\mt{\vjPi\res\ba}$, and apply \ref{pi1}.
To prove \ref{pi3} note that $\zp,\zq$ 
are incompatible in $\fP$, as otherwise 
\ref{pi1} leads to contradiction, but the 
incompatibility in $\fP$ implies being \sad\ 
by Corollary~\ref{lsad}.
\epf

\bte
\lam{fofot}
If\/ 
$\vpi$ is a closed\/ $\cL$-formula in\/ 
$\lsp11\cup\ls12\cup\lp12\cup\ldots\cup
\ls1{\nn-2}\cup\lp1{\nn-2}\cup\ls1{\nn-1}$
and\/ $\zp\in\fP$, 
then\/ $\zp$ \dd\fP forces\/ $\vpi[\uG]$ over\/ $\rL$ 
in the usual sense, if and only if\/ $\zp\fof\vpi$.
\ete
\bpf
Let $\ofo$ denote the usual \dd\fP forcing relation 
over $\rL$.\vom 

{\ubf Part 1\/}: $\vpi$ is a formula in $\lsp11$. 
If $\zp\fof\vpi$ then $\zp\foe{\vjPi\res\ga}\vpi$
for some $\ga<\omi$, and then $\zp\ofo\vpi[\uG]$ 
by Lemma~\ref{lsp} with $\vjqo=\vjPi$ and $\gN=\rL$. 

Suppose now that $\zp\ofo\vpi[\uG]$. 
There is an ordinal $\gao<\omi$ such that 
$\zp\in\fP_{\gao}=\mt{\vjPi\res\gao}$ and $\vpi$ 
belongs to $\fl(\vjPi\res\gao)$. 
(Recall Definition \ref{zflm} on models $\fl(x)\mo\zflm$.)
The set $U$ of all sequences $\vjpi\in\vmf$ such that 
$\gao<\dom\vjpi$ and there is an ordinal $\vt$, 
$\gao<\vt<\dom\vjpi$, such that 
$\vjpi\res\vt\su_{\fl(\vjpi\res\vt)}\vjpi$, is 
dense in $\vmf$ by Lemma~\ref{xisc}\ref{xisc1}, 
and is $\hd 1$. 
Therefore by Corollary~\ref{18a} there is an 
ordinal $\ga<\omi$ such that  
$\vjpi=\vjPi\res\ga\in U$. 
Let this be witnessed by an ordinal $\vt$, 
$\gao<\vt<\ga=\dom\vjpi$ and  
$\vjpi\res\vt\su_{\fl(\vjpi\res\vt)}\vjpi$. 
We claim that $\zp$ \dd{\mt{\vjpi\res\vt}}forces 
$\vpi[\uG]$ over $\fl(\vjpi\res\vt)$ in the 
usual sense --- then 
by definition $\zp\foe\vjpi\vpi$, and we are done. 

To prove the claim, assume otherwise. 
Then there is a \mut\ $\zq\in\mt{\vjPi\res\vt}$, 
$\zq\leq\zp$, which \dd{\mt{\vjpi\res\vt}}forces 
$\neg\:\vpi[\uG]$ over $\fl(\vjpi\res\vt)$. 
Then by definition
(\ref{fo2} in Section \ref{auxA})
$\zq\foe\vjpi\neg\,\vpi$ holds, 
hence $\zq\fof\neg\,\vpi$, and then 
$\zq\ofo\neg\:\vpi[\uG]$ (see above), 
with a contradiction to $\zp\ofo\vpi[\uG]$.\vom 

{\ubf Part 2\/}:  
the step $\lp1n\to\ls1{n+1}$ ($n\ge1$). 
Consider a $\lp1n$ formula $\vpi(x)$. 
Assume $\zp\fof\sus x\,\vpi(x)$. 
By definition there is a small \qn\ $\rc$ such that 
$\zp\fof\vpi(\rc)$. 
By inductive hypothesis, $\zp\ofo\vpi(c)[\uG]$, that is, 
$\zp\ofo\sus x\,\vpi(x)[\uG]$. 
Conversely, assume that $\zp\ofo\sus x\,\vpi(x)[\uG]$. 
As $\fP$ is CCC, there is a small \qn\ $\rc$ (in $\rL$) 
such that $\zp\ofo\vpi(\rc)[\uG]$. 
We have $\zp\fof\vpi(\rc)$ by the inductive hypothesis, 
hence $\zp\fof\sus x\,\vpi(x)$.\vom

{\ubf Part 3\/}:  
the step $\ls1n\to\lp1{n}$ ($2\le n\le\nn-2$). 
Assume that $\vpi$ is a closed $\ls1n$ formula,
and $\zp\fof\vpi^-$.
By Lemma~\ref{pi}\ref{pi4}, there is no \mut\ 
$\zq\in\fP$, $\zq\leq\zp$, with $\zq\fof\vpi$.
This implies $\zp\ofo\vpi^-$ by 
the inductive hypothesis. 

Conversely, let $\zp\ofo\vpi^-$.
There is an ordinal $\gao<\omi$ such that 
$\zp\in\fP_{\gao}=\mt{\vjPi\res\gao}$ and $\vpi$ 
belongs to $\fl(\vjPi\res\gao)$. 
Consider the set $U$ of all 
\muq s 
$\vjpi\in\vmf$ such that 
$\dom\vjpi>\gao$ and 
there is a \mut\ $\zq\in\mt\vjpi$ satisfying 
$\zq\leq\zp$
and $\zq\foe\vjpi\vpi$.
Then 
$U$ belongs to $\hs{n-1}$ 
($\vpi$, $\zpo$ as parameters) 
by Lemma~\ref{deff}, 
hence to $\hs{\nn-3}$, where $\nn\ge4$ by \ref{bass4}. 
Therefore by \ref{keys} 
(and \ref{kyt2} of Theorem~\ref{kyt})
there is $\ga<\omi$ such that 
$\vjPi\res\ga$ \se s $U$. 

{\ubf Case 1\/}: 
$\vjPi\res\ga\in U$. 
Let this be witnessed by a 
\mut\ $\zq\in\mt\vjpi$, so that in particular 
$\zq\leq\zp$ and $\ga>\gao$. 
Thus $\zq\in\mt{\vjPi\res\ga}$, $\zq\leq\zp$, and 
$\zq\foe{\vjPi\res\ga}\vpi$, that is, 
$\zq\ofo\vpi[\uG]$ by the inductive hypothesis, 
contrary to the choice of $\zp$.  
Therefore Case 1 cannot happen, and we have:

{\ubf Case 2\/}: 
no \muq\ in $U$ extends $\vjPi\res\ga$.
We can assume that $\ga>\gao$. 
(If not, replace $\ga$ by $\gao+1$.) 
We claim that $\zp\foe{\ga}\vpi^-$. 
Indeed otherwise by \ref{fo4} 
there is a \muq\ $\vjpi\in\vmf$   
and a \mut\ $\zq\in\mt{\vjpi}$, such that
$\vjPi\res\ga\sq\vjpi$,
$\zq\leq\zp$, and $\zq\foe{\vjqo}\vpi$.
But then $\vjpi$ 
belongs to $U$.  
On the other hand, 
$\vjPi\res\ga\sq\vjpi$, 
contrary to the Case 2 assumption.
Thus indeed $\zp\fof\vpi^-$, as required.  
\epf

\parf{Elementary equivalence theorem}
\las{ee}

\bpf[Theorem~\ref{shon}]
Suppose the contrary. 
Then there is a $\ip1{\nn-2}$ formula $\vpi(r,x)$ with 
$r\in\dn\cap\rL[G\res X]$ as the only parameter, 
and a real $x_0\in\dn\cap\rL[G]$
such that $\vpi(r,x_0)$ is true in $\rL[G]$ 
but there is no $x\in\dn\cap\rL[G\res X]$ 
such that $\vpi(r,x)$ is true in $\rL[G]$. 
By a version of Lemma~\ref{r2n}, 
we have $r=\rco[G]$, where
$\rco\sq\mt{\jPi\res X}\ti\om\ti2$
is a small 
\rn{(\jPi\res X)}. 
\imar{utochnit'}%
(See Section~\ref{kge} on notation.)
And there is a small \rn{\fP}
$\rc\sq\fP\ti\om\ti2$ such that 
$x_0=\rc[\uG]$. 

By Theorem~\ref{fofot}, there is a \mut\ $\zpo\in G$ 
such that 
\ben
\nenu
\itlb{po1}
$\zpo$ \ \dd\fP forces \ 
`$
\vpi(\rc_0[\uG],\rc[\uG])
\,\land\,
\neg\:\sus x\in \rL[\uG\res X]\,\vpi(\rc_0[\uG],x) 
$'
\
over $\rL$; 

\itlb{eet*}\msur
$\zpo\fof \vpi(\rc_0,\rc)$, that is, 
$\zpo\foe{\vjPi\res\gao}\vpi(\rc_0,\rc)$, 
where $\gao<\omi$ --- 
and we can assume that $\zpo\in\mt{\vjPi\res\gao}$ 
as well.   
\een
As $\rc,\rco$ are small names, there is an 
ordinal $\da<\omi$ satisfying 
\ben
\atc
\atc
\nenu
\itlb{eet*2}\msur
$\abc\rco\sq \da\cap X$, $\abc\rc\sq\da$, 
 and $\abc\zpo\sq\da$. 
\een
As $\abc\vjPi=\omi$ 
by Corollary~\ref{doml}, 
we can enlarge $\gao$, if necessary, 
to make sure that 
\ben
\atc
\atc
\atc
\nenu
\itlb{eet*3}\msur
$\da\sq\abc{\vjPi\res\gao}$, 
that is, if $\et<\da$ then 
$\et\in\abc{\jPi_{\al'}}$ for some 
$\al'=\al'(\et)<\gao$.
\een
We start from here towards a contradiction. 
Let $D=\da\bez X$.  

Let $U$ consist of all \muq s $\vjpi\in\vmf$, 
such that 
$\vjPi\res\gao\su\vjpi$, 
and hence $\zpo\in\mt\vjpi$ by \ref{eet*}, 
and 
there is 
$\za<\dom\vjpi$ and 
$\hh\in \aut$ such that 
\ben
\Aenu
\itlb{bab2} 
$\abk\hh\cap(\da\cap X)=\pu$, and 
$\hh$ maps $D$ onto a set $R\sq X\bez\da$;  


\itlb{bab1}\msur
$\gao\le\za<\dom\vjpi$ and  
$(\hh\vjpi)\reg\za=\vjpi\reg\za$, 
that is, $\hh(\vjpi(\al))=\vjpi(\al)$ whenever 
$\za\le\al<\dom\vjpi$.
\een
It holds by routine estimations that $U$ is 
a $\hs{1}$ set 
(with $\vjPi\res\gao$, $\da$ as parameters), 
hence a $\hs{\nn-3}$ set because $\nn\ge4$ by 
\ref{bass4}. 
Therefore by \ref{keys} 
there is an ordinal $\ga<\omi$ such that 
$\vjPi\res\ga$ \ \se s $U$.\vom

{\ubf Case 1\/}: $\vjPi\res\ga\in U$, 
so that  
\ref{bab2}, 
\ref{bab1} 
hold for $\vjpi=\vjPi\res\ga$, via some  
$\za\in[\gao,\ga)$ and $\hh\in\aut$.
In particular, by \ref{bab1}, $\hh(\jPi_\al)=\jPi_\al$ 
whenever $\za\le\al<\ga$.
By Lemma~\ref{mont} and \ref{eet*}, we have 
$\zpo\foe{\vjPi\res\ga}\vpi(\rco,\rc)$. 
Let $\rc'=\hh\rc$, $\zpo'=\hh\zpo$. 
Note that $\hh\rco=\rco$ since $\abc\rco\cap\abk\hh=\pu$ 
by \ref{bab2}. 
Now Theorem~\ref{auT} implies  
$\zpo'\foe{\hh\ap(\vjPi\res\ga)}\vpi(\rco,\rc')$. 
%
Thus  
$\zpo'\foe{\vjPi\res\ga}\vpi(\rco,\rc')$ 
holds by Theorem~\ref{tat} and \ref{bab1}. 
But the common domain $\abc\zpo\cap\abc{\zpo'}$ 
does not intersect  
$\abk\hh$ by \ref{bab2} since $\abc{\zpo}\sq\da$. 
It follows that $\zpo,\zpo'$ are compatible, 
basically $\zp=\zpo\cup\zpo'\in\md$ 
(not necessarily $\in\mt{\vjPi\res\ga}$)  
and $\zp\leq\zpo'$, hence still 
$\zp\foe{\vjPi\res\ga}\vpi(\rco,\rc')$.

Unfortunately Theorem~\ref{fofot} is not applicable 
immediately to conclude that 
$\zp$ \ \dd\fP forces  
$
\vpi(\rco[\uG],\rc'[\uG])
$
over $\rL$, simply because $\zp$ may not belong to 
$\fP$. 
We need an additional argument. 
Recall that  $\zpo\in\mt{\vjPi\res\gao}$, hence 
$\zpo'\in\mt{\hh\ap(\vjPi\res\gao)}$. 
As $\za>\gao$, there is a \mut\ 
$\zqo\in\mt{\hh\ap\jPi_\za}$ 
satisfying $\abs\zqo=\abs{\zpo'}$ and $\zqo\leq\zpo'$. 
Then still 
$\zqo\foe{\vjPi\res\ga}\vpi(\rco,\rc')$ 
(because $\zpo'\foe{\vjPi\res\ga}\vpi(\rco,\rc')$),
and $\zqo\in\mt{\jPi_\za}$ since 
$\hh\ap\jPi_\za=\jPi_\za$. 
Thus $\zqo\in\mt{\vjPi\res\ga}$. 
Moreover $\zqo$ is compatible with $\zpo$ in 
$\mt{\vjPi\res\ga}$ because 
$\abs\zqo=\abs{\zpo'}$ and $\zqo\leq\zpo'$, 
and $\zpo'$ coinsides with $\zpo$ on the common 
domain $\abs\zpo\cap\abs{\zpo'}=\da\cap X$. 
Thus there exists $\zq\in\mt{\vjPi\res\ga}$ with 
$\zq\leq\zpo$, $\zq\leq\zqo$. 
Then $\zq\foe{\vjPi\res\ga}\vpi(\rco,\rc')$ 
holds, and we conclude that
\ben
\nenu
\atc\atc\atc
\atc
\itlb{999}
$\zq$ \ \dd\fP forces  
$
\vpi(\rco[\uG],\rc'[\uG])
$
over $\rL$ 
\een 
by Theorem~\ref{fofot}.  
However $\abc{\rc'}\sq (\da\cap X)\cup R\sq X$ 
by construction, 
and hence $\rc'[\uG]\in \rL[\uG\res X]$ is forced. 
Thus $\zq$ \ \dd\fP forces  
$\sus x\in \rL[\uG\res X]\,\vpi(\rco[\uG],x)$
over $\rL$ by \ref{999}, contrary to \ref{po1}. 
The contradiction closes Case 1.\vom

{\ubf Case 2\/}:
no \muq\ in 
$U$ extends $\vjPi\res\ga$. 
We can assume that $\ga>\gao$. 
(Otherwise replace $\ga$ by $\gao+1$.) 
Pick any set $R\sq X\bez\da$ 
satisfying 
$\card{R}=\card D$. 
However $D\sq\da$,  
hence $D\cap R=\pu$, so there is a 
permutation $\hh\in\aut$,
$\hh:D\onto R$, satisfying 
$\abk\hh=D\cup R$, hence \ref{bab2}.

Pick any ordinal $\la$, $\ga<\la<\omi$. 
Our plan is now to somewhat modify $\vjPi\res\la$ 
in order to fulfill \ref{bab1} as well, 
with $\za=\ga$.
The modification will replace 
the \dd Rpart of $\vjPi\res\la$ above $\ga$ 
by the \dd\hh copy of its \dd Dpart.
To render this in detail, recall that 
$\vjPi\res\la=\sis{\jPi_\al}{\al<\la}$, where 
each $\jPi_\al$ is a small special \muf, whose domain 
$d_\al=\abc{\jPi_\al}\sq\omi$ is countable. 
If $\al<\ga$ then  put $\jpi_\al=\jPi_\al$. 
Suppose that $\ga\le\al<\la$. 
Then $D\sq\abc{\jPi_\al}$ by \ref{eet*3}.
On the base of $\jPi_\al$, 
define a modified \muf\ $\jpi_\al$ such that
\ben
\aenu
\itlb{jpia1}\msur
$\abc{\jpi_\al}=d_\al\cup R$ --- note that 
$D\sq d_\al=\abc{\jPi_\al}\sq\abc{\jpi_\al}$ 
in this case because 
$D\sq\da\sq\abc{\vjPi\res\ga}$ by \ref{eet*3} 
(as $\gao\le\ga$); 

\itlb{jpia2}
if $\xi\in d_\al\bez R$ then  
$\jpi_\al(\xi)=\jPi_\al(\xi)$, 

\itlb{jpia3}
if $\xi\in D$, so  
$\hh(\xi)=\et\in R$, then  
$\jpi_\al(\et)=\jPi_\al(\xi
)$. 
\een
We assert that $\vjpi=\sis{\jpi_\al}{\al<\la}\in\vmf$, 
that is, if $\al<\ba<\la$ then 
$\jpi_\al\bssq\jpi_\ba$. 
This amounts to the following: if 
$\et\in \abc{\jpi_\al}$ then 
$\jpi_\al(\et)\ssq\jpi_\ba(\et)$. 

If $\et\nin R$ then $\jpi_\al(\et)=\jPi_\al(\et)$ 
by construction. 
It remains to check that 
$\jpi_\al(\et)\ssq\jpi_\ba(\et)$ whenever 
$\al<\ba<\la$, 
$\et=\hh(\xi)\in R\cap \abc{\jpi_\al}$, and 
$\xi\in D$.
If now $\al<\ga$ then $R\cap\abc{\jpi_\al}=\pu$ 
by the choice of $R$, so it remains to consider  
the case when $\ga\le\al$. 
Then $\xi,\et\in\abc{\jpi_\al}$ by construction, 
and we have 
$\jpi_\al(\et)=\jPi_\al(\xi)$ and 
$\jpi_\ba(\et)=\jPi_\ba(\xi)$.  
Therefore $\jpi_\al(\xi)\ssq\jpi_\ba(\xi)$, 
and we are done.

We claim that the \muq\ 
$\vjpi=\sis{\jpi_\al}{\al<\la}$ 
satisfies $\vjPi\res\ga\sq\vjpi$ and 
\ref{bab2}, \ref{bab1}. 
Indeed $\vjPi\res\ga\sq\vjpi$ as $\ga\ge\gao$. 
\ref{bab2}  
hold by construction. 
We claim that \ref{bab1} is satisfied with  
$\za=\ga$, that is, if $\ga\le\al<\la$ then 
$\hh\ap\jpi_\al=\jpi_\al$. 
Indeed $D\cup R\sq\abc{\jpi_\al}$ by 
\ref{jpia1}, and hence  
$\hh\ap\jpi_\al=\jpi_\al$ holds 
by \ref{jpia2}, \ref{jpia3}.

Thus $\vjpi\in U$  
and $\vjPi\res\ga\su\vjpi$. 
But this contradicts to the Case 2 assumption.\vom

To conclude, either case leads to a contradiction. 
\epF{Theorem~\ref{shon}} 

\qeD{Theorem \ref{Tsep}, see the end of Section~\ref{nspf}}

\parf{Remarks and problems}
\las{rp}

One may ask what happens with the separation theorem 
at other projective levels $m\ne\nn$ in the model 
of Section~\ref{ns}. 
As for the above levels, it happens that, in the model 
$\rL[G\res\Da_H]$ of Theorem~\ref{Tsep'}, 
there is a ``good'' $\id1{\nn+1}$ wellordering of the 
reals, of length $\omi$. 
(The gaps in $\Da_H$ do not allow the wellorder 
construction of Corollary~\ref{mod3} does not go 
through at level $\nn$!) 
It follows by standard arguments 
that the separation theorem holds for $\fp1m$ 
and fails for $\fs1m$, for all $m>\nn$, in the model 
$\rL[G\res\Da_H]$. 
As for the levels $3\le m<\nn$, we conjecture that 
separation 
holds for $\fp1m$ and fails for $\fs1m$
in $\rL[G\res\Da_H]$, 
but this problem is open.

Let $\fP_\nn$ be the forcing notion $\fP$ defined 
in Section~\ref{bpf} for a given $\nn\ge3$. 
Using a certain amalgamation of all 
$\fP_\nn$, $\nn\ge3$, defined by a rather 
sophisticated product-like construction, first applied 
in \cite[A]{h74} and \cite{k79ian}, a generic 
extension of $\rL$ can be defined, in which the 
separation theorem {\ubf fails} for all classes
$\fs1\nn$, $\fp1\nn$,  
$\nn\ge3$.

And finally, 
it is an interesting and perhaps very difficult problem 
to define a generic extension of $\rL$ in which the 
separation theorem {\ubf holds} for a given class
$\fs1\nn$,  $\nn\ge3$, beginning with say $\fs13$. 
This problem has been open since early years of 
forcing, see \cite[Problem 3029]{matsur}.

\vyk{
We recall that $\fp13$ Separation \rit{holds} in $\rL$. 
Thus Theorem~\ref{Tsep} in fact shows that the $\fp13$ 
Separation principle is destroyed in an appropriate generic 
extension of $\rL$.  
It would be interesting to find a generic extension in which, 
the other way around, the $\fs13$ Separation 
(false in $\rL$) holds.
This can be a difficult problem.
At least, the model used to prove Theorem~\ref{Tsep} does
not help: it was establ \cite{kl28} that any
pair of disjoint $\fs13$ sets, non-separable by
disjoint $\fp13$ sets in $\rL$, remains $\fs13$ and
non-separable by disjoint $\fp13$ sets  in the extension.

A complicated alternative proof of Theorem~\ref{Tsep} 
can be obtained 
with the help of \rit{countable-support} products and 
iterations of Jensen's forcing studied earlier in 
\cite{abr,k79,k83}. 
The \rit{finite-support} approach which we pursue here yields 
a significantly more compact proof, which still uses some 
basic constructions from \cite{h74}. 
Whether countable-support 
products and iterations can lead to  
countable-section non-uniformization results 
remains to be seen.
}


\vyk{
\back
The authors thank Ali Enayat and Mohammad Golshani for 
fruitful discussions, and an anonymous reviewer for valuable
remarks and comments.
\eack
}

\renek{\refname}
{{\large References}\addcontentsline{toc}{subsection}{References}}

\bibliographystyle{plain}
{\small
\bibliography{49,kle}
}

\small

\renek{\indexname}
{%
{\large Index}
\addcontentsline{toc}{subsection}{Index}
}

\printindex

\end{document}